\documentclass[11pt,oneside,reqno]{amsart}
\usepackage{amssymb}
\usepackage{color}
\usepackage{tikz-cd}
\usetikzlibrary{decorations.markings}

\usepackage{xcolor}
\definecolor{deepgreen}{cmyk}{0.99998,0,1,0}

\usepackage{relsize}

\usepackage{hyperref}
\hypersetup{hypertex=true,
	colorlinks=true,
	linkcolor=blue,
	anchorcolor=blue,
	citecolor=blue}
\usepackage{comment}
\usepackage{mathrsfs}

\allowdisplaybreaks[4]

\usepackage{tabularx}

\usepackage{enumitem}

\usepackage{etoolbox}

  \usepackage{bm}
  \usepackage{bbm}
  
  \usepackage[new]{old-arrows}
  
\theoremstyle{definition}
  \newtheorem{defi}{$\mathbf{Definition}$}[section]
  \newtheorem*{pro}{$\mathbf{Proof}$}
  
  \theoremstyle{plain}
  \newtheorem{theo}[defi]{$\mathbf{Theorem}$}
  \newtheorem{lemma}[defi]{$\mathbf{Lemma}$}

  \newtheorem{prop}[defi]{$\mathbf{Proposition}$}

\newcommand{\tro}{\mathrm{Tr}}

\newcommand{\pa}{\partial}

\newcommand{\lk}{\left(}
\newcommand{\rk}{\right)}

\newcommand{\lv}{\left\vert}
\newcommand{\rv}{\right\vert}
\newcommand{\lV}{\left\Vert}
\newcommand{\rV}{\right\Vert}

\newcommand{\bv}{\big\vert}
\newcommand{\Bv}{\Big\vert}
\newcommand{\bbv}{\bigg\vert}

\newcommand{\bV}{\big\Vert}
\newcommand{\BV}{\Big\Vert}
\newcommand{\bbV}{\bigg\Vert}

\makeatletter
\newcommand{\itemtaglabel}[2]{\protected@edef\@currentlabel{[#1]}%
	\label{#2}%
}
\makeatother

\definecolor{pink}{RGB}{249,164,186}
\definecolor{grassgreen}{RGB}{128,255,0}

\numberwithin{equation}{section}

\title{Discrete Mixed Quantization}

\author{Qiaochu Ma}
\address{Department of Mathematics, Texas A\&M University}
\email{qiaochu@tamu.edu}
\date{}

\makeatletter
\newcommand*{\rom}[1]{\expandafter\@slowromancap\romannumeral #1@}
\makeatother

\begin{document}
\clearpage\maketitle

\begin{abstract}
In this paper, we develop a mixed quantization technique for graph vector bundles and apply it to several asymptotic spectral problems, including the Alon-Boppana bound, the Kesten-McKay law, asymptotic determinant, quantum ergodicity, zero divisor convergence, and Ramanujan vector bundles.
\end{abstract}


\section{Introduction}

\subsection{Background}\label{A1''}

In \cite{MR2838248, MR3615411}, Bismut-Ma-Zhang obtained an asymptotic expansion for analytic torsion. This was later extended by Ma \cite{MR4665497} to full asymptotic analytic torsion, and by Puchol \cite{MR4611826} to asymptotic holomorphic torsion. Inspired by the framework and techniques of \cite{MR2838248,MR3615411,MR4665497,MR4611826}, Ma-Ma \cite{MR4808253} and Ben Ovadia-Ma-Rodriguez-Hertz \cite{ovadia2025mixedquantizationpartialhyperbolicity} introduced a \emph{uniform} version of quantum ergodicity (QE) for unitary flat vector bundles over manifolds and studied its stability under small perturbations. In these works, mixed quantization plays a central role, bringing together semiclassical and geometric quantization.

Using the pseudodifferential calculus on homogeneous trees developed by Le Masson \cite{MR3245884}, Anantharaman-Le Masson \cite{MR3322309} established the first discrete version of QE for a family of regular graphs. Brooks-Le Masson-Lindenstrauss \cite{MR3567266} proved QE for averaging operators on the two-sphere.

In this paper, we develop the technique of mixed quantization for vector bundles over regular graphs and use it to study several asymptotic spectral problems, including spectral gap, eigenvalue density, QE, asymptotic determinant, zero divisor convergence, and Ramanujan vector bundles, all exhibiting \emph{a strong uniformity property}. Informally, in our setup, there are two parameters, the scale parameter and the spin parameter. Our main results require only that \emph{one of these two parameters tends to infinity}, with the other parameter allowed to vary arbitrarily. In comparison, the uniform QE of \cite{MR4808253} asserts that eigensections become equidistributed as the frequency tends to infinity, uniformly with respect to the spin parameter.

We now explain in more detail.

\subsection{Geometric setup}

We now introduce the geometric setting that will be used throughout the paper.

\subsubsection{Graph vector bundles}

Let $X$ be a finite regular graph. A unitary vector bundle 
\begin{equation}
F=(X\times \mathbb{C}^\ell,\phi)
\end{equation}
over $X$ with $\ell$-dimensional fibres consists of \emph{the total space} $X\times \mathbb{C}^\ell$ and \emph{a connection} $\phi$. We denote the fibre over $x\in X$ by $F_x$ or $\mathbb{C}^\ell_x$. The connection $\phi$ assigns to each oriented edge $(x,x')$ of $X$ (with head $x$ and tail $x'$) a \emph{unitary map}
\begin{equation}\label{1.2}
\phi_{x,x'}\colon\mathbb{C}^\ell_{x'}\to \mathbb{C}^\ell_{x}
\end{equation}
such that $\phi_{x,x'}=\phi_{x',x}^{-1}$.

Unlike vector bundles over manifolds, we \emph{do not} impose any algebraic cocycle condition such as $\phi_{x'',x'}\phi_{x',x}=\phi_{x'',x}$. We shall also use dimension notation slightly different from the standard convention in differential geometry. Namely, we write $\dim_{\mathbb{C}} F_x$ for the complex dimension of the fibre over $x\in X$, and $\dim_{\mathbb{C}} F$ for the dimension of the total space,
\begin{equation}
	\dim_{\mathbb{C}}F_x=\ell,\quad\dim_{\mathbb{C}}F=\sum_{x\in X}\dim_{\mathbb{C}} F_{x}=\ell\lv X\rv,
\end{equation}
where $\lv X\rv$ denotes the cardinality of the vertex set of $X$.

A section $u$ of $F$ assigns to each $x\in X$ a vector $u(x)\in \mathbb{C}^\ell_x$, and the space of all sections is denoted by $C^\infty(X,\mathbb{C}^\ell)$. Here $C^\infty$ is only a conventional notation, and $X$ has no smooth structure. We define the Laplacian $\Delta^F$ by
\begin{equation}\label{1.3}
	\begin{split}
	&\Delta^F\colon C^\infty(X,\mathbb{C}^\ell)\to C^\infty(X,\mathbb{C}^\ell),\\
	&(\Delta^F u)(x)=\sum_{x'\sim x}\phi_{x,x'}u(x'),
	\end{split}
\end{equation}
where $x\sim x'$ means that $x$ and $x'$ are adjacent.

A loop (or a cycle) in $X$ is a sequence $(x_k,\cdots,x_0)$ such that $x_i\sim x_{i+1}$ and $x_0=x_k$. The fundamental group $\pi_1(X)$ of $X$ is generated by loops based at a fixed vertex, modulo homotopy equivalence. For any loop $\gamma=(x_k,\cdots,x_0)\in\pi_1(X)$, we denote
\begin{equation}\label{1.5n}
\phi_\gamma=\phi_{x_k,x_{k-1}}\cdots \phi_{x_1,x_{0}}\colon \mathbb{C}_{x_0}^\ell\to \mathbb{C}_{x_0}^\ell
\end{equation}
the corresponding holonomy map.

\subsubsection{Complex geometry}\label{s1.2.2}

Let $N$ be a compact complex manifold with $\dim_{\mathbb{C}}N=n$, and $(L,h^L)$ a \emph{positive} Hermitian holomorphic line bundle over $N$. Denote by $c_1(L,h^L)$ the first Chern form of $(L,h^L)$, and by $dv_N=c_1(L,h^L)^n/n!$ the associated volume form.

For ${p\in{\mathbb{N}}}$, let $L^p=L^{\otimes p}$ be the $p$-th tensor power of $L$, and let $H^{(0,0)}(N,L^p)$ denote the (finite dimensional) space of holomorphic sections of $L^p$.

\subsubsection{A double series of vector bundles}

Consider a series of $d$-regular graphs
\begin{equation}\label{1.4}
(X_m)_{m\in\mathbb{N}}.
\end{equation}
For each edge $(x,x')$ of $X_m$, we assign a holomorphic diffeomorphism $\phi_{m,x,x'}$ of $N$, and suppose that it lifts to a holomorphic isometry of $L$. Then, by \eqref{1.5n}, this induces a $\pi_1(X_m)$-action on $(N,L)$. Moreover, for each $p\in\mathbb{N}$, the map $\phi_{m,x,x'}$ acts unitarily on $H^{(0,0)}(N,L^p)$. We denote these actions by
\begin{equation}\label{1.5}
(\phi_{m,x,x'}, \pi_1(X_m))\curvearrowright (N,L),\quad\phi_{m,p,x,x'}\curvearrowright H^{(0,0)}(N,L^p).
\end{equation}
By \eqref{1.2} and \eqref{1.3}, we obtain a unitary vector bundle with its Laplacian
 \begin{equation}\label{1.8}
 	\begin{split}
 		&F_{m,p}=\big(X_m\times H^{(0,0)}(N,L^p),\phi_{m,p,x,x'}\big),\\
 		&\Delta^{F_{m,p}}\colon C^\infty\big(X_m,H^{(0,0)}(N,L^p)\big)\to C^\infty\big(X_m,H^{(0,0)}(N,L^p)\big),\\
 		&(\Delta^{F_{m,p}} u)(x)=\sum_{x'\sim x}\phi_{m,p,x,x'}u(x').
 	\end{split}
 \end{equation}

We list the eigenvalues of $\Delta^{F_{m,p}}$, counted with multiplicity, and choose an associated orthonormal basis of eigensections
\begin{equation}\label{a1}
	\Delta^{F_{m,p}}u_{m,p,i}=\lambda_{m,p,i}u_{m,p,i},\quad\lV u_{m,p,i}\rV^2_{L^2(X_m,F_{m,p})}=1.
\end{equation}

\subsubsection{Some assumptions}

We say that the series $(X_m)_{m\in\mathbb{N}}$ given in \eqref{1.4} converges to the $d$-regular tree $\mathbb{T}_d$ in
the sense of Benjamini-Schramm, denoted by \textbf{[BST]}, if for any $k\in\mathbb{N}$,
\begin{equation}\label{bst}
	\lim_{m\to\infty}\frac{\lv\{x\in X_m\mid \mathrm{inj}_x\leqslant k\}\rv}{\lv X_m\rv}=0,
\end{equation}
where $\mathrm{inj}_x$ denotes the injectivity radius of $X_m$ at $x$.

We say that the action \eqref{1.5} is free, denoted by \textbf{[FREE]}, if for any $m\in\mathbb{N}$, the only element $\gamma_m\in\pi_1(X_m)$ acting trivially on $N$ is $\gamma_m=1$. Equivalently, the induced morphism
\begin{equation}\label{free}
\pi_1(X_m)\to \text{the diffeomorphism group of }N
\end{equation}
is injective. 

We say that the action \eqref{1.5} is dense, denoted by \textbf{[DEN]}, if for any $m\in\mathbb{N}$, there exists $z\in N$ such that the orbit
\begin{equation}\label{den}
	\{\gamma z\mid \gamma\in\pi_1(X_m)\}\subset N
\end{equation}
 is dense.

 We say that the action \eqref{1.5} has a spectral gap, denoted by \textbf{[GAP]}, if there is $\varepsilon>0$ such that for any $m,p\in\mathbb{N}$,
 \begin{equation}\label{sg}
 \mathrm{spec}(\Delta^{F_{m,p}})\subseteq [-d,d-\varepsilon]\cup \{d\}.
 \end{equation}

We say that the action \eqref{1.5} induces an expander, denoted by \textbf{[EXP]}, if there is $\varepsilon>0$ such that for any $m\in\mathbb{N}$, the averaging operator $\Delta^{X_m\times N}$, defined by
\begin{equation}\label{1.13n}
	\begin{split}
&\Delta^{X_m\times N}\colon L^2(X_m\times N)\to L^2(X_m\times N),\\
&\Delta^{X_m\times N}f(x,z)=\sum_{x'\sim x}f(x',\phi_{m,x,x'}^{-1}z),
	\end{split}
\end{equation}
satisfies
\begin{equation}\label{exp}
	\begin{split}
		&\mathrm{spec}(\Delta^{X_m\times N})\subseteq [-d,d-\varepsilon]\cup \{d\},\\
		&\ker\big(\Delta^{X_m\times N}-d\mathrm{Id}\big)=\mathrm{span}\{\mathbbm{1}_{X_m\times N}\}.
	\end{split}
\end{equation}
In other words, $d$ is a simple discrete eigenvalue of $\Delta^{X_m\times N}$, with a uniform spectral gap for all $m\in\mathbb{N}$. In general, \textbf{[EXP]} \eqref{exp} is much stronger than \textbf{[DEN]} \eqref{den}.

\subsection{Main results}

We shall prove several strong uniform limits involving two parameters $(m,p)$. For a double series $(a_{m,p})_{m,p\in\mathbb{N}}$ of complex numbers, consider the condition
\begin{equation}
\lim_{k\to\infty}\big(\sup_{(m,p)\in\mathbb{N}^{2},\mathrm{max}(m,p)\geqslant k}\lv a_{m,p}\rv_\mathbb{C}\big)=0,
\end{equation}
where $\lv\cdot\rv_{\mathbb{C}}$ denotes the modulus on $\mathbb{C}$. We denote this by $\limsup_{\max(m,p)\to\infty}\lv a_{m,p}\rv_\mathbb{C}=0$ for short.

This condition is strictly stronger than requiring
\begin{equation}
	\begin{split}
\lim_{m\to\infty}\lv a_{m,p}\rv_\mathbb{C}=0,&\text{ for any fixed } p,\\
		\lim_{p\to\infty}\lv a_{m,p}\rv_\mathbb{C}=0,&\text{ for any fixed } m,
	\end{split}
\end{equation}
which we refer to as the \emph{large scale limit} and \emph{large spin limit}, respectively. Indeed, it means that for any sequence $(m_j,p_j)_{j\in\mathbb{N}}$, as long as we have $\lim_{j\to\infty}\max(m_j,p_j)=\infty$, then $\lim_{j\to\infty}\lv a_{m_j,p_j}\rv_\mathbb{C}=0$. For example, we may take $p_j=j$, while $m_j$ oscillates wildly. It may rise to $10$, then return to $0$, then rise to $100$, return to $0$, then rise to $1000$, return to $0$, and so on.

\subsubsection{Alon-Boppana bound}

Define the \emph{nontrivial spectral radius} of $F_{m,p}$ by
\begin{equation}
	r_{m,p,\mathrm{nt}}=\max\big\{\lv\lambda_{m,p,i}\rv\mid\lambda_{m,p,i}\neq\pm d\big\}.
\end{equation}
Note that the eigenvalue multiset may or may not contain $\pm d$, and when it does, these values may occur with multiplicity at most $\dim_{\mathbb{C}}H^{(0,0)}(N,L^p)$.

We now state our first result, an Alon-Boppana type bound.
\begin{theo}\label{t1.1}
Assume \textbf{[BST]} \eqref{bst}, \textbf{[FREE]} \eqref{free}, and \textbf{[DEN]} \eqref{den}. Then
		\begin{equation}\label{1.19n}
			\liminf_{\max(m,p)\to\infty}r_{m,p,\mathrm{nt}}\geqslant 2\sqrt{d-1}.
		\end{equation}
\end{theo}

We refer to \S\,\ref{s5n} for the full statement and proof. In particular, in the large scale limit $m\to\infty$ with $p$ fixed, it suffices to assume \textbf{[BST]} \eqref{bst}, and in the large spin limit $p\to\infty$ with $m$ fixed, it suffices to assume \textbf{[FREE]} \eqref{free} and \textbf{[DEN]} \eqref{den}.

\subsubsection{Kesten-McKay law}

We define the Kesten-McKay distribution $d\mu_{\mathrm{KM}}(\lambda)$ by
\begin{equation}\label{1.21n}
	d\mu_{\mathrm{KM}}(\lambda)=\frac{d\sqrt{4(d-1)-\lambda^2}}{2\pi\big(d^2-\lambda^2\big)}\mathbbm{1}_{[-2\sqrt{d-1},2\sqrt{d-1}]}(\lambda)d\lambda.
\end{equation}
Here, the $d$'s in $d\mu_{\mathrm{KM}}(\lambda)$ and $d\lambda$ denote measure notation, and should not be confused with the degree $d$ of the graph.

We now state a Kesten-McKay type eigenvalue density law. 

\begin{theo}\label{t1.2}
Assume \textbf{[BST]} \eqref{bst} and \textbf{[FREE]} \eqref{free}. Then for any interval $I\subseteq \mathbb{R}$,
	\begin{equation}\label{1.20n}
		\limsup_{\max(m,p)\to\infty}\lv\frac{1}{\dim_\mathbb{C}F_{m,p}}\lv\{i\mid\lambda_{m,p,i}\in I\}\rv-\int_{I}d\mu_{\mathrm{KM}}(\lambda)\rv=0.
	\end{equation}
\end{theo}

We refer to \S\,\ref{s6n} for the full statement and proof. In particular, in the large scale limit $m\to\infty$ with $p$ fixed, it suffices to assume \textbf{[BST]} \eqref{bst}, and in the large spin limit $p\to\infty$ with $m$ fixed, it suffices to assume \textbf{[FREE]} \eqref{free}.

\subsubsection{Asymptotic determinant}

The operator $(d-\Delta^{F_{m,p}})$ is nonnegative. We define the normalized log determinant of $(d\mathrm{Id}_F-\Delta^F)$ by
\begin{equation}
	\begin{split}
		\frac{1}{\dim_{\mathbb{C}}F_{m,p}}\ln\det(d-\Delta^{F_{m,p}})&=\sum_{\lambda_{m,p,i}\neq d}\frac{1}{\dim_{\mathbb{C}}F_{m,p}}\ln(d-\lambda_{m,p,i}).
	\end{split}
\end{equation}

We now state the asymptotic log determinants result.
\begin{theo}
Assume \textbf{[BST]} \eqref{bst}, \textbf{[FREE]} \eqref{free}, and \textbf{[GAP]} \eqref{sg}. Then
	\begin{equation}\label{n7.8n}
		\begin{split}
			&\limsup_{\max(m,p)\to\infty}\\
			&\Bv\frac{1}{\dim_\mathbb{C}F_{m,p}}\ln\det(d-\Delta^{F_{m,p}})\\
			&-\Big((d-1)\ln(d-1)-\frac{1}{2}(d-2)\ln(d-2)-\frac{1}{2}(d-2)\ln(d)\Big)\Bv=0.
		\end{split}
	\end{equation}
\end{theo}

We refer to \S\,\ref{det} for the full statement and proof. In particular, in the large scale limit $m\to\infty$ with $p$ fixed, it suffices to assume \textbf{[BST]} \eqref{bst} and \textbf{[GAP]} \eqref{sg}, and in the large spin limit $p\to\infty$ with $m$ fixed, it suffices to assume \textbf{[FREE]} \eqref{free} and \textbf{[GAP]} \eqref{sg}.

\subsubsection{Strong uniform QE}

For any $x\in X_m$ and any eigensection $u_{m,p,i}$ from \eqref{a1}, we have $u_{m,p,i}(x)\in H^{(0,0)}(N,L^p)$. We therefore write $u_{m,p,i}(x,z)=\big(u_{m,p,i}(x)\big)(z)$. This defines the following probability measure on $X_m\times N$,
\begin{equation}\label{1.10}
	\begin{split}
	&\mathscr{Q}_m\in C^\infty(X_m\times N)\\
		&\mapsto\sum_{x\in X_m}\int_{N}\mathscr{Q}_m(x,z)\lV u_{m,p,i}(x,z)\rV_{h^{L_z^p}}^2d{v}_{N}(z)\in\mathbb{C}.
	\end{split}
\end{equation}
We now state QE, asserting that the quantum states in \eqref{1.10} become equidistributed.

\begin{theo}\label{C9'}
Assume \textbf{[BST]} \eqref{bst}, \textbf{[FREE]} \eqref{free}, and \textbf{[EXP]} \eqref{exp}. Then for any series of smooth functions $(\mathscr{Q}_m\in C^\infty(X_m\times N))_{m\in\mathbb{N}}$ satisfying $\lV\mathscr{Q}_m\rV_{L^\infty(X_m\times N)}\leqslant C$, we have
	\begin{equation}\label{9.}
		\begin{split}
	\limsup_{\max(m,p)\to\infty}&\frac{1}{\dim_\mathbb{C}F_{m,p}}\sum_{i=1}^{\dim_\mathbb{C}F_{m,p}}\\
	\bbv&\sum_{x\in X_m}\int_{N}\mathscr{Q}_m(x,z)\lV u_{m,p,i}(x,z)\rV_{h^{L_z^p}}^2d{v}_{N}(z)\\
	&-\frac{1}{\lv X_m\rv\mathrm{Vol}(N)}\sum_{x\in X_m}\int_{N}\mathscr{Q}_m(x,z)dv_{N}(z)\bbv=0.
		\end{split}
	\end{equation}
\end{theo}

We refer to \S\,\ref{s10} for the full statement and proof. In particular, in the large scale limit $m\to\infty$ with $p$ fixed, it suffices to assume \textbf{[BST]} \eqref{bst} and \textbf{[EXP]} \eqref{exp}, and in the large spin limit $p\to\infty$ with $m$ fixed, it suffices to assume \textbf{[FREE]} \eqref{free} and \textbf{[DEN]} \eqref{den}.

\subsubsection{Zero divisors convergence}

For any $x\in X_m$ and any eigensection $u_{m,p,i}$ in \eqref{a1}, let $\mathrm{Div}(u_{m,p,i}(x))$ denote the zero divisor of $u_{m,p,i}(x)\in H^{(0,0)}(N,L^p)$, which is a formal sum of analytic subvarieties of $N$ of complex codimension one.

The following zero divisor convergence result only involves the large spin limit $p\to\infty$, and at present, no large scale analogue is known.
\begin{theo}\label{t1.4}
Assume \textbf{[FREE]} \eqref{free} and \textbf{[DEN]} \eqref{den}. Fix $m\in\mathbb{N}$. Then there exists a series of subsets $\big(I_{m,p}\subseteq \{1,\cdots, \dim_{\mathbb{C}}F_{m,p}\}\big)_{p\in\mathbb{N}}$ of asymptotic density one, in the sense that
\begin{equation}
	\lim_{p\to\infty}\lv I_{m,p}\rv\big/\dim_{\mathbb{C}}F_{m,p}=1,
\end{equation}
such that for any series $(i_p\in I_{m,p})_{p\in\mathbb{N}}$, and for any continuous $(n-1,n-1)$-form $\mathscr{Q}_m^{(n-1,n-1)}$ on $X_m\times N$, we have
		\begin{equation}\label{1.25}
		\begin{split}
			\lim_{p\to\infty}\frac{1}{p}&\sum_{x\in X_m}\int_{\mathrm{Div}(u_{m,p,i_p}(x))}\mathscr{Q}_m^{(n-1,n-1)}(x,z)\\
			=&\sum_{x\in X_m}\int_{N}\mathscr{Q}_m^{(n-1,n-1)}(x,z)\wedge c_1(L,h^L).
		\end{split}
	\end{equation}
\end{theo}

We refer to \S\,\ref{s11} for the full statement and proof.

\subsection{Ramanujan vector bundles}

The assumptions \textbf{[BST]} \eqref{bst}, \textbf{[FREE]} \eqref{free}, and \textbf{[DEN]} \eqref{den} are relatively mild, whereas \textbf{[GAP]} \eqref{sg} and \textbf{[EXP]} \eqref{exp} are quite strong. Thus, the existence of even a single nontrivial example is not a priori clear. Perhaps surprisingly, an example compatible with all of our main results was already constructed nearly forty years ago.

Let $p_0,p_1$ be distinct primes satisfying $p_0,p_1\equiv1\mod{4}$. Define
\begin{equation}\label{1.25n}
	\begin{split}
A=\big\{a=(a_0,a_1,a_2,a_3)\in\mathbb{Z}^4\mid a_0^2+a_1^2+a_2^2+a_3^2=p_0,&\\
a_0 \text{ positive and odd}, a_i \text{ even for }i=1,2,3&\big\}.
	\end{split}
\end{equation}
Then $\lv A\rv=p_0+1$. Fix $b\in\mathbb{Z}$ with $b^2\equiv-1\mod{p_1}$. For any $a\in A$, define $a_{\mathbb{Z}/p_1\mathbb{Z}}\in\mathrm{PSL}_2(\mathbb{Z}/p_1\mathbb{Z})$ by
\begin{equation}\label{1.26}
	a_{\mathbb{Z}/p_1\mathbb{Z}}=\left(\begin{matrix}a_0+ba_1&a_2+ba_3\\-a_2+ba_3&a_0-ba_1\end{matrix}\right).
\end{equation}

If $\lk\frac{p_0}{p_1}\rk=-1$, we define $X_{(p_0,p_1)}$ to be the \emph{Cayley graph} of $\mathrm{PGL}_2(\mathbb{Z}/p_1\mathbb{Z})$ with respect to the generator set $\{a_{\mathbb{Z}/p_1\mathbb{Z}}\mid a\in A\}$. Its vertices are the elements of $\mathrm{PGL}_2(\mathbb{Z}/p_1\mathbb{Z})$, and two vertices $(x,x')$ are connected by an edge if $x=a_{\mathbb{Z}/p_1\mathbb{Z}}x'$
for some generator $a_{\mathbb{Z}/p_1\mathbb{Z}}$. If $\lk\frac{p_0}{p_1}\rk=1$, we define $X_{(p_0,p_1)}$ similarly to be the Cayley graph of $\mathrm{PSL}_2(\mathbb{Z}/p_1\mathbb{Z})$.

We now take the complex geometry data in \S\,\ref{s1.2.2} to be
\begin{equation}\label{1.26n}
		\big(N,L,c_1(L,h^L),H^{(0,0)}(N,L^p)\big)= \big(\mathbb{CP}^1,\mathcal{O}_{\mathbb{CP}^1}(1),\omega_{\mathrm{FS}},\mathrm{Sym}^p(\mathbb{C}^2)\big).
\end{equation}
Here $\mathbb{CP}^1$ is the complex projective line, $\mathcal{O}_{\mathbb{CP}^1}(1)$ is the dual tautological line bundle, $\omega_{\mathrm{FS}}$ is the Fubini-Study form, and $\mathrm{Sym}^p(\mathbb{C}^2)$ is the $p$-th symmetric power of $\mathbb{C}^2$, which can be viewed as the space of homogeneous polynomial functions of degree $p$,
\begin{equation}\label{1.28}
	\mathrm{Sym}^p(\mathbb{C}^2)=\{c_0z_0^p+c_1z_0^{p-1}z_1+\cdots+c_p{z_1}^p\mid c_0,\cdots,c_p\in\mathbb{C}\}.
\end{equation}
This is a geometric realization of irreducible $\mathrm{SU}_2(\mathbb{C})$-representations. For $g\in\mathrm{SU}_2(\mathbb{C})$ and $s\in\mathrm{Sym}^p(\mathbb{C}^2)$, we define $gs\in\mathrm{Sym}^p(\mathbb{C}^2)$ by
\begin{equation}\label{1.29n}
	(gs)\lk\begin{smallmatrix}
		z_0\\
		z_1
	\end{smallmatrix}\rk=s\big(g^{-1}\lk\begin{smallmatrix}
	z_0\\
	z_1
	\end{smallmatrix}\rk\big).
\end{equation}

We now construct a natural $\mathrm{SU}_2(\mathbb{C})$-connection over $X_{(p_0,p_1)}$, that is, we assign an element $\phi_{(p_0,p_1),x,x'}\in \mathrm{SU}_2(\mathbb{C})$ to each oriented edge $(x,x')$ of $X_{(p_0,p_1)}$. Suppose that $x=a_{\mathbb{Z}/p_1\mathbb{Z}}x'$ for some generator $a_{\mathbb{Z}/p_1\mathbb{Z}}$. Looking at \eqref{1.26}, then mathematics explains itself, and the associated element $a_{\mathbb{R}}\in \mathrm{SU}_2(\mathbb{C})$ must be
\begin{equation}\label{1.30nn}
a_{\mathbb{R}}=\frac{1}{\sqrt{p_0}}\left(\begin{matrix}a_0+\sqrt{-1}a_1&a_2+\sqrt{-1}a_3\\-a_2+\sqrt{-1}a_3&a_0-\sqrt{-1}a_1\end{matrix}\right),
\end{equation}
obtained by replacing the mod $p_1$ imaginary unit $b$ with the genuine imaginary unit $\sqrt{-1}$ and normalizing the determinant.  We then define $F_{(p_0,p_1),p}=(X_{(p_0,p_1)}\times\mathrm{Sym}^p(\mathbb{C}^2),a_{\mathbb{R}})$, the \emph{Cayley vector bundle} over $X_{(p_0,p_1)}$. That is, for any oriented edge $(x,x')$ with $x=a_{\mathbb{Z}/p_1\mathbb{Z}}x'$, we set $\phi_{(p_0,p_1),x,x'}=a_{\mathbb{R}}$. 

We denote the above construction by
\begin{equation}\label{1.30n}
	\begin{split}
		&\big(F_{(p_0,p_1),p},X_{(p_0,p_1)}\big)\\
		&=\begin{cases}
			\mathrm{Cay}\Big(\mathrm{PGL}_2(\mathbb{Z}/p_1\mathbb{Z}),\mathrm{Sym}^p(\mathbb{C}^2),(a_{\mathbb{Z}/p_1\mathbb{Z}},a_\mathbb{R})\Big),&\text{if }\lk\frac{p_0}{p_1}\rk=-1,\\
				\mathrm{Cay}\Big(\mathrm{PSL}_2(\mathbb{Z}/p_1\mathbb{Z}),\mathrm{Sym}^p(\mathbb{C}^2),(a_{\mathbb{Z}/p_1\mathbb{Z}},a_\mathbb{R})\Big),&\text{if }\lk\frac{p_0}{p_1}\rk=1.
		\end{cases}
	\end{split}
\end{equation}
Clearly, we have
\begin{equation}\label{1.32n}
	\begin{split}
&\lv X_{(p_0,p_1)}\rv=\begin{cases}
	p_1(p_1^2-1),&\text{if }\lk\frac{p_0}{p_1}\rk=-1,\\
		\frac{1}{2}p_1(p_1^2-1), &\text{if }\lk\frac{p_0}{p_1}\rk=1,
\end{cases}\\
&\dim_\mathbb{C}F_{(p_0,p_1),p}=\lv X_{(p_0,p_1)}\rv(p+1).
	\end{split}
\end{equation}

 By \eqref{1.8}, \eqref{1.25n}, \eqref{1.26}, \eqref{1.26n}, \eqref{1.29n}, \eqref{1.30nn}, and \eqref{1.30n}, the Laplacian $\Delta^{F_{(p_0,p_1),p}}$ is given by
\begin{equation}
	\begin{split}
		&\Delta^{F_{(p_0,p_1),p}}\colon C^\infty\big(X_{(p_0,p_1)},\mathrm{Sym}^p(\mathbb{C}^2)\big)\to C^\infty\big(X_{(p_0,p_1)},\mathrm{Sym}^p(\mathbb{C}^2)\big),\\
		&(\Delta^{F_{(p_0,p_1),p}} u)(x)=\sum_{a\in A}a_\mathbb{R}\cdot u\big(a_{\mathbb{Z}/p_1\mathbb{Z}}^{-1}\cdot x\big),
	\end{split}
\end{equation}
 and it corresponds, in physics language, to \emph{the Schrödinger-Pauli operator of spin $p/2$}. We denote its eigenvalues and orthonormal eigensections by $(\lambda_{(p_0,p_1),p,i},u_{(p_0,p_1),p,i})$ as in \eqref{a1}. Also, by \eqref{1.13n}, the averaging operator $\Delta^{X_{(p_0,p_1)}\times\mathbb{CP}^1}$ is given by
 \begin{equation}\label{1.35n}
 	\begin{split}
 		&\Delta^{X_{(p_0,p_1)}\times\mathbb{CP}^1}\colon L^2\big(X_{(p_0,p_1)}\times\mathbb{CP}^1\big)\to L^2\big(X_{(p_0,p_1)}\times\mathbb{CP}^1\big),\\
 		&\big(\Delta^{X_{(p_0,p_1)}\times\mathbb{CP}^1}f\big)(x,z)=\sum_{a\in A}f\big(a_{\mathbb{Z}/p_1\mathbb{Z}}^{-1}\cdot x,a_\mathbb{R}^{-1}\cdot z\big).
 	\end{split}
 \end{equation}
 We fix $p_0$. The prime $p_1$ plays the role of the large scale parameter $m$.

First, Ramanujan vector bundles achieve the optimal lower bound in \eqref{1.19n}.
\begin{theo}	
	If $p=0$ and $\lambda_{(p_0,p_1),0,i}\neq\pm(p_0+1)$, or if $p\neq 0$, then
	\begin{equation}\label{1.31}
		\lv\lambda_{(p_0,p_1),p,i}\rv\leqslant 2\sqrt{p_0}.
	\end{equation}
\end{theo}

Next, Ramanujan vector bundles satisfy the asymptotic spectral distribution in \eqref{1.20n}.
\begin{theo}
Let $I\subseteq\mathbb{R}$ be a fixed interval. Then
	\begin{equation}
		\begin{split}
		\limsup_{\max(p_1,p)\to\infty}\bbv&\frac{1}{\dim_\mathbb{C}F_{(p_0,p_1),p}}\lv\{i\mid\lambda_{(p_0,p_1),p,i}\in I\}\rv\\
		&-\int_{I}\frac{(p_0+1)\sqrt{4p_0-\lambda^2}}{2\pi\big((p_0+1)^2-\lambda^2\big)}\mathbbm{1}_{[-2\sqrt{p_0},2\sqrt{p_0}]\cap I}(\lambda)d\lambda\bbv=0,
		\end{split}
	\end{equation}
	where $\dim_\mathbb{C}F_{(p_0,p_1),p}$ is given in \eqref{1.32n}.
\end{theo}

Then, Ramanujan vector bundles satisfy the asymptotic log determinant in \eqref{n7.8n}.
\begin{theo}
We have
	\begin{equation}
		\begin{split}
&\limsup_{\max(p_1,p)\to\infty}\\
&\Bv\frac{1}{\dim_\mathbb{C}F_{(p_0,p_1),p}}\ln\det\big(p_0+1-\Delta^{F_{(p_0,p_1),p}}\big)\\
			&-\Big(p_0\ln(p_0)-\frac{1}{2}(p_0-1)\ln(p_0-1)-\frac{1}{2}(p_0-1)\ln(p_0+1)\Big)\Bv=0.
		\end{split}
	\end{equation}
\end{theo}

Moreover, we explain the strong uniform QE \eqref{9.} for Ramanujan vector bundles. Consider the $3$-sphere
\begin{equation}
\mathbb{S}^3=\{z=\lk\begin{smallmatrix}
	z_0\\
	z_1
\end{smallmatrix}\rk\in\mathbb{C}^2\mid\vert z_0\vert^2+\vert z_1\vert^2=1\}.
\end{equation}
By \eqref{1.28}, we have an embedding
\begin{equation}\label{1.10...}
	\mathrm{Sym}^p(\mathbb{C}^2)\hookrightarrow C^\infty(\mathbb{S}^3).
\end{equation}
For any $x\in X_{(p_0,p_1)}$ and any eigensection $u_{(p_0,p_1),p,i}$, we have $u_{(p_0,p_1),p,i}(x)\in \mathrm{Sym}^p(\mathbb{C}^2)$. Thus we can view $u_{(p_0,p_1),p,i}\in C^\infty(X_{(p_0,p_1)}\times\mathbb{S}^3)$ by
\begin{equation}\label{1.14''}
	u_{(p_0,p_1),p,i}\colon (x,z)\in X_{(p_0,p_1)}\times\mathbb{S}^3\mapsto\big(u_{(p_0,p_1),p,i}(x)\big)(z).
\end{equation}
It follows that $\lv u_{(p_0,p_1),p,i}\rv_\mathbb{C}^2\in C^\infty(X_{(p_0,p_1)}\times \mathbb{S}^3)$. Since homogeneity implies $\vert u_{(p_0,p_1),p,i}(x,e^{i\theta}z)\vert_\mathbb{C}^2=\lv u_{(p_0,p_1),p,i}(x,z)\rv_\mathbb{C}^2$ for any unit complex number $e^{i\theta}\in\mathbb{S}^1$, the function $\lv u_{(p_0,p_1),p,i}\rv_\mathbb{C}^2$ descends to $C^\infty(X_{(p_0,p_1)}\times \mathbb{CP}^1)$ via \emph{the Hopf fibration} $\mathbb{CP}^1\cong \mathbb{S}^1\backslash\mathbb{S}^3$. Then \eqref{1.10} becomes the measure
\begin{equation}
	\begin{split}
		&\mathscr{Q}_{(p_0,p_1)}\in C^\infty(X_{(p_0,p_1)}\times \mathbb{CP}^1)\\
		&\mapsto\sum_{x\in X_{(p_0,p_1)}}\int_{\mathbb{CP}^1}\mathscr{Q}_{(p_0,p_1)}(x,z)\lv u_{(p_0,p_1),p,i}(x,z)\rv_{\mathbb{C}}^2\omega_{\mathrm{FS}}(z)\in\mathbb{C}
	\end{split}
\end{equation}
on $X_{(p_0,p_1)}\times \mathbb{CP}^1$.

\begin{theo}
For any series $\big(\mathscr{Q}_{(p_0,p_1)}\in C^\infty(X_{(p_0,p_1)}\times \mathbb{CP}^1)\big)_{p_1}$ satisfying $\lV\mathscr{Q}_{(p_0,p_1)}\rV_{L^\infty(X_{(p_0,p_1)}\times \mathbb{CP}^1)}\leqslant C$, we have
	\begin{equation}
		\begin{split}
			\limsup_{\max(p_1,p)\to\infty}&\frac{1}{\dim_\mathbb{C}F_{(p_0,p_1),p}}\sum_{i=1}^{\dim_\mathbb{C}F_{(p_0,p_1),p}}\\
			&\bbv\sum_{x\in X_{(p_0,p_1)}}\int_{\mathbb{CP}^1}\mathscr{Q}_{(p_0,p_1)}(x,z)\lv u_{(p_0,p_1),p,i}(x,z)\rv_\mathbb{C}^2\omega_{\mathrm{FS}}(z)\\
			&-\frac{1}{\lv X_{(p_0,p_1)}\rv}\sum_{x\in X_{(p_0,p_1)}}\int_{\mathbb{CP}^1}\mathscr{Q}_{(p_0,p_1)}(x,z)\omega_{\mathrm{FS}}(z)\bbv=0,
		\end{split}
	\end{equation}
	where $\vert X_{(p_0,p_1)}\vert$ and $\dim_\mathbb{C}F_{(p_0,p_1),p}$ are given in \eqref{1.32n}.
\end{theo}

In addition, the zero divisor of eigensections of Ramanujan vector bundles satisfy property \eqref{1.25}. In particular, since $\dim_{\mathbb{C}}\mathbb{CP}^1=1$, the zero divisor reduces to the counting measure of zeros with multiplicity.
\begin{theo}
Fix a pair $(p_0,p_1)$. Then there exists a series of subsets $\big(I_{(p_0,p_1),p}\subseteq \{1,\cdots, \dim_{\mathbb{C}}F_{(p_0,p_1),p}\}\big)_{p\in\mathbb{N}}$ of asymptotic density one, in the sense that
\begin{equation}
	\lim_{p\to\infty}\lv I_{(p_0,p_1),p}\rv\big/\dim_{\mathbb{C}}F_{(p_0,p_1),p}=1,
\end{equation}
	where $\dim_\mathbb{C}F_{(p_0,p_1),p}$ is given in \eqref{1.32n}, such that for any series $(i_p\in I_{(p_1,p_2),p})_{p\in\mathbb{N}}$, and for any $\mathscr{Q}_{(p_0,p_1)}\in C^\infty(X_{(p_0,p_1)}\times \mathbb{CP}^1)$, we have
\begin{equation}\label{1.36}
	\begin{split}
	\lim_{p\to\infty}\frac{1}{p}&\sum_{\substack{(x,z)\in X_{(p_0,p_1)}\times \mathbb{CP}^1,\\ u_{(p_0,p_1),p,i_p}(x,z)=0}}\mathscr{Q}_{(p_0,p_1)}(x,z)\\
	=&\sum_{x\in X_{(p_0,p_1)}}\int_{\mathbb{CP}^1}\mathscr{Q}_{(p_0,p_1)}(x,z)\omega_{\mathrm{FS}}(z).
	\end{split}
\end{equation}
\end{theo}

By \eqref{1.36}, we see that $u_{(p_0,p_1),p,i_p}$ has $p\vert X_{(p_0,p_1)}\vert+o(p)$ zeros, and that these zeros become equidistributed on $X_{(p_0,p_1)}\times\mathbb{CP}^1$.

\subsection{Mixed quantization and spin-scale duality}

Our main difficulty lies in controlling eigensections of infinitely many vector bundles over infinitely many graphs simultaneously. To address this, we use the technique of \emph{mixed quantization}, originally developed in the context of heat kernel analysis in \cite{MR2838248,MR3615411,MR4665497,MR4611826}, and later extended to pseudodifferential calculus in \cite{MR4808253,ovadia2025mixedquantizationpartialhyperbolicity}. Mixed quantization combines \emph{Weyl quantization}, which controls the high frequency limit, with \emph{Berezin-Toeplitz quantization}, which controls the large spin limit. For background on Weyl quantization, we refer to Zworski \cite[\S\,4]{MR2952218} and Dyatlov-Zworski \cite[Appendix E]{MR3969938}, and for Berezin-Toeplitz quantization, we refer to Ma-Marinescu \cite[\S\,7]{MR2339952}. In the current discrete setting, Weyl quantization on manifolds is replaced by \emph{kernel operators} for graph vector bundles, following Anantharaman \cite[\S\S\,2, 3]{MR3649482}, to control the large scale limit. There is a precise analogy, that the space of nonbacktracking paths is the graph theoretic counterpart of the cotangent phase space in the manifold setting, and the nonbacktracking walk is the counterpart of the geodesic flow, see \S\S\,\ref{s7.1n} and \ref{s7.3n} for more details.

Although the mixed quantization formalism is developed in a general setting, the most important examples behind Theorems \ref{t1.1}, \ref{t1.2}, \ref{C9'}, and \ref{t1.4} arise from representations of compact Lie groups. For these special cases, we also provide independent proofs using representation theory. In particular, throughout the paper, the large spin limit relies crucially on the off-diagonal exponential decay estimate for Bergman kernels by Ma-Marinescu \cite[Theorem 1]{MR3368102}. This estimate, in turn, depends on the deep analytic localization technique of Bismut-Lebeau \cite[\S\,11]{MR1188532}. We use this general result as a black box. In the special case of compact Lie group representations, we also give an alternative proof, and see \S\,\ref{s3.4} for more details.

The above discussion indicates a parallel between the large scale and large spin limits. Within the framework of mixed quantization, they can be understood through a common underlying mechanism. We refer to this phenomenon as \emph{spin-scale duality}. It may be viewed as a discrete counterpart of \emph{wave-particle duality} in quantum mechanics, and see \S\,\ref{s4.4n} for further discussions.

\subsection{Related works}

 We now turn to results related to our main theorems. 

 \subsubsection{Discrete QE}
 
Let us begin with two extreme cases of \eqref{1.8}. 

In the first case, the fibre is trivial, we take $(N,L)=(\mathrm{pt},\mathbb{C})$, namely a point (a $0$-dimensional complex manifold) with the trivial line bundle, while the base graph remains nontrivial. That is, we take a series $(X_m)_{m\in\mathbb{N}}$ of graphs satisfying \textbf{[BST]} \eqref{bst} and \textbf{[EXP]} \eqref{exp}, and denote this case by
\begin{equation}\label{1.40}
	\big(X_m,F_m\big)=\big(X_m,(\mathbb{C},\phi=1)\big).
\end{equation}
In this case, the vector bundle Laplacian is simply the usual graph Laplacian. 

In the second case, the graph is as simple as possible, a bouquet graph, consisting of a single vertex together with several self-loops $\circlearrowleft_{i=1}^k$, while the fibre is nontrivial. We take the geometric realization \eqref{1.26n}, assign to each $\circlearrowleft_{i}$ an element $g_i\in\mathrm{SU}_2(\mathbb{C})$, and denote this case by
\begin{equation}\label{1.41}
	\big(X,F_p\big)=\Big(\big(\mathrm{pt},\circlearrowleft_{i=1}^k\big),\big(\mathrm{Sym}^p(\mathbb{C}^2),\phi_{\circlearrowleft_i}=g_i\big)\Big).
\end{equation}
In this case, the vector bundle Laplacian $\Delta^{F_p}$ becomes an \emph{averaging operator}
\begin{equation}\label{1.41n}
	(\Delta^{F_p}s)(z)=\sum_{i}\big(s(g_iz)+s\big(g_i^{-1}z\big)\big).
\end{equation}

Applying Theorem \ref{C9'} to \eqref{1.40} and \eqref{1.41}, we obtain the results of Anantharaman-Le Masson \cite[Theorem 1]{MR3322309} and Brooks-Le Masson-Lindenstrauss \cite[Theorem 1]{MR3567266}, corresponding respectively to the large scale limit and the large spin limit. Note that, under the setting \eqref{1.40}, our spectral gap condition \textbf{[EXP]} \eqref{exp} is slightly weaker than the one in Anantharaman-Le Masson \cite[\S\,1]{MR3322309}, in that it does not require a gap at $(-d)$, and therefore applies to additional examples, including bipartite graphs. Moreover, \cite[Theorem 1]{MR3567266} uses a slightly different geometric representation of $\mathrm{SO}_3(\mathbb{R})$ via spherical harmonics, corresponding to even values of $p$ in \eqref{1.28} through the double covering $\mathrm{SU}_2(\mathbb{C})\to\mathrm{SO}_3(\mathbb{R})$, this case can nevertheless be treated in essentially the same way. In addition, \cite[\S\S\,2, 3]{MR3567266} provides a new proof of \cite[Theorem 1]{MR3322309}, and the similarity between the large scale and large spin limits was already observed there. Indeed, placing these two limits within a unified framework is one of our original motivations.

\subsubsection{QE for manifolds}

QE for manifolds was established by Shnirelman \cite{MR0402834}, Zelditch \cite{MR916129}, and Colin de Verdière \cite{MR818831}. 

Using Fourier integral operator techniques, Schrader-Taylor \cite{MR0995750} proved nonuniform QE for vector bundles over manifolds. In \cite{MR4808253,ovadia2025mixedquantizationpartialhyperbolicity}, uniform QE for flat vector bundles over manifolds, as well as its stability under general perturbations, was studied. Subsequently, Cekić-Lefeuvre \cite{cekić2024semiclassicalanalysisprincipalbundles} obtained nonuniform QE for flat vector bundles over manifolds and perturbations arising from principal bundles.

\subsubsection{Asymptotic torsion}


In addition to \cite{MR2838248,MR3615411,MR4665497,MR4611826}, we mention the asymptotic torsion results of Bismut-Vasserot \cite{MR1016875,MR1096593}, Bergeron-Venkatesh \cite{MR3028790}, Calegari-Venkatesh \cite{MR3961523}, Müller \cite{MR3220447}, Finski \cite{MR3864507}, and Liu \cite{MR4746871}.

\subsubsection{Laplacian determinant and Graph vector bundles}

Lyons \cite{MR2160416} estimated the number of spanning trees of large graphs by analyzing asymptotic normalized determinants of graph Laplacians.

Kenyon \cite[\S\,3]{MR2884879} introduced graph vector bundles and used determinants of graph vector bundle Laplacians to study various probabilistic questions, including quantities related to loop-erased random walks that cannot be computed from the standard graph Laplacian.

\subsubsection{Spectral distribution}

The Alon-Boppana bound was established in Alon \cite[Page 95]{MR875835} and \cite[Theorem 1]{MR1124768}. The Kesten-McKay law was obtained by Kesten \cite[Theorem 3]{MR109367} and McKay \cite[Theorem 1.1]{MR629617}. Applying Theorem \ref{t1.2} to \eqref{1.41}, we recover Gamburd-Jakobson-Sarnak \cite[Proposition 4.1]{MR1677685} and \cite[Corollary 2]{MR3567266}. Moreover, Theorem \ref{t1.2} takes a form particularly similar to the eigenvalue density of Hecke operators by Serre \cite[Théorème 1]{MR1396897}, see also related results of Sarnak \cite[Theorem 1.2]{MR1018385} and Conrey-Duke-Farmer \cite[Theorem 1]{MR1438595}.
This similarity suggests that the large spin part of Theorem \ref{C9'} may be viewed as an analogue of the QE for Hecke operators, as in Luo-Sarnak \cite[Theorem 1.1]{MR1990480}. By the strong approximation for $\mathrm{PGL}_2$, we have for a prime $p_0$,
\begin{equation}
	\mathrm{PGL}_2(\mathbb{Z})\big\backslash \mathbb{H}^2\cong	\mathrm{PGL}_2\big(\mathbb{Z}\big[\tfrac{1}{p_0}\big]\big)\big\backslash \big(\mathbb{H}^2\times\big(\mathrm{PGL}_2(\mathbb{Q}_{p_0})\big/\mathrm{PGL}(\mathbb{Z}_{p_0})\big)\big),
\end{equation}
where the right hand side can be viewed as a family of \emph{the Bruhat-Tits trees} $\mathrm{PGL}_2(\mathbb{Q}_{p_0})\big/\mathrm{PGL}(\mathbb{Z}_{p_0})$ over the base $\mathrm{PGL}_2\Big(\mathbb{Z}\Big[\tfrac{1}{p_0}\Big]\Big)\big\backslash\mathbb{H}^2$. Therefore, certain number theoretic results for Hecke operators may be viewed as family versions of graph theoretic results.

\subsubsection{Random holomorphic sections}

In Shiffman-Zelditch \cite[Theorems 1.1, 1.3]{MR1675133} and Nonnenmacher-Voros \cite[Theorem 1]{MR1649013}, it was shown for almost every random series $(u_p\in H^{(0,0)}(N,L^p))_{p\in\mathbb{N}}$, sampled from a suitable statistical ensemble, we have $\lim_{p\to\infty}\mathrm{Div}(u_p)/p=c_1(L,h^L)$. Thus, Theorem \ref{t1.4} supports a familiar principle in quantum chaos, when a quantum system corresponds to a classically \emph{chaotic} dynamical system, its quantum states often behave like \emph{random} states. See Drewitz-Liu-Marinescu \cite{MR4922234,MR4905035} for more results on zeros of random holomorphic sections.

\subsubsection{Spectral gap and Ramanujan graphs}

Examples satisfying \textbf{[EXP]} \eqref{exp} are difficult to construct even for two extreme cases \eqref{1.40} and \eqref{1.41}. For \eqref{1.40}, the spectral gap assumption requires that $(X_m)_{m\in\mathbb{N}}$ are expanders, and the first explicit example was constructed by Margulis \cite{MR484767}. For \eqref{1.41}, the spectral gap is related to the \emph{Banach-Ruziewicz problem}, and the first example was given by Margulis \cite{MR596890}, Sullivan \cite{MR590825}, and Drinfeld \cite{MR757256}. For more constructions, see Gamburd-Jakobson-Sarnak \cite{MR1677685} and Bourgain-Gamburd \cite{MR2358056,MR2966656}.

Ramanujan graphs $X_{(p_0,p_1)}$ in \eqref{1.30n} were introduced by Lubotzky-Phillips-Sarnak \cite{MR963118} and Margulis \cite{MR939574} as optimal expanders. Ramanujan vector bundles $F_{(p_0,p_1),p}$ in \eqref{1.30n}, together with their remarkable spectral gap property \eqref{1.31}, were established in the same series of papers \cite{MR861487,MR890171}, even earlier than Ramanujan graphs! We note that the spectral gap for Ramanujan vector bundles also has a Banach-Ruziewicz type consequence, namely that the natural product probability measure on $X_{(p_0,p_1)}\times \mathbb{CP}^1$ is the unique finitely additive probability measure on all subsets that is invariant under the averaging operator defined in \eqref{1.35n}.

Expander graphs and Ramanujan graphs have become a major subject with numerous applications in mathematics and computer science. We refer to the surveys of Hoory-Linial-Wigderson \cite{MR2247919} and Li \cite{MR4055707} for more information. However, expander graph vector bundles and Ramanujan vector bundles seem to have received less attention and remain far from fully developed.

Marcus-Spielman-Srivastava \cite{MR3374962} proved that any bipartite Ramanujan graph admits a Ramanujan $2$-covering. This was generalized by Hall-Puder-Sawin \cite{MR3725881}, who showed that a Ramanujan graph admits a Ramanujan $k$-covering for any $k\in\mathbb{N}$. We note that the covering graph setting fits naturally into the framework of graph vector bundles, since a $k$-covering may be viewed as a discrete bundle whose fibre consists of $k$ isolated points, with gauge group given by the permutation group $\mathfrak{S}_k$.

\subsection{Organization of the paper}

In \S\,\ref{s1}, we define graph vector bundles and give an equivalent description in terms of flat vector bundles. In \S\,\ref{s3n}, we introduce complex geometry, Bergman kernels, and geometric representations. In \S\,\ref{s4n}, we discuss the large scale, large spin, and strong uniform limits, and explain the spin-scale duality underlying them. In \S\,\ref{s5n}, we state and prove the Alon-Boppana bound. In \S\,\ref{s6n}, we state and prove the Kesten-McKay law. In \S\,\ref{det}, we study asymptotic log determinants of Laplacians. In \S\,\ref{s7n}, we introduce kernel operators on graph vector bundles, which form a discrete analogue of semiclassical quantization on manifolds. In \S\,\ref{s8n}, we introduce Berezin-Toeplitz quantization, a form of geometric quantization. In \S\,\ref{s9n}, we combine \S\S\,\ref{s7n} and \ref{s8n} to develop mixed quantization. In \S\,\ref{s10}, we state and prove quantum ergodicity. In \S\,\ref{s11}, we state and prove zero divisor equidistribution. In \S\,\ref{s12}, we introduce general Ramanujan vector bundles and describe various questions and conjectures motivated by the rest of the paper.

\subsection{Acknowledgments}

I am grateful to Sherry Gong, Zhizhang Xie, and Guoliang Yu for support and understanding. I would like to thank Dean Baskin, Semyon Dyatlov, Yulin Gong, Elena Kim, and Xiaonan Ma for helpful discussions on an early version of this manuscript. I would like to thank Wen-Ching Winnie Li for bringing \cite{MR861487,MR890171,MR963118} to my attention. This work was supported by NSF grants DMS-1546917 and DMS-2247313. Finally, I owe special thanks to my wife and the wonderful gift of our newborn, without whose nighttime cries, I would sleep well each night and would not have had the time to work on mathematics.

\subsection{Declarations} 

The author have no relevant financial or nonfinancial interests to disclose.  Data sharing does not apply to this article as no datasets were generated or analyzed during the current study.

\section{Graph Vector Bundles}\label{s1}

In this section, we discuss vector bundles on graphs. In \S\,\ref{s2.1}, we introduce the twisted Laplacian approach. In \S\,\ref{s2.2}, we present the flat vector bundle approach, analogous to the manifold setting. In \S\,\ref{s2.3}, we prove that these two approaches are equivalent and discuss the advantages of the latter. Throughout \S\S\,\ref{s2.1}, \ref{s2.2}, and \ref{s2.3}, we allow graphs that are not necessarily regular and vector bundles that are not necessarily unitary. In \S\,\ref{s2.4}, we specialize to unitary vector bundles over regular graphs, which is the main setting for the rest of the paper.

For the definition of graph vector bundles, we refer to Kenyon \cite[\S\,3]{MR2884879}. For the fundamental group of graphs, we refer to Hatcher \cite[\S\,1.A]{MR1867354}. For flat vector bundles over manifolds, we refer to Kobayashi \cite[\S\,\uppercase\expandafter{\romannumeral1}.2]{MR909698}.

\subsection{Graph vector bundles and twisted Laplacian}\label{s2.1}

Let $X$ be a graph. An $\ell$-dimensional complex vector bundle $F=(X\times \mathbb{C}^\ell,\phi)$ over $X$ consists of the total space $X\times \mathbb{C}^\ell$ and a connection $\phi$, where $\phi$ assigns to each oriented edge $(x,x')$ of $X$ a linear isomorphism $\phi_{x,x'}\colon\mathbb{C}^\ell_{x'}\to \mathbb{C}^\ell_{x}$, called the parallel transport from $\mathbb{C}^\ell_{x'}$ to $\mathbb{C}^\ell_{x}$, with the consistency condition
\begin{equation}\label{2.1}
	\phi_{x,x'}=\phi_{x',x}^{-1}.
\end{equation}
We define the associated \emph{twisted} Laplacian operator $\Delta^\phi$ by
\begin{equation}\label{2.2}
	\begin{split}
&\Delta^\phi\colon C^\infty(X,\mathbb{C}^\ell)\rightarrow C^\infty(X,\mathbb{C}^\ell),\\
&(\Delta^\phi u)(x)=\sum_{x'\sim x}\phi_{x,x'}u(x').
	\end{split}
\end{equation}

We say that $F$ is trivial if $\phi_{x,x'}=\mathrm{Id}_{\mathbb{C}^\ell}$ for any $x,x'\in X$.

Let $F=(X\times \mathbb{C}^\ell,\phi),F'=(X\times \mathbb{C}^\ell,\phi')$ be two vector bundles over $X$ with the same fibre dimension. We say that they are \emph{gauge equivalent} if there is for each vertex $x\in X$ a linear invertible morphism $\psi_x\in \mathrm{End}(\mathbb{C}^\ell_x)$ such that for any $x,x'\in X$,
\begin{equation}
\psi_{x}\phi_{x,x'}=\phi_{x,x'}'\psi_{x'}.
\end{equation}
In other words, the following diagram
\begin{equation}\label{2.4}
	\begin{tikzcd}[ampersand replacement=\&, column sep=normal,row sep=normal]
		\mathbb{C}^\ell_{x'}\arrow[r,"\phi_{x,x'}"]\arrow[d,"\psi_{x'}"swap]\& \mathbb{C}^\ell_{x} \arrow[d,"\psi_{x}"] \\
		\mathbb{C}^\ell_{x'} \arrow[r,"\phi_{x,x'}'"swap]\& \mathbb{C}^\ell_{x}
	\end{tikzcd}
\end{equation}
is commutative.

\subsection{Universal covering and fundamental group}\label{s2.2}

Fix an arbitrary vertex $x_0 \in X$ as a base point. We construct the associated universal covering $\widetilde{(X,x_0)}$ of $X$ as follows. Each vertex of $\widetilde{(X,x_0)}$ represents a nonbacktracking path that starts from $x_0$, that is, a sequence $y=(x_j,\cdots,x_1,x_0)$ such that $x_i\sim x_{i+1}$ and $x_{i-1}\neq x_{i+1}$ for any $1\leqslant i\leqslant j-1$. Two vertices of $\widetilde{X}$ are adjacent if one is a simple extension of another, that is, $y=(x_j,\cdots,x_1,x_0)$ and $y'=(x_{j+1},x_j,\cdots,x_1,x_0)$ are adjacent. We have a natural projection
\begin{equation}
y=(x_j,\cdots,x_0)\in\widetilde{(X,x_0)}\mapsto x_j \in X.
\end{equation}

We can also describe $\widetilde{(X,x_0)}$ as the space of all paths modulo homotopy equivalence, where two paths $y=(x_j,\cdots,x_{i+1},x_i,x_{i-1},\cdots,x_0)$ and $y'=(x_j,\cdots,x_{i+2},x_{i-1},\cdots,x_0)$ are identified whenever $x_{i-1}=x_{i+1}$, that is, when a backtracking step is added or removed. 

Let us define the loop space
\begin{equation}\label{2.7n}
\Omega(X,x_0)=\big\{\text{paths of }X \text{ starting and ending at } x_0\big\}.
\end{equation}
Then the fundamental group $\pi_1(X,x_0)$ of $X$ based at $x_0$ satisfies
\begin{equation}
\pi_1(X,x_0)=\text{homotopy equivalence}\big\backslash\Omega(X,x_0).
\end{equation}
The group $\pi_1(X,x_0)$ acts on $\widetilde{(X,x_0)}$ as the composition of paths, namely $\gamma(y)=(\gamma,y)$, and we have
\begin{equation}\label{2.5}
\pi_1(X,x_0)\big\backslash \widetilde{(X,x_0)}\cong X.
\end{equation}

The following result describes the group structure of $\pi_1(X,x_0)$.

\begin{lemma}
The fundamental group $\pi_1(X,x_0)$ is a finitely generated free group.
\end{lemma}

\begin{pro}
Choose a maximal tree $T\subseteq X$. Each edge $(x_i,x_i')$ of $X-T$ determines a loop $\gamma_i$ in $X$ that goes first from $x_0$ to one endpoint of $(x_i,x_i')$ by a path in $T$, then across $(x_i,x_i')$, then back to $x_0$ by a path in $T$. Then $\pi_1(X)$ is the free group generated by homotopy classes $[\gamma_i]$ corresponding to all edges $(x_i,x_i')$ of $X-T$.\qed
\end{pro}

Using a maximal tree $T\subseteq X$, we can construct \emph{a fundamental domain} $D(x_0)$
as follows. For any $x\in X$, there is a unique path $(x_j,\cdots,x_0)$ on $T$ from $x_0$ to $x_j=x$. Then let 
\begin{equation}\label{2.10n}
	\begin{split}
D(x_0)=\text{the induced subgraph with these paths as vertices}&\\
\subset \widetilde{(X,x_0)}.&
	\end{split}
\end{equation}
Its vertex set is a set of representatives for the quotient \eqref{2.5}.

In what follows, we shall omit the starting point $x_0$ and simply denote
\begin{equation}\label{2.11n}
(\pi_1(X),\widetilde{X},D)=(\pi_1(X,x_0), \widetilde{(X,x_0)},D(x_0)).
\end{equation}

\subsection{Graph vector bundles are all flat}\label{s2.3}

In differential geometry, a vector bundle is said to be flat if it is equipped with a connection whose curvature vanishes. Flat vector bundles also admit a representation theoretic description, arising from linear representations of the fundamental group of the base manifold.

Analogously, we can define flat vector bundles over graphs. Given a linear representation 
\begin{equation}\label{2.6}
	\rho\colon \pi_1(X)\to \mathrm{End}(\mathbb{C}^\ell),
\end{equation}
we can associate a flat vector bundle $F$ over $X$ by
\begin{equation}\label{2.7}
	\begin{split}
		F&=\pi_1(X)\backslash\big(\widetilde{X}\times\mathbb{C}^\ell\big)\\
		&=\big\{(y,v)\in\widetilde{X}\times\mathbb{C}^\ell\big\}/(y,v)\sim(\gamma y,\rho(\gamma)v)\ \text{for any }\gamma\in\pi_1(X).
	\end{split}
\end{equation}
Here the symbol $\sim$ denotes the above equivalence relation and should not be confused with the adjacency relation between vertices. We recall the dimension conventions
\begin{equation}\label{2.14n}
	\dim_{\mathbb{C}}F_x=\ell,\quad\dim_{\mathbb{C}}F=\ell\lv X\rv.
\end{equation}

A section $u$ of $F$ over $X$ assigns to each $x\in X$ a vector $u(x)\in F_x$, and we denote the space of sections by $C^\infty(X,F)$. Let $C^\infty(\widetilde{X},\mathbb{C}^\ell)$ be the space of $\mathbb{C}^\ell$-valued functions on $\widetilde{X}$, and we have a $\pi_1(X)$-action on it by
\begin{equation}\label{2.8}
	\big(\gamma(u)\big)(y)=\rho(\gamma)u(\gamma^{-1}y), \quad\gamma\in\pi_1(X), u\in C^\infty(\widetilde{X},\mathbb{C}^\ell).
\end{equation}
Then by \eqref{2.7}, the space $C^\infty(X,F)$ of sections of $F$ over $X$ is isomorphic to the $\pi_1(X)$-invariant subspace $C^\infty_{\pi_1(X)}(\widetilde{X},\mathbb{C}^\ell)\subset C^\infty(\widetilde{X},\mathbb{C}^\ell)$, that is,
\begin{equation}\label{2.9}
	C^\infty(X,F)\cong C^\infty_{\pi_1(X)}(\widetilde{X},\mathbb{C}^\ell).
\end{equation}

Let $\Delta^{\mathbb{C}^\ell}$ be the Laplacian operator acting on $C^\infty(\widetilde{X},\mathbb{C}^\ell)$ given by
\begin{equation}\label{2.10}
(\Delta^{\mathbb{C}^\ell}u)(y)=\sum_{y'\sim y}u(y').
\end{equation}
Since $y\sim y'$ if and only if $\gamma y\sim \gamma y'$ for any $\gamma\in\pi_1(X)$, we see that $\Delta^{\mathbb{C}^\ell}$ commutes with the $\pi_1(X)$-action given in \eqref{2.8}, in other words,
\begin{equation}
\gamma\Delta^{\mathbb{C}^\ell}u=\Delta^{\mathbb{C}^\ell}\gamma u.
\end{equation}
Therefore, $\Delta^{\mathbb{C}^\ell}$ defines an operator acting on $C^\infty_{\pi_1(X)}(\widetilde{X},\mathbb{C}^\ell)$. By \eqref{2.9}, $\Delta^{\mathbb{C}^\ell}$ descends to an operator acting on $C^\infty(X,F)$, which we shall denote by $\Delta^F$.

\begin{prop}\label{p2.2}
Modulo gauge equivalence, any graph vector bundle is flat.
\end{prop}

\begin{pro}
Let us choose a maximal tree $T\subseteq X$ and a base vertex $x_0\in T$. A key observation is that, over $T$, any vector bundle is gauge equivalent to a trivial vector bundle. To see this, we use a gauge transform defined as follows. For any $x\in T$, let $(x_j,\cdots,x_0)$ be the unique path on $T$ from $x_0$ to $x_j=x$, then define $\psi_x=\phi_{x_0,x_{1}}\cdots \phi_{x_{j-1},x_j}$. From \eqref{2.1}, the diagram \eqref{2.4} extends to
	\begin{equation}
		\begin{tikzcd}[ampersand replacement=\&, column sep=large,row sep=normal]
			\mathbb{C}^\ell_{x_0}\arrow[r,"\phi_{x_1,x_0}"]\arrow[d,"\mathrm{Id}_{\mathbb{C}^\ell}"swap]\& \mathbb{C}^\ell_{x_1} \arrow[d,"\phi_{x_0,x_{1}}"swap] \arrow[r,"\phi_{x_2,x_1}"] \& \mathbb{C}^\ell_{x_2} \arrow[d,"\phi_{x_0,x_{1}}\phi_{x_1,x_2}"swap] \arrow[r,"\phi_{x_3,x_2}"]\& \cdots \arrow[r,"\phi_{x_j,x_{j-1}}"]  \& \mathbb{C}^\ell_{x_j} \arrow[d,"\phi_{x_0,x_{1}}\cdots \phi_{x_{j-1},x_j}"swap]  \\
			\mathbb{C}^\ell_{x_0} \arrow[r,"\mathrm{Id}_{\mathbb{C}^\ell}"swap]\& \mathbb{C}^\ell_{x_1}\arrow[r,"\mathrm{Id}_{\mathbb{C}^\ell}"swap] \&  \mathbb{C}^\ell_{x_2}\arrow[r,"\mathrm{Id}_{\mathbb{C}^\ell}"swap] \& \cdots\arrow[r,"\mathrm{Id}_{\mathbb{C}^\ell}"swap] \& \mathbb{C}^\ell_{x_j}
		\end{tikzcd}
	\end{equation}
	and we obtain a trivialization over $T$. We denote this new connection 
	by $\phi'$.
	
Recall that each edge $(x_i,x_i')$ of $X-T$ determines a loop $\gamma_i$ in $X$, and $\pi_1(X)$ is the free group generated by homotopy classes $[\gamma_i]$ corresponding to all edges $(x_i,x_i')$ of $X-T$. Then the required representation \eqref{2.6} can be defined by $\rho([\gamma_i])=\phi'_{x_i,x_i'}$.\qed
\end{pro}

This may also be interpreted from a differential geometric viewpoint. Curvature is antisymmetric, whereas a graph is one-dimensional, so no nontrivial curvature can appear. Hence, every graph vector bundle is flat.

In the graph vector bundle formulation $(C^\infty(X,\mathbb{C}^\ell),\Delta^\phi)$ discussed in \S\,\ref{s2.1}, the section space $C^\infty(X,\mathbb{C}^\ell)$ is flat, while the Laplacian $\Delta^\phi$ defined in \eqref{2.2} is twisted. In contrast, the flat vector bundle formulation $(C^\infty(X,F),\Delta^F)$ presented in this subsection features a flat Laplacian $\Delta^F$ given in \eqref{2.10}, with the twisting instead absorbed into the section space $C^\infty(X,F)$. Proposition \ref{p2.2} shows that these two formulations are equivalent. The advantage of absorbing the twisting into the section space $C^\infty(X,F)$ is that, later, we will need to compose several operators. Treating these operators as flat simplifies the computations.

\subsection{Unitary vector bundles over regular graphs}\label{s2.4}

In what follows, we shall focus on the case where $X$ is a $d$-regular graph and $F$ is a unitary vector bundle over $X$. 

By discussions in \S\S\,\ref{s2.2} and \ref{s2.3}, we now have $\widetilde{X}\cong \mathbb{T}_d$, the $d$-regular tree. Let us denote $\Gamma=\pi_1(X)$ in accordance with the standard notation in hyperbolic geometry. Then \eqref{2.5}, \eqref{2.10n}, \eqref{2.11n}, \eqref{2.6}, and \eqref{2.7} become
\begin{equation}\label{2.19n}
	\begin{split}
	\big(\pi_1(&X), \widetilde{X}, X, D, \rho\big)=\big(\Gamma,\mathbb{T}_d,\Gamma\backslash \mathbb{T}_d, D, \rho\colon \Gamma\to \mathrm{U}_\ell(\mathbb{C})\big),\\
	F=&\Gamma\backslash\big(\mathbb{T}_d\times\mathbb{C}^\ell\big)\\
	=&\big\{(y,v)\in\mathbb{T}_d\times\mathbb{C}^\ell\big\}/(y,v)\sim(\gamma y,\rho(\gamma)v)\ \text{for any }\gamma\in\Gamma.
\end{split}
\end{equation}
Also, \eqref{2.8}, \eqref{2.9}, and \eqref{2.10} take the form
\begin{equation}\label{2.20n}
	\begin{split}
	&C^\infty(X,F)\\
	&=C^\infty_{\Gamma}\big(\mathbb{T}_{d},\mathbb{C}^\ell\big)\subset C^\infty(\mathbb{T}_{d},\mathbb{C}^\ell)\\
	&=\big\{u \mid u(y)=\rho(\gamma)u(\gamma^{-1}y)\text{ for any } y\in\mathbb{T}_{d}, \gamma\in\Gamma\big\},\\
	&(\Delta^Fu)(y)=\sum_{y'\sim y}u(y').
	\end{split}
\end{equation}

Since $\rho$ is unitary, the standard Hermitian norm $\lV\cdot\rV_{\mathbb{C}^\ell}$ on $\mathbb{C}^\ell$ induces a norm $\lV\cdot\rV_F$ on $F$. Let $\lV\cdot\rV_{L^2(X,F)}$ be the $L^2$-norm on $C^\infty(X,F)$ induced by $\lV\cdot\rV_F$, in other words,
\begin{equation}\label{2.21n}
\lV u\rV_{L^2(X,F)}^2=\sum_{x\in X}\lV u(x)\rV_{F}^2=\sum_{y\in D}\lV u(y)\rV_{\mathbb{C}^\ell}^2.
\end{equation}
Moreover, the Laplacian $\Delta^F$ is selfadjoint with respect to $L^2(X,F)$. Hence we can list all the eigenvalues, counted with multiplicity, and choose an orthonormal basis of eigensections
\begin{equation}\label{3.27}
	\Delta^{F}u_{F,i}=\lambda_{F,i}u_{F,i},\quad\lV u_{F,i}\rV^2_{L^2(X,F)}=1.
\end{equation}

We study these eigenvalues and eigensections in more detail in later sections. We begin here with a simple, fundamental trace formula that will be used repeatedly in the analysis of eigenvalues.
\begin{lemma}
	For a flat vector bundle $F$ over $X$ of the form \eqref{2.19n}, we have for any $k\in\mathbb{N}$,
	\begin{equation}\label{2.23n}
		\mathrm{Tr}^{F}\big[\big(\Delta^{F}\big)^k\big]=\sum_{i=1}^{\dim_\mathbb{C}F}\lambda_{F,i}^k=\sum_{x\in X,\gamma\in \Omega(X,x),\lv\gamma\rv=k}\mathrm{Tr}^{F_{x}}\big[\rho(\gamma)\big],
	\end{equation}
	where $\Omega(X,x)$ is the loop space defined in \eqref{2.7n} and $\lv\gamma\rv$ denotes the length of $\gamma$.
\end{lemma}

\begin{pro}
By \eqref{2.20n}, we have
\begin{equation}
(\Delta^F)^ku(y)=\sum_{\text{paths } (y_k,\cdots,y_0), y_0=y}u(y_k).
\end{equation} 
To compute the trace, it suffices to consider closed paths on $X$, equivalently, $y_k=\gamma y$ for some $\gamma\in \Gamma$. Using \eqref{2.20n} again, we have $u(\gamma y)=\rho(\gamma)u(y)$. Hence, we sum over the fundamental domain $D$ to obtain
	\begin{equation}
	\mathrm{Tr}^{F}\big[\big(\Delta^{F}\big)^k\big]=\sum_{\substack{y\in D, \gamma\in\Gamma\\ \text{paths } (y_k,\cdots,y_0),\  y_0=y, y_k=\gamma y}}\mathrm{Tr}^{\mathbb{C}^\ell}\big[\rho(\gamma)\big],
\end{equation}
which is exactly \eqref{2.23n}.\qed
\end{pro}

\section{Bergman Kernels and Properties}\label{s3n}

In this section, we introduce Bergman kernels and discuss their properties. In \S\,\ref{s3.1}, we define the Bergman kernel associated with a holomorphic line bundle over a compact complex manifold. In \S\,\ref{s3.2}, we recall the off-diagonal exponential decay of Bergman kernels, which will be fundamental for the large spin limit used later. In \S\,\ref{s3.3}, we prove a zero-one law for asymptotic traces. In \S\,\ref{s3.4}, we present concrete examples coming from geometric representations of compact Lie groups and rederive the main results of this section in this setting using lowest weight vectors.

For background in complex geometry, we refer to Huybrechts \cite{MR2093043}. For Bergman kernels, we refer to Ma-Marinescu \cite{MR2339952,MR2393271,MR3368102}. For basic results in the representation theory of compact Lie groups, we refer to Bröcker-tom Dieck \cite{MR0781344}. For the Borel-Weil-Bott theorem, we refer to Bott \cite{MR89473}, Berline-Getzler-Vergne \cite[\S\,8]{MR2273508}, and Etingof \cite[\S\,27]{etingof2024representationsliegroups}. Highest and lowest weight vectors in compact Lie group representations are closely related to coherent states in quantum physics, and we refer to Perelomov \cite{MR858831,MR363209}.

Much of this section serves as preparation for the eigensection results in \S\S\,\ref{s10} and \ref{s11}. Accordingly, for the eigenvalue results in \S\S\,\ref{s5n}, \ref{s6n}, and \ref{det}, it suffices to read only \S\S\,\ref{s3.4.1n} and \ref{s3.4.6n}, and most of the remaining material can be deferred until the discussion of eigensections.

\subsection{Bergman kernels}\label{s3.1}

Let $(N,J)$ be a compact complex manifold with $\dim_{\mathbb{C}}N=n$, and let $(L,h^L)$ be a \emph{positive} holomorphic line bundle over $N$. The first Chern form $c_1(L,h^L)$ induces a Kähler metric $g^{T_\mathbb{R}N}$ on $T_\mathbb{R}N$ and a volume form $dv_N$ on $N$ by
\begin{equation}\label{3.1n}
	g^{T_\mathbb{R}N}(\cdot,\cdot)=c_1(L,h^L)(\cdot,J\cdot),\quad dv_N=\frac{c_1(L,h^L)^n}{n!}.
\end{equation}

For ${p\in{\mathbb{N}}}$, let $L^p=L^{\otimes p}$ be the $p$-th tensor power of $L$, and let $H^{(0,0)}(N,L^p)$ denote the space of holomorphic sections of $L^p$. By the Kodaira vanishing theorem and the Hirzebruch-Riemann-Roch theorem, see \cite[Proposition 3.72, Theorem 4.10]{MR2273508}, we have the following asymptotic dimension formula
\begin{equation}\label{3.2}
	\dim_{\mathbb{C}} H^{(0,0)}(N,L^p)=\mathrm{Vol}(N)p^n+O(p^{n-1}),
\end{equation}
where $\mathrm{Vol}(N)$ is the volume of $N$ with respect to $dv_N$.

Let $\langle\cdot,\cdot\rangle_{H^{(0,0)}(N,L^p)}$ be the $L^2$-inner product on $H^{(0,0)}(N,L^p)$ induced by $(dv_{N},h^{L})$, and let
\begin{equation}\label{3.3n}
P_p\colon L^2(N,L^p)\to H^{(0,0)}(N,L^p)
\end{equation}
be the associated orthogonal projection. Then $P_p$ has a smooth Schwartz kernel $P_p(z,z')\in L_z^p\otimes (L_{z'}^*)^p$ such that for any $s\in C^\infty(N,L^p)$,
\begin{equation}\label{3.3}
	(P_ps)(z)=\int_NP_p(z,z')s(z')dv_N(z').
\end{equation}
This kernel is called \emph{the Bergman kernel}. Let $(s_i)_{i=1}^{\dim_{\mathbb{C}}H^{(0,0)}(N,L^p)}$ be an orthonormal basis of $H^{(0,0)}(N,L^p)$, then
\begin{equation}\label{3.5n}
	P_p(z,z')=\sum_{i=1}^{\dim_{\mathbb{C}}H^{(0,0)}(N,L^p)}s_i(z)\otimes s_i(z')^*,
\end{equation}
where $s_i(z')^*\in (L_{z'}^*)^p$ denotes the metric dual of $s_i(z')$ with respect to $h^{L^p}$.

On the diagonal, we have $P_{p}(z,z)\in \mathrm{End}(L^p)=\mathbb{C}$, and the Bergman kernel admits the following asymptotic expansion \cite[Theorem 4.1.1]{MR2339952}.
\begin{prop}
	There is $C>0$ such that for any $p\in\mathbb{N}$,
\begin{equation}\label{3.6n}
	\lv p^{-n}P_p(z,z)-1\rv_\mathbb{C}\leqslant Cp^{-1}.
\end{equation}
\end{prop}

\subsection{Off-diagonal exponential decay}\label{s3.2}

We have the following exponential decay estimate for Bergman kernels \cite[Theorem 1]{MR3368102}.
\begin{prop}\label{p3.1}
	There exist $C,c>0$ such that for any $p\in\mathbb{N}^*$ and $z,z'\in N$,
	\begin{equation}\label{3.5}
		\Vert P_p(z,z')\Vert_{L_z^p\otimes (L_{z'}^*)^p}\leqslant Cp^{n}e^{-c\sqrt{p}d_N(z,z')},
	\end{equation}
	where $d_N(\cdot,\cdot)$ is the Riemannian geodesic distance on $N$ induced by $g^{T_\mathbb{R}N}$.
\end{prop}

Proposition \ref{p3.1} will play a central role in the analysis of the large spin limit. We now give a first application.

\subsection{A zero-one law for asymptotic traces}\label{s3.3}

Suppose that $\gamma\colon N\to N$ is a holomorphic diffeomorphism that lifts to a holomorphic isometry $\gamma\colon L\to L$. Let $\gamma$ act on $C^\infty(N,L^p)$ by
\begin{equation}\label{3.7n}
	(\gamma s)(z)=\gamma s(\gamma^{-1}z).
\end{equation}
Then this action preserves $H^{(0,0)}(N,L^p)$.
\begin{prop}\label{p3.2n}
The normalized trace satisfies the following dichotomy, as $p\to\infty$,
	\begin{equation}\label{3.7}
		\begin{split}
	\frac{1}{\dim_\mathbb{C}H^{(0,0)}(N,L^p)}&\mathrm{Tr}^{H^{(0,0)}(N,L^p)}[\gamma]\\
	=&\begin{cases}
	c^p,&\text{if}\ \gamma=c\cdot\mathrm{Id}_L \text{ with } \lv c\rv_\mathbb{C}=1,\\
	o(1),&\text{if}\ \gamma\neq c\cdot\mathrm{Id}_L \text{ with } \lv c\rv_\mathbb{C}=1.
\end{cases}
		\end{split}
	\end{equation}
\end{prop}

\begin{pro}
	First, the action of $\gamma$ on $H^{(0,0)}(N,L^p)$ is represented by the kernel $\gamma P_p(\gamma^{-1}z,z')$. Indeed, for $s\in H^{(0,0)}(N,L^p)$, using \eqref{3.3} and \eqref{3.7n}, we compute that
	\begin{equation}\label{3.10n}
		\begin{split}
		(\gamma s)(z)=(\gamma P_ps)(z)=&\gamma(P_ps)(\gamma^{-1}z)\\
		=&\int_N\gamma P_p(\gamma^{-1}z,z')s(z')dv_N(z').
		\end{split}
	\end{equation}
Consequently,
	\begin{equation}\label{3.9}
		\mathrm{Tr}^{H^{(0,0)}(N,L^p)}[\gamma]=\int_N\gamma P_p(\gamma^{-1}z,z)dv_N(z).
	\end{equation}

	By \eqref{3.2}, \eqref{3.5}, and \eqref{3.9}, the limit $p\to\infty$ in \eqref{3.7} is determined by the contribution of the fixed point set of $\gamma$. This is a submanifold of $N$, possibly with several connected components of different dimensions. Since the measure of a proper submanifold is zero unless it coincides with $N$ itself, only the case in which $\gamma$ fixes all of $N$ can contribute. In that case, $\gamma$ induces a holomorphic function on $N$, which must be constant because $N$ is compact. Hence $\gamma=c\cdot\mathrm{Id}_L$ for some $c\in\mathbb{C}$ with $\lv c\rv_\mathbb{C}=1$, as $\gamma$ preserves the Hermitian metric.\qed
\end{pro}

\subsection{Geometric representations of compact Lie groups}\label{s3.4}

\subsubsection{Weyl character formula}\label{s3.4.1n}

Let $G$ be a compact connected Lie group with Lie algebra $\mathfrak{g}$. We fix a maximal torus $T\subseteq G$ with Lie algebra $\mathfrak{t}$. Define the integral lattice $\Lambda\subset \mathfrak{t}$ to be the kernel of the exponential map $\exp\colon\mathfrak{t}\to T$, and the weight lattice $\Lambda^*\subset \mathfrak{t}^*$ by $\Lambda^*=\mathrm{Hom}(\Lambda,2\pi\mathbb{Z})$. We choose a set of positive roots $\Phi^+\subset\Lambda^*$, an open Weyl chamber ${C}^+$ with its closure $\overline{C^+}$, and the associated Weyl group $W$. The Weyl vector $\varrho$ is defined as half the sum of the positive roots
\begin{equation}
\varrho=\frac{1}{2}\sum_{\beta\in \Phi^+}\beta.
\end{equation}

Now we recall the classical Weyl character formula \cite[Theorem \uppercase\expandafter{\romannumeral6}.1.7]{MR0781344}.
\begin{theo}
Irreducible representations $\mathrm{Irr}(G)$ of $G$ are parameterized by
\begin{equation}\label{3.12n}
\mathrm{Irr}(G)=\Lambda^*\cap\overline{C^+}.
\end{equation}
For $\alpha\in \Lambda^*\cap\overline{C^+}$, let $V_\alpha$ denote the unique irreducible representation of $G$ with highest weight $\alpha$. Then
\begin{equation}\label{3.11}
	\dim_\mathbb{C}V_{\alpha}=\prod_{\beta\in\Phi^+}\frac{\lk\alpha+\varrho,\beta\rk}{\lk\varrho,\beta\rk}.
\end{equation}
More generally, for any $H\in \mathfrak{t}$,
\begin{equation}\label{3.10}
	\mathrm{Tr}^{V_{\alpha}}[e^H]=\frac{\sum_{w\in W}\det(w)e^{\sqrt{-1}(w(\alpha+\varrho),H)}}{\prod_{\beta\in\Phi^+}(e^{\sqrt{-1}(\beta,H)/2}-e^{-\sqrt{-1}(\beta,H)/2})}.
\end{equation}
\end{theo}

\subsubsection{Borel-Weil-Bott}\label{s3.4.2n}

We now introduce a geometric realization of $V_\alpha$. 

We have decompositions
\begin{equation}\label{3.15n}
	\begin{split}
&\mathfrak{g}=\mathfrak{t}\oplus \mathfrak{r},\quad\mathfrak{r}=[\mathfrak{t},\mathfrak{g}],\\
&\mathfrak{r}\otimes_\mathbb{R}\mathbb{C}=\mathfrak{r}^+\oplus\mathfrak{r}^-,
	\end{split}
\end{equation}
where $\mathfrak{r}^\pm$ is the sum of the root spaces of $\Phi^\pm$. Hence, we can view
\begin{equation}
\Lambda^*\cap\overline{\mathcal{C}^+}\subset \mathfrak{t}^*\subset \mathfrak{g}^*=\mathfrak{t}^*\oplus \mathfrak{r}^*.
\end{equation}

For any $\alpha\in \Lambda^*\cap\overline{C^+}$, let $\mathcal{O}_\alpha$ be the coadjoint orbit of $\alpha$, and let $G_\alpha\subseteq G$ be the stabilizer subgroup. Then
\begin{equation}
\mathcal{O}_\alpha=G\cdot\alpha\cong G/G_\alpha.
\end{equation}
Let $(L_\alpha,h^{L_\alpha})$ be the canonical prequantum line bundle over $\mathcal{O}_\alpha$. If $\alpha$ is regular, we have
\begin{equation}\label{3.18n}
	\begin{split}
\mathcal{O}_\alpha&\cong G/T,\\
 L_\alpha&\cong G\times_T\mathbb{C}\\
 &=\{(g,a)\in G\times\mathbb{C}\}\big/(g,a)\sim \big(gt,\chi_\alpha(t)a\big) \text{ for any }t\in T,
	\end{split}
\end{equation}
where the character map $\chi_\alpha$ is given by
\begin{equation}\label{3.19n}
\chi_\alpha\colon e^H\in T\mapsto e^{\sqrt{-1}\alpha(H)}\in\mathbb{C}.
\end{equation}

We equip $(\mathcal{O}_\alpha,L_\alpha)$ with a $G$-invariant holomorphic structure such that the holomorphic and antiholomorphic tangent bundles of $\mathcal{O}_\alpha$ satisfy
\begin{equation}\label{3.20n}
T^{(1,0)}\mathcal{O}_\alpha=G\times_T\mathfrak{r}^+,\quad T^{(0,1)}\mathcal{O}_\alpha=G\times_T\mathfrak{r}^-.
\end{equation}
The natural holomorphic $G$-action on $(\mathcal{O}_\alpha,L_\alpha)$ induces a $G$-representation on $H^{(0,0)}(\mathcal{O}_\alpha,L_\alpha)$. The Borel-Weil-Bott theorem \cite[Theorem \rom{4}]{MR89473} yields the following identification.
\begin{theo}
There is an isomorphism of $G$-representations
\begin{equation}\label{3.12}
	H^{(0,0)}(\mathcal{O}_\alpha,L_\alpha)\cong V_{\alpha}.
\end{equation}
Moreover, for any $p\in\mathbb{N}$, since $\mathcal{O}_{p\alpha}\cong \mathcal{O}_\alpha$ and $L_{p\alpha}\cong L_\alpha^p$, it follows that
\begin{equation}\label{3.13}
	H^{(0,0)}(\mathcal{O}_\alpha,L_\alpha^p)\cong V_{p\alpha}.
\end{equation}
\end{theo}

\subsubsection{Lowest weight vector and Bergman kernels}\label{s3.4.3}

For the irreducible representation $V_\alpha$, there exists a unit lowest weight vector $v_{-\alpha}\in V_\alpha$ such that
\begin{equation}\label{3.23n}
	\begin{split}
&tv_{-\alpha}=\chi_{-\alpha}(t)v_{-\alpha} \text{ for any }t\in T,\\
&\mathfrak{r}^-v_{-\alpha}=0,\quad\lV v_{-\alpha}\rV_{V_\alpha}=1.
	\end{split}
\end{equation}

We can describe the isomorphism \eqref{3.12} explicitly as
	\begin{equation}\label{3.24n}
	v\in V_\alpha\mapsto s_v\in H^{(0,0)}(\mathcal{O}_\alpha,L_\alpha),\quad s_v(g)=\langle v, gv_{-\alpha}\rangle_{V_\alpha}.
	\end{equation}
To see this, note first that the first identity in \eqref{3.23n} gives
\begin{equation}
s_v(gt)=\chi_\alpha(t)s_v(g).
\end{equation}
Hence, by \eqref{3.18n} and \eqref{3.19n}, we have $s_v(g)\in C^\infty(\mathcal{O}_\alpha,L_\alpha)$. Then, by \eqref{3.20n}, the second identity in \eqref{3.23n} is equivalent to the holomorphicity condition $\overline{\partial}s_v=0$. Finally, injectivity of the map \eqref{3.24n} follows from the fact that $\mathrm{span}\{gv_{-\alpha}\}_{g\in G}=V_\alpha$.

\begin{lemma}
The map in \eqref{3.24n} is a scalar multiple of an isometry. Indeed,
\begin{equation}\label{3.25n}
	\lV s_v\rV_{L^2(\mathcal{O}_\alpha,L_\alpha)}^2=\frac{\mathrm{Vol}(\mathcal{O}_\alpha)}{\dim_\mathbb{C}V_\alpha}\lV v\rV_{V_\alpha}^2.
\end{equation}
\end{lemma}

\begin{pro}

This follows by
	\begin{equation}\label{3.26n}
		\begin{split}
		\lV s_v\rV_{L^2(\mathcal{O}_\alpha,L_\alpha)}^2&=\int_{\mathcal{O}_\alpha}\langle v, gv_{-\alpha}\rangle_{V_\alpha}\overline{\langle v, gv_{-\alpha}\rangle_{V_\alpha}}dv_{\mathcal{O}_\alpha}(g)\\
		&=\Big\langle v, \int_{\mathcal{O}_\alpha}\langle v, gv_{-\alpha}\rangle_{V_\alpha} gv_{-\alpha}dv_{\mathcal{O}_\alpha}(g)\Big\rangle_{V_\alpha}.
		\end{split}
	\end{equation}
Observe that the map
	\begin{equation}
		v\in V_\alpha\mapsto\int_{\mathcal{O}_\alpha}\langle v, gv_{-\alpha}\rangle_{V_\alpha} gv_{-\alpha}dv_{\mathcal{O}_\alpha}(g)\in V_\alpha
	\end{equation}
	is $G$-equivariant. By Schur's lemma, it must be of the form $C_\alpha\mathrm{Id}_{V_\alpha}$. Taking its trace, we get
	\begin{equation}
		\begin{split}
	C_\alpha\dim_\mathbb{C}V_\alpha&=\int_{\mathcal{O}_\alpha}\mathrm{Tr}^{V_\alpha}\big[\langle\cdot, gv_{-\alpha}\rangle_{V_\alpha} gv_{-\alpha}\big]dv_{\mathcal{O}_\alpha}(g)\\
	&=\int_{\mathcal{O}_\alpha}\lV gv_{-\alpha}\rV^2dv_{\mathcal{O}_\alpha}(g)=\mathrm{Vol}(\mathcal{O}_\alpha)\lV v_{-\alpha}\rV_{V_\alpha}^2=\mathrm{Vol}(\mathcal{O}_\alpha),
		\end{split}
	\end{equation}
which, together with \eqref{3.26n}, implies \eqref{3.25n}.\qed
\end{pro}

We now give an explicit formula for the Bergman kernel \eqref{3.5n}.

\begin{prop}
For $(N,L)=(\mathcal{O}_\alpha,L_\alpha)$, we have the following formula for the Bergman kernel
		\begin{equation}\label{3.29n}
			P_{p}(g,g')=\frac{\dim_\mathbb{C}V_{p\alpha}}{\mathrm{Vol}(\mathcal{O}_\alpha)}\langle g'v_{-\alpha}, gv_{-\alpha}\rangle_{V_\alpha}^p.
\end{equation}
\end{prop}

\begin{pro}
Let $(v_i)$ be an orthonormal basis of $V_\alpha$. By \eqref{3.24n} and \eqref{3.25n}, the sections $(s_{v_i}/\sqrt{C_\alpha})$ form an orthonormal basis of $H^{(0,0)}(\mathcal{O}_\alpha,L_\alpha)$. From \eqref{3.5n}, the Bergman kernel $P_\alpha(g,g')$ is given by
\begin{equation}\label{3.30n}
	\begin{split}
		P_\alpha(g,g')&=\sum_{i}\frac{1}{C_\alpha}s_{v_i}(g)s_{v_i}(g')^*=\sum_{i}\frac{1}{C_\alpha}\langle v_i, gv_{-\alpha}\rangle_{V_\alpha}\overline{\langle v_i, g'v_{-\alpha}\rangle_{V_\alpha}}\\
		&=\frac{\dim_\mathbb{C}V_\alpha}{\mathrm{Vol}(\mathcal{O}_\alpha)}\langle g'v_{-\alpha},gv_{-\alpha}\rangle_{V_\alpha}.
	\end{split}
\end{equation}	
Since $V_{p\alpha}\subset V_\alpha^{\otimes p}$ and $v_{-p\alpha}=v_{-\alpha}^{\otimes p}$, by \eqref{3.30n} we get
\begin{equation}
	P_{p\alpha}(g,g')=\frac{\dim_\mathbb{C}V_{p\alpha}}{\mathrm{Vol}(\mathcal{O}_\alpha)}\langle g'v_{-\alpha}^{\otimes p},gv_{-\alpha}^{\otimes p}\rangle_{V_\alpha^{\otimes p}},
\end{equation}
which is exactly \eqref{3.29n}.\qed
\end{pro}

In particular, we have the following exact version of \eqref{3.6n}.
\begin{prop}
The diagonal Bergman kernel is given by
	\begin{equation}\label{3.33nn}
	P_{p}(g,g)=\frac{\dim_\mathbb{C}V_{p\alpha}}{\mathrm{Vol}(\mathcal{O}_\alpha)}.
\end{equation}
\end{prop}

\subsubsection{Complex projective line}

We now give a concrete example. Take $(N,L)=\big(\mathbb{CP}^1,\mathcal{O}_{\mathbb{CP}^1}(1)\big)$, where $\mathbb{CP}^1$ is the complex projective line and $\mathcal{O}_{\mathbb{CP}^1}(1)$ is the dual tautological line bundle.

Let us denote $\mathcal{O}_{\mathbb{CP}^1}(p)=\mathcal{O}_{\mathbb{CP}^1}(1)^{\otimes p}$. Then the space of holomorphic sections of $\mathcal{O}_{\mathbb{CP}^1}(p)$ identifies with the space of homogeneous polynomials of degree $p$ in two variables, that is,
\begin{equation}
	\begin{split}
		H^{(0,0)}(\mathbb{CP}^1,\mathcal{O}_{\mathbb{CP}^1}(p))=&\mathrm{Sym}^p(\mathbb{C}^2)\\
		=&\{c_0z_0^p+c_1z_0^{p-1}z_1+\cdots+c_p{z_1}^p\mid c_i\in\mathbb{C}\}.
	\end{split}
\end{equation}

Moreover, with respect to the Fubini-Study form $\omega_\mathrm{FS}$ on $\mathbb{CP}^1$, the space $\mathrm{Sym}^p(\mathbb{C}^2)$ has an orthonormal basis
\begin{equation}\label{3.33n}
	\left\{s_{p,i}(z)=\big((p+1)\tbinom{p}{i}\big)^{1/2}z_0^iz_1^{p-i}\mid 0\leqslant i\leqslant p\right\}.
\end{equation}
Therefore, by \eqref{3.5n}, we have
\begin{equation}\label{3.34n}
	\begin{split}
		P_p(z,z')&=\sum_{i=0}^{p}s_{p,i}(z)s_{p,i}(z)^*=\sum_{i=0}^{p}(p+1)\tbinom{p}{i}z_0^iz_1^{p-i}\cdot \overline{z_0'}^i\overline{z_1'}^{p-i}\\
		&=(p+1)\langle z,z'\rangle_{\mathbb{C}^2}^p.
	\end{split}
\end{equation}

We now rederive \eqref{3.34n} using the results of \S\,\ref{s3.4.3}. Consider $\mathrm{SU}_2(\mathbb{C})$, its maximal torus $\mathbb{S}^1$, and the $p$-th character $\chi_p$ of $\mathbb{S}^1$ by
\begin{equation}
	\begin{split}
		\mathrm{SU}_2(\mathbb{C})=\Big\{&g=g_z=\lk\begin{smallmatrix}
			z_0& -\overline{z_1}\\
			z_1& \overline{z_0}
		\end{smallmatrix}\rk\mid z=\lk\begin{smallmatrix}
			z_0\\
			z_1
		\end{smallmatrix}\rk, \lv z_0\rv_{\mathbb{C}}^2+\lv z_1\rv_{\mathbb{C}}^2=1\Big\},\\
		\mathbb{S}^1=\Big\{&g_\theta=\lk\begin{smallmatrix}e^{\sqrt{-1}\theta}& 0\\
			0& e^{-\sqrt{-1}\theta}
		\end{smallmatrix}\rk\mid\theta\in\mathbb{R}\Big\},\\
		\chi_p(&g_\theta)=e^{\sqrt{-1}p\theta}.
	\end{split}
\end{equation}
The Hopf fibration identifies
\begin{equation}
\big(\mathbb{CP}^1,\mathcal{O}_{\mathbb{CP}^1}(p)\big)\cong\big(\mathrm{SU}_2(\mathbb{C})/\mathbb{S}^1,\mathrm{SU}_2(\mathbb{C})\times_{\mathbb{S}^1,\chi_p}\mathbb{C}\big),
\end{equation}
which matches \eqref{3.18n}.

For $g'=g_{z'}\in \mathrm{SU}_2(\mathbb{C})$ (or $g_\theta\in\mathbb{S}^1$), it acts on $s\in\mathrm{Sym}^p(\mathbb{C}^2)$ by
\begin{equation}\label{3.37n}
	\begin{split}
g'(s)(z)&=s((g')^{-1}z)=s\lk\begin{smallmatrix}
	\overline{z_0'}z_0+\overline{z_1'}z_1\\
	-z_1'z_0+z_0'z_1
\end{smallmatrix}\rk,\\
g_\theta(s)(z)&=s\lk\begin{smallmatrix}
	e^{-\sqrt{-1}\theta}z_0\\
	e^{\sqrt{-1}\theta}z_1
\end{smallmatrix}\rk.
	\end{split}
\end{equation}
Hence the unit lowest weight vector is $(p+1)^{1/2}z_0^p$, the monomial of highest $z_0$-degree. By \eqref{3.37n}, in the case $p=1$ and $g'=g_{z'}, g''=g_{z''}$, we have
\begin{equation}\label{3.38n}
\begin{split}
&\big\langle g''(\sqrt{2}z_0), g'(\sqrt{2}z_0)\big\rangle_{\mathrm{Sym}^1(\mathbb{C}^2)}\\
&=\big\langle \sqrt{2}( \overline{z_0''}z_0+\overline{z_1''}z_1),\sqrt{2}(\overline{z_0'}z_0+\overline{z_1'}z_1)\big\rangle_{\mathrm{Sym}^1(\mathbb{C}^2)}=\langle z',z''\rangle_{\mathbb{C}^2},
\end{split}
\end{equation}
Taking together \eqref{3.29n} and \eqref{3.38n}, we recover \eqref{3.34n}.

\subsubsection{Kodaira embedding and off-diagonal exponential decay}

We now explain Proposition \ref{p3.1} in the Borel-Weil-Bott setting.

Let us start with the simplest case $\big(\mathbb{CP}^1,\mathcal{O}_{\mathbb{CP}^1}(p)\big)$. The Fubini-Study distance on $\mathbb{CP}^1$ is
\begin{equation}\label{3.39n}
	d_{\mathrm{FS}}(z,z')=\arccos\lv\langle z,z'\rangle_{\mathbb{C}^2}\rv.
\end{equation}
Combining \eqref{3.34n} and \eqref{3.39n}, we obtain
\begin{equation}
	\begin{split}
		\lV P_p(z,z')\rV_{\mathcal{O}_{\mathbb{CP}^1}(p)_{z'}\otimes\mathcal{O}_{\mathbb{CP}^1}(-p)_{z'}}=&(p+1)\big(\cos d_{\mathrm{FS}}(z,z')\big)^p\\
		\approx&(p+1)\big(1-\tfrac{d_{\mathrm{FS}}(z,z')^2}{2}\big)^p\\
		\approx&(p+1)e^{-pd_{\mathrm{FS}}(z,z')^2/2},
	\end{split}
\end{equation}
which is indeed stronger than \eqref{3.5}. 

We next reduce the general coadjoint orbit case to the projective space case. The \emph{Kodaira map} is defined by
	\begin{equation}\label{3.41n}
		\iota_\alpha\colon g\in\mathcal{O}_\alpha\mapsto [gv_{-\alpha}]\in\mathbb{P}(V_\alpha),
	\end{equation}
	where $\mathbb{P}(V_\alpha)$ denotes the complex projective space of $V_\alpha$. To see that \eqref{3.41n} agrees with the general definition in \cite[Definition 2.2.3]{MR2339952}, note that $\iota_\alpha g$ represents the one-dimensional subspace of $V_\alpha$ such that $\langle \iota_\alpha g, v\rangle_{V_\alpha}=0$
	for every $v$ with $s_v(g)=0$, which is equivalent to \eqref{3.41n} by \eqref{3.24n}.
	
\begin{prop}\label{p3.6n}
The Kodaira map $\iota_\alpha$ is a holomorphic embedding.
\end{prop}	
	
\begin{pro}

Holomorphicity is immediate. First, we show that $\iota_\alpha$ is injective. Suppose that
\begin{equation}\label{3.42n}
gv_{-\alpha}=cv_{-\alpha},\quad\lv c\rv_\mathbb{C}=1.
\end{equation}
We claim that $g\in G_\alpha$, the stabilizer subgroup of $\alpha$. Let us form
\begin{equation}
\begin{split}
\langle \cdot v_{-\alpha},v_{-\alpha}\rangle_{V_\alpha}&\in \mathfrak{g}^*\colon A\in \mathfrak{g}\mapsto \langle Av_{-\alpha},v_{-\alpha}\rangle_{V_\alpha},\\
\langle \cdot gv_{-\alpha},gv_{-\alpha}\rangle_{V_\alpha}&\in \mathfrak{g}^*\colon A\in \mathfrak{g}\mapsto \langle Agv_{-\alpha},gv_{-\alpha}\rangle_{V_\alpha}.
\end{split}
\end{equation}
By \eqref{3.42n}, we clearly have
\begin{equation}
\langle \cdot v_{-\alpha},v_{-\alpha}\rangle_{V_\alpha}=\langle \cdot gv_{-\alpha},gv_{-\alpha}\rangle_{V_\alpha},
\end{equation}
while by definition
\begin{equation}
	\langle \cdot gv_{-\alpha},gv_{-\alpha}\rangle_{V_\alpha}=\mathrm{Ad}^*_g(\langle \cdot v_{-\alpha},v_{-\alpha}\rangle_{V_\alpha}).
\end{equation}
It therefore suffices to show that $\langle \cdot v_{-\alpha},v_{-\alpha}\rangle_{V_\alpha}\in\mathrm{span}\{\alpha\}$. By \eqref{3.23n},
\begin{equation}\label{3.47n}
\langle \cdot v_{-\alpha},v_{-\alpha}\rangle_{V_\alpha}\big|_{\mathfrak{t}}=-\sqrt{-1}\alpha,\quad\langle \cdot v_{-\alpha},v_{-\alpha}\rangle_{V_\alpha}\big|_{\mathfrak{r}^-}=0.
\end{equation}
Moreover, by \eqref{3.15n} we have
\begin{equation}\label{3.48n}
	\langle \cdot v_{-\alpha},v_{-\alpha}\rangle_{V_\alpha}\big|_{\mathfrak{r}^+}=0
\end{equation}
since $\mathfrak{r}^+(v_{-\alpha})$ consists of higher weight vectors, which are orthogonal to $v_{-\alpha}$. Combining \eqref{3.47n} and \eqref{3.48n}, we obtain the claim.

Next, we show that $\iota_\alpha$ is an immersion. Using \eqref{3.20n}, let $(g,A)\in G\times_T\mathfrak{r}=T\mathcal{O}_\alpha$. If $(\iota_\alpha)_*(g,A)=0$, then $gAv_{-\alpha}\in\mathrm{span}\{gv_{-\alpha}\}$, hence $Av_{-\alpha}\in\mathrm{span}\{v_{-\alpha}\}$. Consequently, for any $\varepsilon\in\mathbb{R}$, $e^{\varepsilon A}v_{-\alpha}\in\mathrm{span}\{v_{-\alpha}\}$. By the injectivity part, we have $e^{\varepsilon A}\in G_\alpha$, hence $A\in TG_\alpha\cap T\mathcal{O}_\alpha=0$. This completes the proof.\qed
\end{pro}

We now prove Proposition \ref{p3.1} in the Borel-Weil-Bott setting.

\begin{prop}\label{p3.7n}
For $(N,L)=(\mathcal{O}_\alpha,L_\alpha)$, there exist $C,c>0$ such that for any $p\in\mathbb{N}$,
\begin{equation}\label{3.50n}
\lV P_p(g,g')\rV_{L^{p}_{\alpha,g}\otimes(L^{*}_{\alpha,g'})^p}\leqslant Cp^{\dim_{\mathbb{C}}\mathcal{O}_\alpha}e^{-cpd_{\mathcal{O}_\alpha}(g,g')^2}.
\end{equation}
\end{prop}	
	
\begin{pro}
First, by \eqref{3.2} or \eqref{3.11}, we have for some $C>0$,
\begin{equation}\label{3.51n}
	\dim_\mathbb{C}V_{p\alpha}=\prod_{\beta\in\Phi^+}\frac{\lk p\alpha+\varrho,\beta\rk}{\lk\varrho,\beta\rk}=Cp^{\vert\Phi^+\vert}+O(p^{\vert\Phi^+\vert-1}).
\end{equation}
Since $\vert\Phi^+\vert=\dim_{\mathbb{C}}\mathcal{O}_\alpha$, this matches the factor $Cp^{\dim_{\mathbb{C}}\mathcal{O}_\alpha}$ in \eqref{3.50n}.

Next, the Fubini-Study distance for the projective space $\mathbb{P}(V_\alpha)$ is similar to \eqref{3.39n}, in particular, for $g,g'\in\mathcal{O}_\alpha$,
\begin{equation}
d_{\mathrm{FS}}(\iota_\alpha g,\iota_\alpha g')=\arccos\lv\langle\iota_\alpha g,\iota_\alpha g'\rangle_{V_\alpha}\rv_\mathbb{C}.
\end{equation}
Then by \eqref{3.41n} we have
\begin{equation}\label{3.52n}
\lv\langle gv_{-\alpha},g'v_{-\alpha}\rangle_{V_\alpha}\rv_\mathbb{C}=\cos d_{\mathrm{FS}}(\iota_\alpha g,\iota_\alpha g')\approx 1-\frac{d_{\mathrm{FS}}(\iota_\alpha g,\iota_\alpha g')^2}{2}.
\end{equation}
Combining \eqref{3.29n}, \eqref{3.51n}, and \eqref{3.52n}, we obtain \eqref{3.50n} with $d_{\mathrm{FS}}(\iota_\alpha g,\iota_\alpha g')$ in place of $d_{\mathcal{O}_\alpha}(g,g')$.

Finally, by Proposition \ref{p3.6n}, the map $\iota_\alpha$ is an embedding. Since $\mathcal{O}_\alpha$ is compact, the intrinsic distance $d_{\mathcal{O}_\alpha}(g,g')$ and the pullback Fubini-Study distance $d_{\mathrm{FS}}(\iota_\alpha g,\iota_\alpha g')$ are bi-Lipschitz equivalent. This is analogous to the familiar fact that, for the unit sphere embedded in Euclidean space, the arc length distance and the chord length distance are not equal but are comparable.\qed
\end{pro}

\subsubsection{(Semi)simple property and asymptotic traces zero-one law}\label{s3.4.6n}

Now we explain Proposition \ref{p3.2n} in the Borel-Weil-Bott setting.

By \cite[Theorem \uppercase\expandafter{\romannumeral4}.1.6]{MR0781344}, any two maximal torus in $G$ are conjugate, and every element of $G$ is contained in a maximal tori. By \cite[Theorem \uppercase\expandafter{\romannumeral4}.2.3]{MR0781344} the center $Z(G)$ of $G$ is the intersection of all maximal tori, namely
\begin{equation}\label{3.14}
	Z(G)=\bigcap_{g\in G}gT g^{-1}.
\end{equation}
We call $G$ semisimple if it possesses no nontrivial abelian connected normal subgroup, or equivalently, its center $Z(G)$ is finite, see \cite[Theorem \uppercase\expandafter{\romannumeral5}.3.9]{MR0781344}. Furthermore, $G$ is called simple if it does not have a nontrivial connected normal subgroup.

\begin{prop}\label{p3.8n}
	Suppose that $G$ is compact connected simple and $\alpha$ is any highest weight, or that $G$ is compact connected and $\alpha$ is any regular highest weight. We have
	\begin{equation}\label{3.15}
\frac{1}{\dim_\mathbb{C}V_{p\alpha}}\mathrm{Tr}^{V_{p\alpha}}[g]=\begin{cases}
			\chi_\alpha(g)^p,&\text{if }g\in Z(G),\\
			o(1),&\text{if }g\notin Z(G).
		\end{cases}
	\end{equation}
\end{prop}

\begin{pro}
	From \eqref{3.13}, we can now plug
	\begin{equation}
		(N,L,\gamma)=(\mathcal{O}_\alpha,L_\alpha,g)
	\end{equation}
	into \eqref{3.7} to obtain \eqref{3.15}. It suffices to check that if $g$ fixes $\mathcal{O}_{\alpha}\cong G/G_\alpha$, then $g\in Z(G)$. Clearly $gg'G_\alpha=g'G_\alpha$ if and only if $g\in g'G_\alpha (g')^{-1}$, therefore
	\begin{equation}\label{3.17}
		g\in \bigcap_{g'\in G}g'G_\alpha (g')^{-1}.
	\end{equation}
	
	If $G$ is simple, then $\cap_{g'\in G}g'G_\alpha (g')^{-1}$ is discrete since it is a normal subgroup. We claim that any normal discrete subgroup $G'\subset G$ is contained in $Z(G)$. Indeed, suppose that there exist $g'\in G'$ and $g\in G$ such that $gg'g^{-1}\neq g'$. Choose a path $g(t)$ from $e$ to $g$. Since $G'$ is normal, we have $g(t)g'g(t)^{-1}\in G'$ for all $t$, contradicting the discreteness of $G'$. 
	
	If $\alpha$ is regular, then by \eqref{3.14} and \eqref{3.17} we see that $g\in Z(G)$.
	
	We can also verify \eqref{3.15} using the Weyl character formula. For simplicity, we only consider the case when $\alpha$ is regular. By \eqref{3.10}, we have
	\begin{equation}\label{3.19}
		\mathrm{Tr}^{V_{p\alpha}}[e^H]=\frac{\sum_{w\in W}\det(w)e^{\sqrt{-1}(w(p\alpha+\varrho),H)}}{\prod_{\beta\in\Phi^+}(e^{\sqrt{-1}(\beta,H)/2}-e^{-\sqrt{-1}(\beta,H)/2})}.
	\end{equation}
	The numerator of \eqref{3.19} is uniformly bounded for $p\in\mathbb{N}$, hence, if the denominator is nonzero, the limit \eqref{3.15} vanishes by \eqref{3.51n}. If the denominator is zero, then each differentiation with respect to $H$ in L'Hôpital's rule produces a factor of $p$ in the numerator. Therefore, the limit \eqref{3.15} can be nonzero only if each factor $(e^{\sqrt{-1}(\beta,H)/2}-e^{-\sqrt{-1}(\beta,H)/2})$ vanishes, to match the power $\vert\Phi^+\vert$ of $p$ in \eqref{3.51n}. Now the condition $(\beta,H)\in 2\pi\mathbb{Z}$ for every $\beta\in \Phi^+$ implies that $e^H\in Z(G)$.\qed
\end{pro}
Note that the same argument proves \eqref{3.7} if we replace the Weyl character formula \eqref{3.19} by the Atiyah-Segal-Singer fixed point theorem \cite[Theorems 2.12, 3.3]{MR236951}.

\section{Spin-scale Duality}\label{s4n}

In this section, we describe several limit settings used throughout the paper. In \S\,\ref{s4.1n}, we set up the large scale limit. In \S\,\ref{s4.2n}, we set up the large spin limit. In \S\,\ref{s4.3n}, we introduce the strong uniform limit and explain how the separate limits can be strengthened. In \S\,\ref{s4.4n}, we present an interpretation of spin-scale duality that links the large spin and large scale limits.

In the large scale limit, convergence is understood in the sense of Benjamini-Schramm \cite{MR1873300}. The spin-scale duality is inspired by Anantharaman's discussion \cite[\S\,5]{MR4477342} of Brooks-Le Masson-Lindenstrauss \cite[Theorem 1]{MR3567266}.

\subsection{Benjamini-Schramm convergence}\label{s4.1n}

Let $X$ be a $d$-regular graph of the form \eqref{2.19n}. We write $\mathscr{F}(X)$ for the class of all unitary flat vector bundles over $X$, allowing arbitrary fibre dimension,
\begin{equation}\label{4.1n}
\mathscr{F}(X)=\big\{F\mid F \text{ is a flat unitary vector bundle of the form } \eqref{2.19n}\big\}.
\end{equation}

Consider a series
\begin{equation}\label{4.2n}
(X_m=\Gamma_m\backslash\mathbb{T}_d)_{m\in\mathbb{N}}
\end{equation}
of $d$-regular graphs. In later sections, we shall prove uniform large scale limits of the form
\begin{equation}\label{4.2nn}
\lim_{m\to \infty}\big(\sup_{F_m\in\mathscr{F}(X_m)} \text{a quantity associated with }F_m)=0,
\end{equation}
where the uniformity is with respect to the supremum over all unitary flat vector bundles. This uniformity will be used later to obtain a strong joint limit involving both the large scale parameter and the large spin parameter.

For $x\in X$, define the injectivity radius $\mathrm{inj}_x$ by
\begin{equation}\label{4.4nn}
	\mathrm{inj}_x=\sup\{r\in\mathbb{N}\mid \text{the projection } B_{\mathbb{T}_d}(y,r)\to X\text{ is injective}\}
\end{equation}
where $y\in\mathbb{T}_d$ is any lift of $x$ and $B_{\mathbb{T}_d}(y,r)$ denotes the ball of radius $r$ in $\mathbb{T}_d$ centered at $y$. We say that the series $(X_m)_{m\in\mathbb{N}}$ \eqref{4.2n} converge to the $d$-regular tree $\mathbb{T}_d$ in
the sense of Benjamini-Schramm, denoted by \textbf{[BST]} for short, if for any $k\in\mathbb{N}$,
\begin{equation}\label{4.3nn}
	\lim_{m\to\infty}\frac{\lv\{x\in X_m\mid \mathrm{inj}_x\leqslant k\}\rv}{\lv X_m\rv}=0.
\end{equation}
In particular, \textbf{[BST]} holds if the girth, that is, the length of the shortest cycle, goes to infinity.

\subsection{A family of vector bundles}\label{s4.2n}

Let $X$ be a $d$-regular graph of the form \eqref{2.19n}, with fundamental group $\Gamma$. Let $(N,L)$ be a compact Kähler manifold and a positive line bundle, as in \eqref{3.1n}.

Assume that $\Gamma$ acts on $N$ as holomorphic diffeomorphisms, and that this action lifts to holomorphic isometries of $L$. By \eqref{3.7n}, this action induces, for any $p\in \mathbb{N}$, a unitary representation of $\Gamma$ on $H^{(0,0)}(N,L^p)$. We summarize these actions by
\begin{equation}\label{4.4n}
	\Gamma\curvearrowright \big(N,L,H^{(0,0)}(N,L^p)\big).
\end{equation}
Accordingly, for each $p\in\mathbb{N}$, we obtain an associated unitary flat vector bundle
\begin{equation}\label{4.5n}
F_p=\Gamma\backslash\big(\mathbb{T}_d\times H^{(0,0)}(N,L^p)\big)\in \mathscr{F}(X).
\end{equation}
In later sections, we shall prove pointwise large spin limits of the form
\begin{equation}\label{4.6n}
	\lim_{p\to \infty} (\text{a quantity associated with }F_p)=0,
\end{equation}
where pointwise means that $X$ is fixed.

We say that the action \eqref{4.4n} is free, denoted by \textbf{[FREE]}, if the only element $\gamma\in\Gamma$ acting trivially on $N$ is $\gamma=1$. Equivalently, the induced morphism
\begin{equation}\label{4.8n}
	\Gamma\to \text{the diffeomorphism group of }N
\end{equation}
is injective. Moreover, we say that the action \eqref{4.4n} is dense, denoted by \textbf{[DEN]}, if there exists $z\in N$ such that the orbit
\begin{equation}\label{4.9n}
	\{\gamma z\}_{\gamma\in\Gamma}\subset N
\end{equation}
is dense.

In the Borel-Weil-Bott case $(N,L)=(\mathcal{O}_\alpha,L_\alpha)$ as in \eqref{3.18n}, the action \eqref{4.4n} typically arises from a representation
 \begin{equation}\label{4.10n}
\Gamma\to G,
 \end{equation}
 and by \eqref{3.12} and \eqref{3.13}, we denote
 \begin{equation}\label{4.11nn}
F_{p\alpha}=\Gamma\backslash\big(\mathbb{T}_d\times V_{p\alpha}\big).
 \end{equation}
 In this situation, \textbf{[FREE]} is the injectivity of \eqref{4.10n}, and \textbf{[DEN]} is the denseness of the image of \eqref{4.10n}.

\subsection{Upgrading to the strong uniform limit}\label{s4.3n}

Now we combine the settings of \S\S\,\ref{s4.1n} and \ref{s4.2n}. Let $(X_m)_{m\in\mathbb{N}}$ be a series of $d$-regular graphs as in \eqref{4.2n}, and assume that each fundamental group $\Gamma_m$ acts on $(N,L)$ as in \eqref{4.4n}. This produces a two parameter family of unitary flat vector bundles $F_{m,p}$ as in \eqref{4.5n},
\begin{equation}\label{4.11n}
	\begin{split}
&\Gamma_m\curvearrowright (N,L),\\ &F_{m,p}=\Gamma_m\backslash\big(\mathbb{T}_d\times H^{(0,0)}(N,L^p)\big)\in \mathscr{F}(X_m).
	\end{split}
\end{equation}

In later sections, we shall prove strong uniform limits of the form
\begin{equation}\label{4.13nn}
	\limsup_{\max(m,p)\to\infty}(\text{a quantity associated with }F_{m,p})=0.
\end{equation}

In the Borel-Weil-Bott case $(N,L)=(\mathcal{O}_\alpha,L_\alpha)$ as in \eqref{4.10n} and \eqref{4.11nn}, the construction \eqref{4.11n} becomes
\begin{equation}\label{4.20nn}
	\begin{split}
	&\Gamma_m\to G,\\ &F_{m,p\alpha}=\Gamma_m\backslash\big(\mathbb{T}_d\times V_{p\alpha}\big).
	\end{split}
\end{equation}

The mechanism for upgrading the uniform large scale limit \eqref{4.2nn} and the pointwise large spin limit \eqref{4.6n} to the strong uniform limit \eqref{4.13nn} relies on the following elementary observation.


\begin{lemma}
	For a double series $(a_{m,p})_{(m,p)\in \mathbb{N}^2}$ of real numbers, we have
	\begin{equation}\label{4.13n}
		\limsup_{\max(m,p)\to\infty}a_{m,p}=\sup\Big\{\limsup_{m\to\infty}\big(\sup_{p\in\mathbb{N}}a_{m,p}\big),\sup_{m\in\mathbb{N}}\big(\limsup_{p\to\infty}a_{m,p}\big)\Big\}.
	\end{equation}
\end{lemma}

\begin{pro}
	To see this, for any $k\in\mathbb{N}$, we have
	\begin{equation}
		\begin{split}
			\sup_{\max(m,p)\geqslant k}a_{m,p}&=\sup\Big\{\sup_{m\geqslant k}\big(\sup_{p\in\mathbb{N}}a_{m,p}\big),\sup_{0\leqslant m\leqslant k-1}\big(\sup_{p\geqslant k}a_{m,p}\big)\Big\}\\
			&\geqslant \sup\Big\{\sup_{m\geqslant k}\big(\sup_{p\in\mathbb{N}}a_{m,p}\big),\sup_{0\leqslant m\leqslant k-1}\big(\limsup_{p\to \infty}a_{m,p}\big)\Big\},
		\end{split}
	\end{equation}
and this yields the inequality $\geqslant$ in \eqref{4.13n} by letting $k\to\infty$. To prove the reverse inequality, for any $\varepsilon>0$, we choose $k_0\in\mathbb{N}$ such that
	\begin{equation}\label{1.29}
		\sup_{m\geqslant k_0}\big(\sup_{p\in\mathbb{N}}a_{m,p}\big)\leqslant \limsup_{m\to\infty}\big(\sup_{p\in\mathbb{N}}a_{m,p}\big)+\varepsilon.
	\end{equation}
	Then we choose $k_1\in\mathbb{N}$ such that for any $0\leqslant m\leqslant k_0-1$, 
	\begin{equation}\label{1.30}
		\sup_{p\geqslant k_1}a_{m,p}\leqslant \limsup_{p\to \infty}a_{m,p}+\varepsilon.
	\end{equation}
	From \eqref{1.29} and \eqref{1.30}, we get
	\begin{equation}
		\begin{split}
			&\sup_{\max(m,p)\geqslant\max(k_0,k_1)}a_{m,p}\\
			&\leqslant\sup\Big\{\sup_{m\geqslant k_0}\big(\sup_{p\in\mathbb{N}}a_{m,p}\big),\sup_{0\leqslant m\leqslant k_0-1}\big(\sup_{p\geqslant k_1}a_{m,p}\big)\Big\}\\
			&\leqslant \sup\Big\{\limsup_{m\to\infty}\big(\sup_{p\in\mathbb{N}}a_{m,p}\big)+\varepsilon,\sup_{0\leqslant m\leqslant k_0-1}\big(\limsup_{p\to \infty}a_{m,p}\big)+\varepsilon\Big\},
		\end{split}
	\end{equation}
	and this yields the inequality $\leqslant$ in \eqref{4.13n} by letting $\varepsilon\to0$.\qed
\end{pro}

An immediate consequence of \eqref{4.13n} is the following equivalence.
\begin{lemma}
		For a double series $(a_{m,p})_{(m,p)\in \mathbb{N}^2}$ of complex numbers, the following two conditions are equivalent,
	\begin{equation}\label{4.18n}
	\limsup_{\max(m,p)\to\infty}\lv a_{m,p}\rv_\mathbb{C}=0,\quad\begin{cases}
\lim_{m\to\infty}\sup_{p\in\mathbb{N}}\lv a_{m,p}\rv_{\mathbb{C}}=0,\\
\sup_{m\in\mathbb{N}}\lim_{p\to\infty}\lv a_{m,p}\rv_{\mathbb{C}}=0.
			\end{cases}
	\end{equation}
\end{lemma}

By \eqref{4.18n}, to obtain the strong uniform limit, it suffices to establish the large scale limit $m\to\infty$ uniformly for $p\in\mathbb{N}$, together with the pointwise large spin limit $p\to\infty$ for each fixed $m$. This also clarifies why the uniformity over all unitary flat vector bundles in \eqref{4.2nn} is essential.

\subsection{Spin-scale duality}\label{s4.4n}

Before proving our main results, let us revisit the limits \eqref{4.2nn}, \eqref{4.6n}, and \eqref{4.13nn}. In these statements, the large scale limit $m\to\infty$ and the large spin limit $p\to\infty$ play parallel roles. It is hard not to ask whether these two limits share a common underlying mechanism. We offer a viewpoint that gives a conceptual explanation. Although the discussion below is heuristic, it clarifies the similarities between the two limits that will show up repeatedly in later sections.

Recall first that, in quantum mechanics, high energy corresponds to high frequency, hence to the large eigenvalue limit $\lambda\to\infty$. In microlocal analysis on manifolds, the large eigenvalue limit $\lambda\to\infty$ is equivalent to the semiclassical limit $h\to 0$, where $h$ denotes the Planck constant. Since the wavelength is of order $h$, the limit $h\to 0$ may be viewed as a small wavelength limit, in which the background geometry is fixed while the wavelength tends to zero. Moreover, in the discrete setting, the large scale limit $m\to\infty$ \eqref{4.2nn} may be viewed as keeping the wavelength fixed while the environment size tends to infinity. Formally, the two descriptions are dual.

We now consider the large spin limit $p\to\infty$ \eqref{4.6n}. By \eqref{3.5n}, the Bergman kernel $P_p(z,z')$ contains all the information about $H^{(0,0)}(N,L^p)$. By Propositions \ref{p3.1} and \ref{p3.7n}, we see that $P_p(z,z')$ decays rapidly when $(z,z')$ are separated by more than the scale $O(p^{-1/2})$. Thus, intuitively, we may view $P_p(z,z')$ as a \emph{discretization} of $N$, where each ball of radius $O(p^{-1/2})$ is identified with the center of the ball. Since such a ball has volume $O(p^{-n})$, where $n=\dim_{\mathbb{C}}N=\tfrac{1}{2}\dim_{\mathbb{R}}N$, the manifold $N$ is replaced by $O(p^n)$ points.
This gives the following \emph{spin-scale duality},
\begin{equation}\label{4.23n}
H^{(0,0)}(N,L^p) \approx O(p^n)\text{ points}.
\end{equation}
In other words, the space $H^{(0,0)}(N,L^p)$ is approximately represented by $O(p^n)$ localized states. This is a form of \emph{wave-particle duality} in quantum mechanics. Substituting \eqref{4.23n} into \eqref{4.5n}, we get
\begin{equation}
	F_p\approx \Gamma\backslash \big(\mathbb{T}_d\times \{O(p^n)\ \mathrm{points}\}\big),
\end{equation}
which is a fibre bundle over $X$ with discrete fibre, equivalently, a covering space. Therefore, the large spin limit $p\to\infty$ along the fibre behaves like a large $O(p^n)$-sheeted covering limit over the base graph, that is, like a large scale limit with size $O(\lv X\rv p^n)$.

\section{Alon-Boppana Bound}\label{s5n}

In this section, we prove several versions of the Alon-Boppana bound. In \S\,\ref{s5.1n}, we introduce the nontrivial spectral radius. In \S\,\ref{s5.2n}, we establish a uniform large scale bound. In \S\,\ref{s5.3n}, we prove a large spin bound. In \S\,\ref{s5.4n}, we derive a strong uniform bound. In \S\,\ref{s5.5n}, we interpret these results in the Borel-Weil-Bott setting.

For the large scale Alon-Boppana bound, we refer to Alon \cite[Page 95]{MR875835} and \cite[Theorem 1]{MR1124768}.

\subsection{Nontrivial spectral radius}\label{s5.1n}

Let $(X,F)$ be a $d$-regular graph and a unitary flat vector bundle of the form \eqref{2.19n}.

By \eqref{2.20n} and \eqref{2.21n}, we have
\begin{equation}\label{5.1n}
	\begin{split}
		\lV \Delta^Fu\rV_{L^2(X,F)}^2=&\sum_{y\in D}\BV\sum_{y'\sim y}u(y')\BV_{\mathbb{C}^\ell}^2\\
		\leqslant& d\sum_{y\in D}\sum_{y'\sim y}\bV u(y')\bV_{\mathbb{C}^\ell}^2=d^2\lV u\rV_{L^2(X,F)}^2.
	\end{split}
\end{equation}
We briefly justify the last equality in \eqref{5.1n}. Fix $y''\in D$ and a neighbor $y'''\sim y''$. There is a unique $\gamma\in\Gamma$ such that $\gamma y'''\in D$. Then $\gamma y''\sim \gamma y'''\in D$ and by \eqref{2.19n} and \eqref{2.20n},
\begin{equation}
\bV u(\gamma y'')\bV_{\mathbb{C}^\ell}^2=\bV \rho(\gamma)u(y'')\bV_{\mathbb{C}^\ell}^2=\bV u(y'')\bV_{\mathbb{C}^\ell}^2.
\end{equation}
Thus $\bV u(y'')\bV_{\mathbb{C}^\ell}^2$ is counted once in the summation for the pair $(y,y')=(\gamma y''',\gamma y'')$. We will use similar counting arguments repeatedly.

Consequently,
\begin{equation}
\mathrm{spec}(\Delta^F)\subset [-d,d].
\end{equation}

Define the \emph{nontrivial spectral radius} $r_{F,\mathrm{nt}}$ of $F$ by
\begin{equation}\label{5.4n}
	r_{F,\mathrm{nt}}=\max\big\{\lv\lambda_{F,i}\rv\mid\lambda_{F,i}\neq\pm d\big\}.
\end{equation}
A subtlety is that the eigenvalue multiset may or may not contain $\pm d$, and when it does, these values may occur with multiplicity at most $\dim_{\mathbb{C}}F_x$.

\subsection{Uniform large scale bound}\label{s5.2n}

We now state a uniform large scale Alon-Boppana bound.
\begin{theo}
In the setting of \eqref{4.2n}, assume \textbf{[BST]} \eqref{4.3nn}. Then
\begin{equation}\label{4.6}
\liminf_{m\to\infty}\big(\inf_{F_m\in \mathscr{F}(X_m)}r_{F_m,\mathrm{nt}}\big)\geqslant 2\sqrt{d-1},
\end{equation}
where, for each $m\in\mathbb{N}$, the infimum is taken over all unitary flat vector bundles $F_m$ over $X_m$, see \eqref{4.1n} and \eqref{4.2nn}.
\end{theo}

\begin{pro}
First, let $(X,F)$ be a $d$-regular graph and a unitary flat vector bundle of the form \eqref{2.19n}. By \eqref{5.4n}, for any $k\in\mathbb{N}$, the first summation in \eqref{2.23n} satisfies
\begin{equation}\label{4.7}
	\sum_{i=1}^{\dim_\mathbb{C}F}\lambda_{F,i}^{2k}\leqslant 2\dim_{\mathbb{C}}F_x\cdot d^{2k}+(\dim_{\mathbb{C}}F-2\dim_{\mathbb{C}}F_x)r_{F,\mathrm{nt}}^{2k}.
\end{equation}

Then we split the second summation in \eqref{2.23n} into two parts according to the injectivity radius,
\begin{equation}\label{4.8}
	\begin{split}
	\sum_{\substack{x\in X,\\ \gamma\in \Omega(X,x),\lv\gamma\rv=2k}}\mathrm{Tr}^{F_{x}}\big[\rho(\gamma)\big]=&\sum_{\substack{x\in X, \mathrm{inj}_x\geqslant k,\\
		\gamma\in \Omega(X,x),\lv\gamma\rv=2k}}\mathrm{Tr}^{F_{x}}\big[\rho(\gamma)\big]\\
	&+\sum_{\substack{x\in X, \mathrm{inj}_x<k-1,\\
		\gamma\in \Omega(X,x),\lv\gamma\rv=2k}}\mathrm{Tr}^{F_{x}}\big[\rho(\gamma)\big].
	\end{split}
\end{equation}
Let us fix a root point $o\in\mathbb{T}_d$. If $\mathrm{inj}_x\geqslant k$, then a loop $\gamma$ at $x$ of length $2k$ lifts to a loop at $o$ of length $2k$ in $\mathbb{T}_d$. Therefore, the corresponding holonomy is trivial, $\rho(\gamma)=\mathrm{Id}$, and
\begin{equation}\label{4.9}
	\sum_{\gamma\in \Omega(X,x),\lv\gamma\rv=2k}\mathrm{Tr}^{F_{x}}\big[\rho(\gamma)\big]=\bv\{\gamma\in\Omega(\mathbb{T}_d,o)\mid\lv\gamma\rv=2k\}\bv\cdot \dim_{\mathbb{C}}F_x.
\end{equation}
If $\mathrm{inj}_x\leqslant k-1$, then, since $X$ is $d$-regular and $\rho(\gamma)$ is unitary, we get
\begin{equation}
	\lv\{\gamma\in \Omega(X,x)\mid\lv\gamma\rv=2k\}\rv\leqslant d^{2k},\quad\lv\mathrm{Tr}^{F_x}\big[\rho(\gamma)\big]\rv\leqslant \dim_\mathbb{C}F_x,
\end{equation}
and it follows that
\begin{equation}\label{4.11}
	\Bv\sum_{\gamma\in \Omega(X,x),\lv\gamma\rv=2k}\mathrm{Tr}^{F_{x}}\big[\rho(\gamma)\big]\Bv\leqslant d^{2k}\dim_{\mathbb{C}}F_x.
\end{equation}

Combining \eqref{2.23n}, \eqref{4.7}, \eqref{4.8}, \eqref{4.9}, and \eqref{4.11}, we get
\begin{equation}\label{4.13}
	\begin{split}
	r_{F,\mathrm{nt}}^{2k}\geqslant&\frac{(\dim_{\mathbb{C}}F-2\dim_{\mathbb{C}}F_x)r_{F,\mathrm{nt}}^{2k}}{\dim_\mathbb{C}F}\\
	\geqslant&\frac{\bv\{x\in X\mid \mathrm{inj}_x\geqslant k\}\bv}{\lv X\rv}\cdot\bv\{\gamma\in\Omega(\mathbb{T}_d,o)\mid\lv\gamma\rv=2k\}\bv\\
		&-\frac{\lv\{x\in X\mid \mathrm{inj}_x\leqslant k-1\}\rv}{\lv X\rv}d^{2k}-\frac{2d^{2k}}{\lv X\rv},
	\end{split}
\end{equation}
where we recall that $\dim_{\mathbb{C}}F=\dim_{\mathbb{C}}F_x\cdot\lv X\rv$. Since the right hand side of \eqref{4.13} depends only on $X$, we can replace $r_{F,\mathrm{nt}}^{2k}$ with $\inf_{F\in \mathscr{F}(X)}r_{F,\mathrm{nt}}^{2k}$. Applying this to each $X_m$, and using \textbf{[BST]} \eqref{4.3nn}, we obtain
\begin{equation}\label{4.14}
	\begin{split}
		&\liminf_{m\to\infty}\inf_{F_m\in \mathscr{F}(X_m)}r_{F_m,\mathrm{nt}}^{2k}\\
		&\geqslant \bv\{\gamma\in \Omega(\mathbb{T}_d,o)\mid\lv\gamma\rv=2k\}\bv\\
		&\geqslant \bv\{\gamma\in \Omega(\mathbb{T}_d,o)\mid\lv\gamma\rv=2k,\ \text{ending at\ } o \text{\ for the first time}\}\bv\\
		&=\frac{1}{k}\binom{2(k-1)}{k-1}d(d-1)^{k-1},
	\end{split}
\end{equation}
which implies \eqref{4.6} by taking the $2k$-th root and letting $k\to\infty$.\qed
\end{pro}

\subsection{Large spin bound}\label{s5.3n}

We now state a large spin bound.
\begin{theo}
In the setting of \eqref{4.5n}, assume \textbf{[FREE]} \eqref{4.8n} and \textbf{[DEN]} \eqref{4.9n}. Then
\begin{equation}\label{4.15}
	\liminf_{p\to\infty}r_{F_p,\mathrm{nt}}\geqslant 2\sqrt{d-1}.
\end{equation}
\end{theo}

\begin{pro}
We first claim that for $p\geqslant 1$,
\begin{equation}\label{5.15n}
\pm d\notin\mathrm{spec}(\Delta^{F_p}).
\end{equation}
Indeed, suppose $\Delta^{F_p}u=\pm du$. By \eqref{2.20n}, \eqref{4.5n}, and a maximal principle argument, for $u(o)\in H^{(0,0)}(N,L^p)$ satisfying
\begin{equation}
\lV u(o)\rV_{L^2(N,L^p)}=\max_{y\in D}\lV u(y)\rV_{L^2(N,L^p)},
\end{equation}
the eigenvalue equation forces
\begin{equation}\label{5.16n}
\begin{split}
	&u(y)=(\pm1)^{d_{\mathbb{T}_d}(y,o)}u(o),\quad\text{for any }y\in\mathbb{T}_d,\\
	&\gamma (u(o))=u(\gamma o)=(\pm1)^{d_{\mathbb{T}_d}(\gamma o,o)}u(o),\quad\text{for any }\gamma\in\Gamma,
\end{split}
\end{equation}
where we view $o\in \mathbb{T}_d$ as a root and $d_{\mathbb{T}_d}(y,o)$ is the tree distance. Combining \eqref{3.7n} and \eqref{5.16n}, we take pointwise norms on $N$ to get
\begin{equation}
	\begin{split}
\bV\big(\gamma u(o)\big)(z)\bV_{L_z^p}=&\bV\gamma\cdot \big(u(o)(\gamma^{-1}z)\big)\bV_{L_z^p}\\
=&\bV(u(o)(\gamma^{-1}z)\bV_{L_{\gamma^{-1}z}^p}=\bV u(o)(z)\bV_{L_z^p}.
	\end{split}
\end{equation}
By \textbf{[DEN]} \eqref{4.9n}, the function $z\mapsto\Vert u(o)(z)\Vert_{L_z^p}$ is constant. If $u(o)\neq0$, then $u(o)$ is a nowhere vanishing holomorphic section of $L^p$, forcing $L^p$ to be holomorphically trivial, contradicting positivity for $p\geqslant 1$.

By \eqref{5.4n} and \eqref{5.15n}, the first sum in \eqref{2.23n} satisfies
\begin{equation}\label{4.16}
	\sum_{i=1}^{\dim_\mathbb{C}F_p}\lambda_{p,i}^{2k}\leqslant \dim_{\mathbb{C}}F_p\cdot r_{F_p,\mathrm{nt}}^{2k}.
\end{equation}

We split the second sum in \eqref{2.23n} into two parts according to the homotopy class $[\gamma]\in \Gamma$,
\begin{equation}\label{4.17}
	\begin{split}
	\sum_{\substack{x\in X,\gamma\in \Omega(X,x),\\ \lv\gamma\rv=2k}}\mathrm{Tr}^{F_{p,x}}\big[\rho(\gamma)\big]=&\sum_{\substack{x\in X,\gamma\in \Omega(X,x),\\
		\lv\gamma\rv=2k,[\gamma]=1\in\Gamma}}\mathrm{Tr}^{F_{p,x}}\big[\rho(\gamma)\big]\\
	&+\sum_{\substack{x\in X,\gamma\in \Omega(X,x), \\
		\lv\gamma\rv=2k, [\gamma]\neq1\in\Gamma }}\mathrm{Tr}^{F_{p,x}}\big[\rho(\gamma)\big].
	\end{split}
\end{equation}
If $[\gamma]=1\in\Gamma$, that is, the loop $\gamma$ is homotopy to a trivial loop, then $\gamma$ corresponds to a loop from $o$ to $o$ in $\mathbb{T}_d$, therefore
\begin{equation}\label{4.18}
\sum_{\substack{\gamma\in \Omega(X,x),\\
		\lv\gamma\rv=2k,[\gamma]=1\in\Gamma}}\mathrm{Tr}^{F_{p,x}}\big[\rho(\gamma)\big]=\bv\{\gamma\in\Omega(\mathbb{T}_d,o)\mid\lv\gamma\rv=2k\}\bv\cdot \dim_{\mathbb{C}}F_{p,x}.
\end{equation}
If $[\gamma]\neq1\in\Gamma$, then by \textbf{[FREE]} \eqref{4.8n}, the action of $[\gamma]$ on $L$ is not scalar multiplication, hence \eqref{3.7} gives
\begin{equation}\label{4.19}
\lim_{p\to \infty}\bbv\frac{1}{\dim_{\mathbb{C}}F_{p,x}}\mathrm{Tr}^{F_{p,x}}\big[\rho(\gamma)\big]\bbv=0.
\end{equation}

Combining \eqref{2.23n}, \eqref{4.16}, \eqref{4.17}, \eqref{4.18}, and \eqref{4.19}, we obtain
\begin{equation}\label{5.22}
\liminf_{p\to\infty}r_{F_p,\mathrm{nt}}^{2k}\geqslant\bv\{\gamma\in \Omega(\mathbb{T}_d,o),\lv\gamma\rv=2k\}\bv,
\end{equation}
from which \eqref{4.15} follows exactly as in \eqref{4.14}.\qed
\end{pro}

Note that, without \eqref{5.15n}, one could only use \eqref{4.7} instead of \eqref{4.14}. This would introduce an extra term $-2d^{2k}/\lv X\rv$
in \eqref{5.22}, which would in general make the lower bound trivial.

\subsection{Strong uniform bound}\label{s5.4n}

We now combine the large scale and the large spin bounds to get a strong uniform bound.

\begin{theo}\label{t5.3}
In the setting of \eqref{4.11n}, assume \textbf{[BST]} \eqref{4.3nn}, and assume \textbf{[FREE]} \eqref{4.8n} and \textbf{[DEN]} \eqref{4.9n} for every $\Gamma_m$. Then
	\begin{equation}\label{4.21}
		\liminf_{\max(m,p)\to\infty}r_{F_{m,p},\mathrm{nt}}\geqslant 2\sqrt{d-1}.
	\end{equation}
\end{theo}

\begin{pro}
Applying \eqref{4.6} to $(F_{m,p})_{p\in\mathbb{N}}\subset \mathscr{F}(X_m)$, we get
\begin{equation}\label{4.22}
	\liminf_{m\to\infty}\big(\inf_{p\in\mathbb{N}}r_{F_{m,p},\mathrm{nt}}\big)\geqslant 2\sqrt{d-1}.
\end{equation}
Applying \eqref{4.15} to $(F_{m,p})_{p\in\mathbb{N}}$, we obtain for every $m\in\mathbb{N}$,
	\begin{equation}\label{4.23}
\liminf_{p\to\infty} r_{F_{m,p},\mathrm{nt}}\geqslant 2\sqrt{d-1}.
	\end{equation}
Taking $a_{m,p}=-r_{F_{m,p},\mathrm{nt}}$ in \eqref{4.13n}, we get
\begin{equation}
	\begin{split}
	&\liminf_{\max(m,p)\to\infty}r_{F_{m,p},\mathrm{nt}}\\
	&=\inf\Big\{\liminf_{m\to\infty}\big(\inf_{p\in\mathbb{N}}r_{F_{m,p},\mathrm{nt}}\big),\inf_{m\in\mathbb{N}}\big(\liminf_{p\to\infty}r_{F_{m,p},\mathrm{nt}}\big)\Big\},
	\end{split}
\end{equation}
which, combined with \eqref{4.22} and \eqref{4.23}, yields \eqref{4.21}.\qed
\end{pro}

\subsection{Alon-Boppana bound in the Borel-Weil-Bott case}\label{s5.5n}

Now we state the main results of this section in the Borel-Weil-Bott case.

\begin{theo}
In the setting of \eqref{4.20nn}, assume \textbf{[BST]} \eqref{4.3nn}, and suppose that $G$ is compact connected. Then
	\begin{equation}
		\liminf_{m\to\infty}\big(\inf_{\alpha\in \Lambda^*\cap\overline{C^+}}r_{F_{m,\alpha},\mathrm{nt}}\big)\geqslant 2\sqrt{d-1},
	\end{equation}
	where, for each $m\in\mathbb{N}$, the infimum is taken over all irreducible representations of $G$, see \eqref{3.12n}.
\end{theo}

\begin{theo}
In the setting of \eqref{4.11nn}, assume \textbf{[FREE]} and \textbf{[DEN]} \eqref{4.10n}. Suppose, in addition, that either $G$ is compact connected simple and $\alpha$ is an arbitrary highest weight, or that $G$ is compact connected and $\alpha$ is regular. Then
	\begin{equation}
		\liminf_{p\to\infty}r_{F_{p\alpha},\mathrm{nt}}\geqslant 2\sqrt{d-1}.
	\end{equation}
\end{theo}

\begin{theo}\label{t5.6}
	In the setting of \eqref{4.20nn}, assume \textbf{[BST]} \eqref{4.3nn}, and assume \textbf{[FREE]} and \textbf{[DEN]} \eqref{4.10n} for every $\Gamma_m$. Suppose, in addition, that either $G$ is compact connected simple and $\alpha$ is an arbitrary highest weight, or that $G$ is compact connected and $\alpha$ is regular. Then
	\begin{equation}
	\liminf_{\max(m,p)\to\infty}r_{F_{m,p\alpha},\mathrm{nt}}\geqslant 2\sqrt{d-1}.
	\end{equation}

\end{theo}

\section{Kesten-McKay Law}\label{s6n}

In this section, we prove several versions of the Kesten-McKay law. In \S\,\ref{s6.1n}, we introduce spectral distributions. In \S\,\ref{s6.2n}, we establish a uniform large scale law. In \S\,\ref{s6.3n}, we prove a large spin law. In \S\,\ref{s6.4n}, we derive a strong uniform law. In \S\,\ref{s6.5n}, we interpret these results in the Borel-Weil-Bott setting.

For the large scale Kesten-McKay law, we refer to Kesten \cite[Theorem 3]{MR109367} and McKay \cite[Theorem 1.1]{MR629617}. For the large spin Kestin-McKay law, we refer to Gamburd-Jakobson-Sarnak \cite[Proposition 4.1]{MR1677685} and Brooks-Le Masson-Lindenstrauss \cite[Corollary 2]{MR3567266}.

\subsection{Spectral distributions}\label{s6.1n}

Let $(X,F)$ be a $d$-regular graph and a unitary flat vector bundle of the form \eqref{2.19n}. Let $d\mu_{F}(\lambda)$ denote the normalized counting measure of the eigenvalues of $\Delta^F$,
\begin{equation}\label{4.24}
	d\mu_{F}(\lambda)=\frac{1}{\dim_\mathbb{C}F}\sum_{i=1}^{\dim_\mathbb{C}F}\delta_{\lambda_{F,i}}(\lambda).
\end{equation}

We note that the Kesten-McKay distribution defined in \eqref{1.21n} can be characterized by the following moment identity
\begin{equation}
	\int_\mathbb{R}\lambda^kd\mu_{\mathrm{KM}}(\lambda)=\bv\{\gamma\in \Omega(\mathbb{T}_d,o)\mid\lv\gamma\rv=k\}\bv
\end{equation}
for any $k\in\mathbb{N}$.


\subsection{Uniform large scale law}\label{s6.2n}

We now state the uniform large scale eigenvalue distribution.
\begin{theo}
In the setting of \eqref{4.2n}, assume \textbf{[BST]} \eqref{4.3nn}. Then for any interval $I\subseteq\mathbb{R}$,
	\begin{equation}\label{4.25}
		\limsup_{m\to\infty}\sup_{F_m\in \mathscr{F}(X_m)}\lv\frac{1}{\dim_\mathbb{C}F_m}\lv\{i\mid\lambda_{F_m,i}\in I\}\rv-\int_{I}d\mu_{\mathrm{KM}}(\lambda)\rv=0,
	\end{equation}
where, for each $m\in\mathbb{N}$, the supremum is taken over all unitary flat vector bundles $F_m$ over $X_m$, see \eqref{4.1n} and \eqref{4.2nn}.
\end{theo}

\begin{pro}
	By \eqref{2.23n} and \eqref{4.24}, we obtain	for any $k\in\mathbb{N}$,
\begin{equation}\label{4.26}
	\int \lambda^kd\mu_{F}=\frac{1}{\dim_\mathbb{C}F}\sum_{x\in X,\gamma\in \Omega(X,x),\lv\gamma\rv=k}\mathrm{Tr}^{F_{x}}\big[\rho(\gamma)\big],
\end{equation}
which, together with \eqref{4.8}, \eqref{4.9}, and \eqref{4.11}, implies
\begin{equation}\label{4.27}
	\begin{split}
		\bbv\int\lambda^kd\mu_{F}(\lambda)-\frac{\bv\{x\in X\mid \mathrm{inj}_x\geqslant k\}\bv}{\lv X\rv}\cdot\bv\{\gamma\in \Omega(\mathbb{T}_d,o)\mid\lv\gamma\rv=k\}\bv\bbv\\
		\leqslant \frac{\lv\{x\in X\mid \mathrm{inj}_x\leqslant k-1\}\rv}{\lv X\rv}d^k.
	\end{split}
\end{equation}
Since the right hand side of \eqref{4.27} depends only on $X$, we can replace the left hand side with $\sup_{F\in \mathscr{F}(X)}$. Applying this to each $X_m$, and using \textbf{[BST]} \eqref{4.3nn}, we get for any $k\in\mathbb{N}$,
\begin{equation}\label{4.28}
\limsup_{m\to\infty}\sup_{F_m\in \mathscr{F}(X_m)}\bbv\int\lambda^kd\mu_{F_m}(\lambda)-\bv\{\gamma\in \Omega(\mathbb{T}_d,o)\mid\lv\gamma\rv=k\}\bv\bbv=0.
\end{equation}
Using an approximation argument, we see that \eqref{4.28} implies \eqref{4.25}.\qed
\end{pro}

\subsection{Large spin law}\label{s6.3n}

We now state the large spin eigenvalue distribution.
\begin{theo}
	In the setting of \eqref{4.5n}, assume \textbf{[FREE]} \eqref{4.8n}. Then for any interval $I\subseteq\mathbb{R}$,
\begin{equation}\label{4.29}
	\limsup_{p\to\infty}\lv\frac{1}{\dim_\mathbb{C}F_p}\lv\{i\mid\lambda_{F_p,i}\in I\}\rv-\int_{I}d\mu_{\mathrm{KM}}(\lambda)\rv=0.
\end{equation}
\end{theo}

\begin{pro}
By \textbf{[FREE]} \eqref{4.8n}, \eqref{4.17}, \eqref{4.18}, \eqref{4.19}, and \eqref{4.26}, we get for any $k\in\mathbb{N}$,
\begin{equation}
	\limsup_{p\to \infty}\bbv\int \lambda^kd\mu_{F_p}(\lambda)-\bv\{\gamma\in \Omega(\mathbb{T}_d,o)\mid\lv\gamma\rv=k\}\bv\bbv=0,
\end{equation}	
	and \eqref{4.29} follows by an approximation argument.\qed
\end{pro}

\subsection{Strong uniform law}\label{s6.4n}

We now state the strong uniform eigenvalue distribution.
\begin{theo}\label{t6.3}
In the setting of \eqref{4.11n}, assume \textbf{[BST]} \eqref{4.3nn}, and assume \textbf{[FREE]} \eqref{4.8n} for every $\Gamma_m$. Then for any interval $I\subseteq \mathbb{R}$,
	\begin{equation}\label{4.31}
		\limsup_{\max(m,p)\to\infty}\lv\frac{1}{\dim_\mathbb{C}F_{m,p}}\lv\{i\mid\lambda_{F_{m,p},i}\in I\}\rv-\int_{I}d\mu_{\mathrm{KM}}(\lambda)\rv=0.
	\end{equation}
\end{theo}

\begin{pro}
By applying \eqref{4.25} to $(F_{m,p})_{p\in\mathbb{N}}\subset \mathscr{F}(X_m)$ and replacing $(F_p)_{p\in\mathbb{N}}$ with $(F_{m,p})_{p\in\mathbb{N}}$ in \eqref{4.29}, we obtain \eqref{4.31} from \eqref{4.18n}.\qed
\end{pro}

\subsection{Kesten-McKay law in the Borel-Weil-Bott case}\label{s6.5n}

We now state the main results of this section in the Borel-Weil-Bott setting.

\begin{theo}
	In the setting of \eqref{4.20nn}, assume \textbf{[BST]} \eqref{4.3nn}, and suppose that $G$ is compact connected. Then for any interval $I\subseteq\mathbb{R}$,
	\begin{equation}
		\limsup_{m\to\infty}\sup_{\alpha\in \Lambda^*\cap\overline{C^+}}\lv\frac{1}{\dim_\mathbb{C}F_{m,\alpha}}\lv\{i\mid\lambda_{F_{m,\alpha},i}\in I\}\rv-\int_{I}d\mu_{\mathrm{KM}}(\lambda)\rv=0,
	\end{equation}
	where, for each $m\in\mathbb{N}$, the supremum is taken over all irreducible representations of $G$, see \eqref{3.12n}.
\end{theo}

\begin{theo}
	In the setting of \eqref{4.11nn}, assume \textbf{[FREE]} \eqref{4.10n}. Suppose, in addition, that either $G$ is compact connected simple and $\alpha$ is an arbitrary highest weight, or that $G$ is compact connected and $\alpha$ is regular. Then for any interval $I\subseteq\mathbb{R}$,
\begin{equation}
	\limsup_{p\to\infty}\lv\frac{1}{\dim_\mathbb{C}F_{p\alpha}}\lv\{i\mid\lambda_{F_{p\alpha},i}\in I\}\rv-\int_{I}d\mu_{\mathrm{KM}}(\lambda)\rv=0.
\end{equation}
\end{theo}

\begin{theo}\label{t6.6}
	In the setting of \eqref{4.20nn}, assume \textbf{[BST]} \eqref{4.3nn},  and assume \textbf{[FREE]} \eqref{4.10n} for every $\Gamma_m$. Suppose, in addition, that either $G$ is compact connected simple and $\alpha$ is an arbitrary highest weight, or that $G$ is compact connected and $\alpha$ is regular. Then for any interval $I\subseteq\mathbb{R}$,
	\begin{equation}
	\limsup_{\max(m,p)\to\infty}\lv\frac{1}{\dim_\mathbb{C}F_{m,p\alpha}}\lv\{i\mid\lambda_{F_{m,p\alpha},i}\in I\}\rv-\int_{I}d\mu_{\mathrm{KM}}(\lambda)\rv=0.
\end{equation}
\end{theo}

\section{Asymptotic Log Determinants of Laplacians}\label{det}

In this section, we study asymptotic log determinants of Laplacians for graph vector bundles. In \S\,\ref{ns7.1nn}, we introduce the normalized log determinant of the Laplacian. In \S\,\ref{ns7.1n}, we introduce a spectral gap condition. In \S\,\ref{ns7.2n}, we obtain the uniform large scale asymptotic log determinant. In \S\,\ref{ns7.3n}, we obtain the large spin asymptotic log determinant. In \S\,\ref{ns7.4n}, we obtain the strong uniform asymptotic log determinant. In \S\,\ref{ns7.5n}, we interpret these results in the Borel-Weil-Bott setting.

For the large scale asymptotic determinant, we refer to Lyons \cite{MR2160416}.

\subsection{Log determinant}\label{ns7.1nn}

The operator $(d-\Delta^F)$ is nonnegative. We define the normalized log determinant of $(d-\Delta^F)$ by
\begin{equation}\label{n7.1}
	\begin{split}
\frac{1}{\dim_{\mathbb{C}}F}\ln\det(d-\Delta^F)&=\sum_{\lambda_{F,i}\neq d}\frac{1}{\dim_{\mathbb{C}}F}\ln(d-\lambda_{F,i})\\
&=\int_{[-d,d)}\ln(d-\lambda)d\mu_F(\lambda).
	\end{split}
\end{equation}

The generating function of the moments of the Kesten-Mckay distribution \eqref{1.21n} is given by
\begin{equation}
	f_{\mathrm{KM}}(z)=\sum_{k=0}^\infty\Big(\int_\mathbb{R}\lambda^kd\mu_{\mathrm{KM}}(\lambda)\Big)z^k=\frac{2(d-1)}{d-2+d\sqrt{1-4(d-1)z^2}}.
\end{equation}
Then we compute that
\begin{equation}
	\begin{split}
&\int_\mathbb{R}\ln(d-\lambda)d\mu_{\text{KM}}(\lambda)\\
&=\ln(d)-\int_0^{1/d}(f_{\text{KM}}(z)-1)\frac{dz}{z}\\
		&=(d-1)\ln(d-1)-\frac{d-2}{2}\ln(d-2)-\frac{d-2}{2}\ln d.
	\end{split}
\end{equation}

\subsection{Spectral gap condition}\label{ns7.1n}

This integral \eqref{n7.1} has singularity near $d$, which motivates us to introduce the following condition.

We define the \emph{nontrivial top eigenvalue} of $\Delta^F$, an one sided analogue of \eqref{5.4n}, by
\begin{equation}\label{7.90n}
	\lambda_{F,\mathrm{nt}}=\max\{\lambda_{F,i}\mid \lambda_{F,i}\neq d\}.
\end{equation}
Recall the class of unitary flat vector bundles $\mathscr{F}(X)$ defined in \eqref{4.1n}. For $\varepsilon>0$, we define a subclass $\mathscr{F}(X,\varepsilon)$ by
\begin{equation}\label{n7.3}
	\begin{split}
		\mathscr{F}(X,\varepsilon)&=\big\{F\in\mathscr{F}(X)\mid\lambda_{F,\mathrm{nt}}\leqslant d-\varepsilon\big\}\\
		&=\big\{F\in\mathscr{F}(X)\mid\mathrm{spec}(\Delta^F)\subset [-d,d-\varepsilon]\cup\{d\}\big\}.
	\end{split}
\end{equation}
We say that $F$ has an $\varepsilon$-spectral gap, denoted by \textbf{[GAP]}, if $F\in\mathscr{F}(X,\varepsilon)$.

\subsection{Uniform large log determinant}\label{ns7.2n}

We now state the uniform large scale asymptotic log determinant.
\begin{theo}
	In the setting of \eqref{4.2n}, assume \textbf{[BST]} \eqref{4.3nn} and \textbf{[GAP]} \eqref{n7.3}. Then 
	\begin{equation}\label{n7.6}
		\begin{split}
		\limsup_{m\to\infty}\sup_{F_m\in \mathscr{F}(X_m,\varepsilon)}\Bv&\frac{1}{\dim_{\mathbb{C}}F_m}\ln\det(d-\Delta^{F_m})\\
		&-\int_\mathbb{R}\ln(d-\lambda)d\mu_{\mathrm{KM}}(\lambda)\Bv=0,
		\end{split}
	\end{equation}
	where, for each $m\in\mathbb{N}$, the supremum is taken over all unitary flat vector bundles $F_m$ with \textbf{[GAP]}.
\end{theo}

\begin{pro}
Under the spectral gap assumption $F_m\in \mathscr{F}(X_m,\varepsilon)$, the nontrivial spectrum is contained in $[-d,d-\varepsilon]$. Hence the function $\ln(d-\lambda)$ is continuous on the relevant spectral support. We can apply \eqref{4.25} to \eqref{n7.1} to get \eqref{n7.6}.\qed
\end{pro}

\subsection{Large spin log determinant}\label{ns7.3n}

We now state the large spin asymptotic log determinant.
\begin{theo}
	In the setting of \eqref{4.5n}, assume \textbf{[FREE]} \eqref{4.8n}, and assume \textbf{[GAP]} \eqref{n7.3} for every $F_p$. Then
	\begin{equation}\label{n7.7}
	\limsup_{p\to\infty}\lv\frac{1}{\dim_\mathbb{C}F_p}\ln\det(d-\Delta^{F_p})-\int_\mathbb{R}\ln(d-\lambda)d\mu_{\mathrm{KM}}(\lambda)\rv=0.
	\end{equation}
\end{theo}

\begin{pro}
	Under \textbf{[GAP]} \eqref{n7.3} for every $F_p$, we can apply \eqref{4.29} to \eqref{n7.1} to get \eqref{n7.7}.\qed
\end{pro}

\subsection{Strong uniform asymptotic log determinant}\label{ns7.4n}

We now state the strong uniform asymptotic log determinant.
\begin{theo}
	In the setting of \eqref{4.11n}, assume \textbf{[BST]} \eqref{4.3nn}, and assume \textbf{[FREE]} \eqref{4.8n} for every $\Gamma_m$, and assume \textbf{[GAP]} \eqref{n7.3} for all $F_{m,p}$. Then
	\begin{equation}\label{n7.8}
		\begin{split}
		\limsup_{\max(m,p)\to\infty}\Bv&\frac{1}{\dim_\mathbb{C}F_{m,p}}\ln\det(d-\Delta^{F_{m,p}})\\
		&-\int_\mathbb{R}\ln(d-\lambda)d\mu_{\mathrm{KM}}(\lambda)\Bv=0.
		\end{split}
	\end{equation}
\end{theo}

\begin{pro}
By applying \eqref{n7.6} to $(F_{m,p})_{p\in\mathbb{N}}\subset \mathscr{F}(X_m,\varepsilon)$ and replacing $(F_p)_{p\in\mathbb{N}}$ with $(F_{m,p})_{p\in\mathbb{N}}$ in \eqref{n7.7}, we obtain \eqref{n7.8} from \eqref{4.18n}.\qed
\end{pro}

\subsection{Asymptotic determinant in the Borel-Weil-Bott case}\label{ns7.5n}

We now state the main results of this section in the Borel-Weil-Bott setting.

\begin{theo}
	In the setting of \eqref{4.20nn}, assume \textbf{[BST]} \eqref{4.3nn} and \textbf{[GAP]} \eqref{n7.3}, and suppose that $G$ is compact connected. Then
	\begin{equation}
		\begin{split}
	\limsup_{m\to\infty}\sup_{\substack{\alpha\in \Lambda^*\cap\overline{C^+},\\
		F_{m,\alpha}\in \mathscr{F}(X_m,\varepsilon)}}\Bv&\frac{1}{\dim_{\mathbb{C}}F_{m,\alpha}}\ln\det(d-\Delta^{F_{m,\alpha}})\\
&-\int_\mathbb{R}\ln(d-\lambda)d\mu_{\mathrm{KM}}(\lambda)\Bv=0,
		\end{split}
	\end{equation}
	where, for each $m\in\mathbb{N}$, the supremum is taken over all irreducible representations of $G$ such that $F_{m,\alpha}$ satisfies \textbf{[GAP]}.
\end{theo}

\begin{theo}
	In the setting of \eqref{4.11nn}, assume \textbf{[FREE]} \eqref{4.10n}, and assume \textbf{[GAP]} \eqref{n7.3} for every $F_{p\alpha}$. Suppose, in addition, that either $G$ is compact connected simple and $\alpha$ is an arbitrary highest weight, or that $G$ is compact connected and $\alpha$ is regular. Then
	\begin{equation}
	\limsup_{p\to\infty}\lv\frac{1}{\dim_\mathbb{C}F_{p\alpha}}\ln\det(d-\Delta^{F_{p\alpha}})-\int_\mathbb{R}\ln(d-\lambda)d\mu_{\mathrm{KM}}(\lambda)\rv=0.
	\end{equation}
\end{theo}

\begin{theo}\label{nt6.6}
	In the setting of \eqref{4.20nn}, assume \textbf{[BST]} \eqref{4.3nn}, assume \textbf{[FREE]} \eqref{4.10n} for every $\Gamma_m$, and assume \textbf{[GAP]} \eqref{n7.3} for every $F_{m,p\alpha}$. Suppose, in addition, that either $G$ is compact connected simple and $\alpha$ is an arbitrary highest weight, or that $G$ is compact connected and $\alpha$ is regular. Then
	\begin{equation}
		\begin{split}
		\limsup_{\max(m,p)\to\infty}\Bv&\frac{1}{\dim_{\mathbb{C}}F_{m,p\alpha}}\ln\det(d-\Delta^{F_{m,p\alpha}})\\
&-\int_\mathbb{R}\ln(d-\lambda)d\mu_{\mathrm{KM}}(\lambda)\Bv=0.
		\end{split}
	\end{equation}
\end{theo}

\section{Kernel Operators}\label{s7n}

In this section, we study kernel operators on graph vector bundles, which play a role analogous to that of pseudodifferential operators on manifolds. In \S\,\ref{s7.1nn}, we introduce the endomorphism bundle. In \S\,\ref{s7.1n}, we define kernel operators and express an important family of them in terms of Chebyshev polynomials. In \S\,\ref{s7.2n}, we introduce several natural norms on kernel operators. In \S\,\ref{s7.3n}, we introduce several useful operators acting on kernel operators, including the commutator, gradient, and nonbacktracking operators. In \S\,\ref{s7.4n}, we study the spectral analysis of nonbacktracking operators. In \S\,\ref{s7.5n}, we study the spectral analysis of the commutator and gradient operators. In \S\,\ref{s7.6n}, we introduce a spectral gap condition.
 
Throughout this section, our main reference is Anantharaman
\cite[\S\S\,2, 3, 5]{MR3649482}, whose arguments we generalize to the vector bundle setting. For Chebyshev polynomials and the graph Laplacian, we refer to Einsiedler-Wald \cite[Lemma 3.20]{Arizona}. For the spectral analysis of nonbacktracking operators, we refer to Lubetzky-Peres \cite[\S\S\,3, 4]{MR3558308}, whose analysis we extend to the vector bundle setting. As a consequence, our spectral gap condition is slightly weaker than the one in Anantharaman-Le Masson \cite[\S\,1]{MR3322309}, in that it does not require a gap at $(-d)$, and therefore applies to additional examples, including bipartite graphs. Note that the space of oriented nonbacktracking paths is the graph theoretic counterpart of the cotangent phase space in the manifold setting, and the nonbacktracking walk is the counterpart of the geodesic flow.

\subsection{Endomorphism vector bundles}\label{s7.1nn}

Let $(X,F)$ be a $d$-regular graph and a unitary flat vector bundle of the form \eqref{2.19n}. 

Define the endomorphism bundle $\mathrm{End}(F)$ by
\begin{equation}\label{7.58}
	\begin{split}
		\mathrm{End}(F)&=\Gamma\backslash\big(\mathbb{T}_d\times\mathrm{End}(\mathbb{C}^\ell)\big)\\
		&=\big\{(y,A)\in \mathbb{T}_d\times\mathrm{End}(\mathbb{C}^\ell) \big\}\big/(y,A)\sim(\gamma y,\rho(\gamma)A\rho(\gamma)^{-1})\text{ for any }\gamma\in\Gamma.
	\end{split}
\end{equation}
Then, as in \eqref{2.20n}, we have
\begin{equation}\label{7.59}
	\begin{split}
		&C^\infty\big(X, \mathrm{End}(F)\big)\\
		&= C^\infty_{\Gamma}\big(\mathbb{T}_{d},\mathrm{End}(\mathbb{C}^\ell)\big)\subset C^\infty\big(\mathbb{T}_{d},\mathrm{End}(\mathbb{C}^\ell)\big)\\
		&=\big\{Q\mid Q(\gamma y)=\rho(\gamma)Q(y)\rho(\gamma)^{-1}\text{ for any } y\in\mathbb{T}_{d}, \gamma\in\Gamma\big\},\\
		&\Delta^{\mathrm{End}(F)}(Q)(y)=\sum_{y'\sim y}Q(y').
	\end{split}
\end{equation}

The space $C^\infty\big(X, \mathrm{End}(F)\big)$ is too small for the analysis below, in particular, it is not stable under composition with $\Delta^F$. We therefore enlarge the class of operators under consideration.

\subsection{Kernel vector bundles}\label{s7.1n}

Define \emph{the kernel space $K$} and \emph{the kernel vector bundle $K^{\mathrm{End}(F)}$} by
\begin{equation}\label{7.1n}
	\begin{split}
K=&\Gamma\backslash\mathbb{T}_d^2\\
=&\big\{(y,y')\in \mathbb{T}_d^2\big\}\big/(y,y')\sim(\gamma y,\gamma y')\ \text{for any }\gamma\in\Gamma,\\
K^{\mathrm{End}(F)}=&\Gamma\backslash\big(\mathbb{T}_d^2\times \mathrm{End}(\mathbb{C}^\ell)\big)\\
=&\big\{(y,y',A)\in \mathbb{T}_d^2\times\mathrm{End}(\mathbb{C}^\ell)\big\}\big/(y,y',A)\sim(\gamma y,\gamma y',\rho(\gamma)A\rho(\gamma)^{-1})\text{ for any }\gamma\in\Gamma.
	\end{split}
\end{equation}
We call a section $Q(y,y')\in C^\infty(K, K^{\mathrm{End}(F)})$ a \emph{kernel operator}. As in \eqref{2.20n}, we can identify $C^\infty(K, K^{\mathrm{End}(F)})$ with the space of $\Gamma$-equivariant $\mathrm{End}(\mathbb{C}^\ell)$-valued functions on $\mathbb{T}_d^2$, namely
\begin{equation}\label{7.2n}
	\begin{split}
	&C^\infty(K, K^{\mathrm{End}(F)})\\
	&= C^\infty_{\Gamma}\big(\mathbb{T}_{d}^2,\mathrm{End}(\mathbb{C}^\ell)\big)\subset C^\infty\big(\mathbb{T}_{d}^2,\mathrm{End}(\mathbb{C}^\ell)\big)\\
		&=\big\{Q \mid Q(\gamma y,\gamma  y')=\rho(\gamma)Q(y,y')\rho(\gamma)^{-1}\text{ for any } (y,y')\in\mathbb{T}_{d}^2, \gamma\in\Gamma\big\}.
	\end{split}
\end{equation}

A nonbacktracking path $(y_k,\cdots,y_1,y_0)$ in $\mathbb{T}_d$ is uniquely determined by its initial and terminal vertices $(y_k,y_0)\in\mathbb{T}_d^2$. Hence, for any $k\in\mathbb{N}$, we have an identification
\begin{equation}
	\begin{split}
		\{\text{nonbacktracking paths in }\mathbb{T}_d \text{ of length }k&\}\\
\cong\{(y,y')\in\mathbb{T}_d^2\mid d_{\mathbb{T}_d}(y,y')=k&\}.
	\end{split}
\end{equation}
From this point of view, we may regard $K$ defined in \eqref{7.1n} as the space of nonbacktracking paths in $X$. We denote
\begin{equation}\label{7.4}
	K_{(k)}=\{(y,y')\in K\mid d_{\mathbb{T}_d}(y,y')=k\},\quad K_{(\leqslant k)}=\mathsmaller{\bigcup}_{i=0}^{k}K_{(i)}.
\end{equation}
Accordingly, we have the direct sum decomposition
\begin{equation}\label{7.5nnn}
	\begin{split}
C^\infty(K, K^{\mathrm{End}(F)})&=\mathsmaller{\bigoplus}_{k=0}^\infty C^\infty(K_{(k)}, K^{\mathrm{End}(F)}),\\
Q&=Q_{(0)}+Q_{(1)}+\cdots,
	\end{split}
\end{equation}
where for $Q\in C^\infty(K, K^{\mathrm{End}(F)})$, we write $Q_{(k)}$ for its component in $C^\infty(K_{(k)}, K^{\mathrm{End}(F)})$. In particular, we have
\begin{equation}\label{7.10nn}
C^\infty(X,\mathrm{End}(F))\cong C^\infty(K_{(0)}, K^{\mathrm{End}(F)}).
\end{equation}

In what follows, we shall use both notations $Q(y,y')$ and $Q(y_k,\cdots,y_0)$ for a kernel operator, depending on whether we view it as a function of two endpoints or as a function of a nonbacktracking path.

We now introduce an important family of kernel operators $\mathrm{Id}_{(k)}\in C^\infty(K_{(k)},K^{\mathrm{End}(F)})$ for $k\in\mathbb{N}$ by
\begin{equation}\label{7.5nn}
	\mathrm{Id}_{(k)}(y,y')=\mathrm{Id}_{\mathbb{C}^\ell}\mathbbm{1}_{\{d_{\mathbb{T}_d}(y,y')=k\}}(y,y').
\end{equation}
Using Chebyshev polynomials \cite[Lemma 3.20]{Arizona}, we can express $\mathrm{Id}_{(k)}$ as an explicit function of $\Delta^F$ by
\begin{equation}\label{7.6nn}
	\begin{split}
	\mathrm{Id}_{(1)}&=\Delta^{F},\quad\mathrm{Id}_{(k)}=h_{k}(\Delta^F),\\
h_k(\lambda)&=\frac{\lv\pa B_{\mathbb{T}_d}(o,k)\rv}{(d-1)^{k/2}}\bigg(\frac{2}{d}\cos(k\theta)+\frac{(d-2)}{d}\frac{\sin\big((k+1)\theta\big)}{\sin(\theta)}\bigg),\\
\lambda&=2(d-1)^{1/2}\cos(\theta),
	\end{split}
\end{equation}
where $\partial B_{\mathbb{T}_d}(o,k)$ denotes the sphere of radius $k$ centered at $o$ and $\theta\in\mathbb{C}$.

By \eqref{2.20n} and \eqref{7.2n}, a kernel operator $Q\in C^\infty(K, K^{\mathrm{End}(F)})$ acts on a section $u\in C^\infty(X,F)$ by
\begin{equation}\label{7.5n}
	\begin{split}
(Qu)(y)=&\sum_{y'\in\mathbb{T}_d}Q(y,y')u(y')=\sum_{y'\in D,\gamma\in\Gamma}Q(y,\gamma y')u(\gamma y')\\
=&\sum_{y'\in D,\gamma\in\Gamma}Q(y,\gamma y')\rho(\gamma)u(y').
	\end{split}
\end{equation}
For $Q,Q'\in C^\infty(K,K^{\mathrm{End}(F)})$ and $u\in C^\infty(X,F)$,
\begin{equation}
	\begin{split}
\big(QQ'u\big)(y)&=\sum_{y'\in\mathbb{T}_d}Q(y,y')(Qu)(y')\\
&=\sum_{y',y''\in\mathbb{T}_d}Q(y,y')Q(y',y'')u(y'').
	\end{split}
\end{equation}
Therefore, the composition $QQ'$ is again a kernel operator such that
\begin{equation}\label{7.6n}	(QQ')(y,y')=\sum_{y''\in\mathbb{T}_d}Q(y,y'')Q'(y'',y').
\end{equation}

Finally, by \eqref{2.21n} and \eqref{7.5n}, we have
\begin{equation}\label{7.14nn}
\begin{split}
\big\langle Qu,u'\big\rangle_{L^2(X,F)}=&\sum_{y\in D}\big\langle(Qu)(y),u'(y)\big\rangle_{\mathbb{C}^\ell}\\
=&\sum_{y,y'\in D,\gamma\in\Gamma}\big\langle Q(y,\gamma y')\rho(\gamma)u(y'), u'(y)\big\rangle_{\mathbb{C}^\ell}.
\end{split}
\end{equation}

\subsection{Comparison of norms}\label{s7.2n}

By \eqref{7.2n}, we define normalized $L^{2}$- and $L^\infty$-norms on $C^\infty(K,K^{\mathrm{End}(F)})$ by
\begin{equation}\label{7.7n}
	\begin{split}
		\lV Q\rV_{L^2(K,K^{\mathrm{End}(F)})}^2&=\frac{1}{\dim_{\mathbb{C}}F}\sum_{y\in D,y'\in \mathbb{T}_d}\lV Q(y,y')\rV_{\mathrm{HS}(\mathbb{C}^\ell)}^2,\\
		\lV Q\rV_{L^\infty(K,K^{\mathrm{End}(F)})}^2&=\frac{1}{\dim_{\mathbb{C}}F_x}\max_{y\in D,y'\in \mathbb{T}_d}\lV Q(y,y')\rV_{\mathrm{HS}(\mathbb{C}^\ell)}^2,
	\end{split}
\end{equation}
where we use the dimension conventions in \eqref{2.14n}. Moreover, by \eqref{7.14nn}, we define the normalized Hilbert-Schmidt norm on $C^\infty(K,K^{\mathrm{End}(F)})$ by
\begin{equation}\label{7.9n}
	\lV Q\rV_{\mathrm{HS}(X,F)}^2=\frac{1}{\dim_{\mathbb{C}}F}\sum_{y,y'\in D}\BV\sum_{\gamma\in\Gamma}Q(y,\gamma y')\rho(\gamma)\BV_{\mathrm{HS}(\mathbb{C}^\ell)}^2,
\end{equation}
where $\lV\cdot \rV_{\mathrm{HS}(\mathbb{C}^\ell)}$ denotes the Hilbert-Schmidt norm on $\mathrm{End}(\mathbb{C}^\ell)$, given by
\begin{equation}
	\lV Q(y,y')\rV_{\mathrm{HS}(\mathbb{C}^\ell)}^2=\mathrm{Tr}^{\mathbb{C}^\ell}\big[Q(y,y')^*\cdot Q(y,y')\big].
\end{equation}

As we will see later, the norms $\lV \cdot\rV_{L^{2}(K,K^{\mathrm{End}(F)})}$ and
$\lV \cdot\rV_{L^\infty(K,K^{\mathrm{End}(F)})}$ are easy to compute, whereas $\lV \cdot\rV_{\mathrm{HS}(X,F)}$ is the norm that must be controlled. Hence, our basic strategy is to bound $\lV \cdot\rV_{\mathrm{HS}(X,F)}$ in terms of $\lV \cdot\rV_{L^2(K,K^{\mathrm{End}(F)})}$ and $\lV \cdot\rV_{L^\infty(K,K^{\mathrm{End}(F)})}$ under suitable assumptions.

\begin{lemma}
Let $Q\in C^\infty(K, K^{\mathrm{End}(F)})$. Then
\begin{equation}\label{7.10n}
	\begin{split}
		\lV Q\rV_{\mathrm{HS}(X,F)}^2=&\lV Q\rV_{L^2(K,K^{\mathrm{End}(F)})}^2\\
		&+\sum_{\substack{y,y'\in D,\  \gamma,\gamma'\in\Gamma\\ \gamma\neq \gamma'}}\frac{\mathrm{Tr}^{\mathbb{C}^\ell} }{\dim_{\mathbb{C}}F}\Big[\rho(\gamma')^{-1}Q(y,\gamma'y')^*Q(y,\gamma y')\rho(\gamma)\Big].
	\end{split}
\end{equation}
\end{lemma}

\begin{pro}
By direct computation, for any $(y,y')\in\mathbb{T}_d^2$,
	\begin{equation}\label{7.11n}
		\begin{split}
			\BV\sum_{\gamma\in\Gamma}Q(y,\gamma y')\rho(\gamma)\BV_{\mathrm{HS}(\mathbb{C}^\ell)}^2=&\sum_{\gamma\in\Gamma}\BV Q(y,\gamma y')\BV_{\mathrm{HS}(\mathbb{C}^\ell)}^2\\
			&+\sum_{\gamma\neq \gamma'}\mathrm{Tr}^{\mathbb{C}^\ell}\Big[\rho(\gamma')^{-1}Q(y,\gamma'y')^*Q(y,\gamma y')\rho(\gamma)\Big].
		\end{split}
	\end{equation}
Summing \eqref{7.11n} over $y,y'\in D$, by \eqref{7.7n} and \eqref{7.9n}, we get \eqref{7.10n}.\qed
\end{pro}

\begin{lemma}\label{L7.2}
	Let $Q\in C^\infty(K_{(\leqslant k)},K^{\mathrm{End}(F)})$. Then
	\begin{equation}\label{7.17}
		\lV Q\rV_{L^2(K_{(\leqslant k)},K^{\mathrm{End}(F)})}^2\leqslant	\lv B_{ \mathbb{T}_d}(o,k)\rv\lV Q\rV_{L^\infty(K,K^{\mathrm{End}(F)})}^2.
	\end{equation}
\end{lemma}

\begin{pro}
This follows immediately from \eqref{7.7n}.\qed
\end{pro}

\begin{prop}\label{p7.3}
Let $Q\in C^\infty(K_{(\leqslant k)},K^{\mathrm{End}(F)})$. If $k\leqslant\inf_{x\in X}\mathrm{inj}_x$, then
\begin{equation}\label{7.13n}
	\lV Q\rV_{\mathrm{HS}(X,F)}^2=\lV Q\rV_{L^2(K,K^{\mathrm{End}(F)})}^2.
\end{equation}
In general, we have
\begin{equation}\label{7.14n}
	\begin{split}
		\lV Q\rV_{\mathrm{HS}(X,F)}^2\leqslant& \lV Q\rV_{L^2(K,K^{\mathrm{End}(F)})}^2\\
		&+\frac{\lv\{x\in X\mid \mathrm{inj}_x\leqslant k\}\rv\lv B_{ \mathbb{T}_d}(o,k)\rv^2}{\lv X\rv}\lV Q\rV_{L^\infty(K,K^{\mathrm{End}(F)})}^2.
	\end{split}
\end{equation}
\end{prop}

\begin{pro}
If $k\leqslant \mathrm{inj}_x$, then in the sum $\sum_{\gamma\in\Gamma}Q(y,\gamma y')\rho(\gamma)$ there is at most one nonzero term. Hence,
	\begin{equation}\label{7.15n}
		\begin{split}
			\sum_{y'\in D}\BV\sum_{\gamma\in\Gamma}Q(y,\gamma y')\rho(\gamma)\BV_{\mathrm{HS}(\mathbb{C}^\ell)}^2&=\sum_{y'\in D,\gamma\in\Gamma}\lV Q(y,\gamma y')\rho(\gamma)\rV_{\mathrm{HS}(\mathbb{C}^\ell)}^2\\
			&=\sum_{y'\in \mathbb{T}_d}\lV Q(y,y')\rV_{\mathrm{HS}(\mathbb{C}^\ell)}^2,
		\end{split}
	\end{equation}
	where we use the unitarity of $\rho(\gamma)$.  Summing \eqref{7.15n} over $y\in D$, by \eqref{7.7n} and \eqref{7.9n}, we get \eqref{7.13n}.

If $k>\mathrm{inj}_x$, we have
\begin{equation}\label{7.16n}
\BV\sum_{\gamma\in\Gamma}Q(y,\gamma y')\rho(\gamma)\BV_{\mathrm{HS}(\mathbb{C}^\ell)}^2\leqslant \dim_{\mathbb{C}}F_x\lv B_{ \mathbb{T}_d}(o,k)\rv^2\lV Q\rV_{L^\infty(K,K^{\mathrm{End}(F)})}^2.
\end{equation}
 Combining \eqref{7.15n} and \eqref{7.16n}, we obtain \eqref{7.14n}.\qed
\end{pro}

\subsection{Useful operators}\label{s7.3n}

We now introduce several operators acting on $C^\infty(K, K^{\mathrm{End}(F)})$. Our main focus will be on the operator norm induced by $\lV\cdot\rV_{L^2(K, K^{\mathrm{End}(F)})}$, and the norm $\lV\cdot\rV_{L^\infty(K, K^{\mathrm{End}(F)})}$ will appear only once, in \eqref{7.38}. Unless otherwise specified, the adjoint of such an operator is always taken with respect to the inner product induced by $\lV\cdot\rV_{L^2(K, K^{\mathrm{End}(F)})}$, rather than by $\lV\cdot\rV_{\mathrm{HS}(X,F)}$.

\subsubsection{Cutting and orientation reversing operators}

For $k<k'$, we define the cutting operator
\begin{equation}\label{7.17n}
	\begin{split}
		&\mathrm{C}_{k',k}\colon C^{\infty}(K_{(k)},K^{\mathrm{End}(F)}) \to C^{\infty}(K_{(k')},K^{\mathrm{End}(F)}),\\
		&\mathrm{C}_{k',k}(Q)(y_{k'},\dots,y_0)=Q(y_k,\dots,y_0).
	\end{split}
\end{equation}

\begin{lemma}
The cutting operator is a scalar multiple of an isometry, that is, for $Q\in C^{\infty}(K_{(k)},K^{\mathrm{End}(F)})$,
\begin{equation}\label{7.18n}
	\begin{split}
		&\lV \mathrm{C}_{k',k}Q\rV_{L^2(K_{(k')},K^{\mathrm{End}(F)})}^2\\
		&=\begin{cases}
			(d-1)^{k'-k}\lV Q\rV_{L^2(K_{(k)},K^{\mathrm{End}(F)})}^2, &\text{if }k\neq 0,\\
				d(d-1)^{k'-1}\lV Q\rV_{L^2(K_{(k)},K^{\mathrm{End}(F)})}^2,&\text{if }k=0.
		\end{cases}
	\end{split}
\end{equation}
\end{lemma}

\begin{pro}
By \eqref{7.7n}, extending a nonbacktracking path from length $k$ to length $k'$ is equivalent to taking $(k'-k)$ further steps. If $k=0$, then the first step has $d$ choices, and each subsequent step has $(d-1)$ choices. If $k\neq 0$, then each step has $(d-1)$ choices. These explain the factors in \eqref{7.18n}.\qed
\end{pro}

We also define the orientation reversing operator
\begin{equation}\label{7.19n}
	\begin{split}
		&\mathrm{O}\colon C^\infty(K_{(k)},K^{\mathrm{End}(F)}) \to C^\infty(K_{(k)},K^{\mathrm{End}(F)}),\\
		&\mathrm{O}(Q)(y_{k},\dots,y_0)=Q(y_0,\dots,y_k).
	\end{split}
\end{equation}
Clearly we have the following result.

\begin{lemma}
The orientation reversing operator is an isometry, namely that for $Q\in C^{\infty}(K_{},K^{\mathrm{End}(F)})$,
\begin{equation}\label{7.20n}
	\lV\mathrm{O}(Q)\rV_{L^2(K,K^{\mathrm{End}(F)})}=\lV Q\rV_{L^2(K,K^{\mathrm{End}(F)})}.
\end{equation}
\end{lemma}

\subsubsection{Left and right multiplication operators}

By \eqref{7.5nn}, \eqref{7.6nn}, and \eqref{7.6n}, we define the left and right $\Delta^F$-multiplication operators by
\begin{equation}\label{7.21n}
	\begin{split}
		&\mathrm{L}_{\Delta^F}\colon C^{\infty}(K_{(k)},K^{\mathrm{End}(F)}) \to C^{\infty}(K_{(k-1)}\cup K_{(k+1)},K^{\mathrm{End}(F)}),\\
		&\mathrm{L}_{\Delta^F}(Q)(y,y')=\sum_{y''\sim y}Q(y'',y'),\\
			&\mathrm{R}_{\Delta^F}\colon C^{\infty}(K_{(k)},K^{\mathrm{End}(F)}) \to C^{\infty}(K_{(k-1)}\cup K_{(k+1)},K^{\mathrm{End}(F)}),\\
		&\mathrm{R}_{\Delta^F}(Q)(y,y')=\sum_{y''\sim y'}Q(y,y'').
	\end{split}
\end{equation}
Their two components can be written explicitly by
\begin{equation}\label{7.22n}
	\begin{split}
				&\mathrm{L}_{\Delta^F}(Q)(y_{k-1},\cdots,y_0)=\sum_{y_k\sim y_{k-1}, y_k\neq y_{k-2}}Q(y_{k},\cdots,y_0),\\
		&\mathrm{L}_{\Delta^F}(Q)(y_{k+1},y_k,\cdots,y_0)=Q(y_{k},\cdots,y_0),\\
			&\mathrm{R}_{\Delta^F}(Q)(y_k,\cdots,y_1)=\sum_{y_0\sim y_{1}, y_0\neq y_2}Q(y_k,\cdots,y_0),\\
		&\mathrm{R}_{\Delta^F}(Q)(y_{k},\cdots,y_0,y_{-1})=Q(y_{k},\cdots,y_0).
	\end{split}
\end{equation}

\begin{lemma}
Both $\mathrm{L}_{\Delta^F}$ and $\mathrm{R}_{\Delta^F}$ are selfadjoint, that is,
\begin{equation}\label{7.23n}
	\mathrm{L}_{\Delta^F}^*=\mathrm{L}_{\Delta^F},\quad \mathrm{R}_{\Delta^F}^*=\mathrm{R}_{\Delta^F}.
\end{equation}
\end{lemma}

\begin{pro}
To verify this, we compute by \eqref{7.7n} and \eqref{7.21n} that
\begin{equation}\label{7.28nn}
	\begin{split}
		&\big\langle \mathrm{L}_{\Delta^F}(Q),Q'\big\rangle_{L^2(K,K^{\mathrm{End}(F)})}\\
		&=\frac{1}{\dim_{\mathbb{C}}F}\sum_{y\in D,y'\in\mathbb{T}_d}\sum_{y''\sim y}\big\langle Q(y'',y'),Q'(y,y')\big\rangle_{\mathrm{HS}(\mathbb{C}^\ell)}\\
		&=\frac{1}{\dim_{\mathbb{C}}F}\sum_{y''\in D,y'\in\mathbb{T}_d}\sum_{y\sim y''}\big\langle Q(y'',y'),Q'(y,y')\big\rangle_{\mathrm{HS}(\mathbb{C}^\ell)}\\
		&=\big\langle Q,\mathrm{L}_{\Delta^F}(Q')\big\rangle_{L^2(K,K^{\mathrm{End}(F)})}.
	\end{split}
\end{equation}
The second identity follows by the same counting argument as in \eqref{5.1n}. More precisely, for any $y''\in \mathbb{T}_d$, there is a unique $\gamma\in \Gamma$ such that $\gamma y''\in D$. By \eqref{7.2n},
\begin{equation}\label{7.25n}
	\begin{split}
		&\big\langle Q(\gamma y'',y'),Q'(\gamma y,y')\big\rangle_{\mathrm{HS}(\mathbb{C}^\ell)}\\
		&=\big\langle \rho(\gamma)Q(y'',\gamma^{-1}y')\rho(\gamma)^{-1},\rho(\gamma)Q'(y,\gamma^{-1}y')\rho(\gamma)^{-1}\big\rangle_{\mathrm{HS}(\mathbb{C}^\ell)}\\
		&=\big\langle Q(y'',\gamma^{-1}y'),Q'(y,\gamma^{-1}y')\big\rangle_{\mathrm{HS}(\mathbb{C}^\ell)}.
	\end{split}
\end{equation}
Thus $(y\in D,y'\in\mathbb{T}_d, y''\sim y)$ is counted as $(\gamma y''\in D,\gamma y'\in\mathbb{T}_d, \gamma y\sim \gamma y'')$ on both sides. The selfadjointness of $\mathrm{R}_{\Delta^F}$ can be proved in the same way.\qed
\end{pro}

\subsubsection{Commutator and gradient operators}

Let $Q\in C^{\infty}(K,K^{\mathrm{End}(F)})$ be a kernel operator. Then its commutator $[\Delta^F,Q]$ with $\Delta^F$ is again a kernel operator. This commutator plays an important role in the later discussion of QE, as it is analogous to the commutator appearing in Egorov's theorem \cite[Theorem 15.2]{MR2952218}. Much of the remainder of this section is devoted to analyzing its properties.

By \eqref{7.21n}, we define the commutator map
\begin{equation}\label{7.27n}
	\begin{split}
		&\mathrm{ad}_{\Delta^F}\colon C^{\infty}(K_{(k)},K^{\mathrm{End}(F)}) \to C^{\infty}(K_{(k-1)}\cup K_{(k+1)},K^{\mathrm{End}(F)}),\\
		&\mathrm{ad}_{\Delta^F}(Q)(y,y')=\mathrm{L}_{\Delta^F}(Q)-\mathrm{R}_{\Delta^F}(Q)\\
		&\qquad\qquad\qquad\ =\sum_{y''\sim y}Q(y'',y')-\sum_{y''\sim y'}Q(y,y'').
	\end{split}
\end{equation}
As in \eqref{7.22n}, we can decompose $\mathrm{ad}_{\Delta^F}$ into two components,
\begin{equation}\label{7.33nn}
	\mathrm{ad}_{\Delta^F}=\nabla+\nabla^*.
\end{equation}
Here $\nabla$ denotes the gradient operator
\begin{equation}\label{7.33n}
	\begin{split}
	&\nabla\colon C^{\infty}(K_{(k)},K^{\mathrm{End}(F)}) \to C^{\infty}(K_{(k+1)},K^{\mathrm{End}(F)}),\\
		&\nabla(Q)(y_{k+1},\cdots,y_{0})=Q(y_{k+1},\cdots,y_{1})-Q(y_k,\cdots,y_0),
	\end{split}
\end{equation}
and $\nabla^*$ denotes its adjoint
\begin{equation}\label{7.34n}
	\begin{split}
&\nabla^*\colon C^{\infty}(K_{(k)},K^{\mathrm{End}(F)}) \to C^{\infty}(K_{(k-1)},K^{\mathrm{End}(F)}),\\
		&\nabla^*(Q)(y_k,\cdots,y_1)=\sum_{y_0\sim y_1,y_0\neq y_2}Q(y_k,\cdots,y_{0})\\
		&\qquad\qquad\qquad\qquad\quad-\sum_{y_{k+1}\sim y_k,y_{k+1}\neq y_{k-1}}Q(y_{k+1},\cdots,y_{1}).
	\end{split}
\end{equation}

\begin{lemma}
The commutator is selfadjoint,
\begin{equation}\label{7.36}
	\mathrm{ad}_{\Delta^F}^*=\mathrm{ad}_{\Delta^F}.
\end{equation}
\end{lemma}

\begin{pro}
This follows immediately from \eqref{7.23n} or \eqref{7.33nn}.\qed
\end{pro}

\begin{lemma}
Let $Q\in C^\infty(K_{(k)},K^{\mathrm{End}(F)})$. Then
\begin{equation}\label{7.37}
	\begin{split}
&\lV\nabla Q\rV_{L^2(K_{(k+1)},K^{\mathrm{End}(F)})}\\
&\leqslant\begin{cases}
	2(d-1)^{1/2}\lV Q\rV_{L^2(K_{(k)},K^{\mathrm{End}(F)})}, &\text{if }k\neq 0,\\
	2d^{1/2}\lV Q\rV_{L^2(K_{(k)},K^{\mathrm{End}(F)})}, &\text{if }k=0.
\end{cases}
	\end{split}
\end{equation}
\end{lemma}
\begin{pro}
This follows directly from \eqref{7.7n}, \eqref{7.33n}, the Cauchy-Schwarz inequality, and a counting argument similar to that in \eqref{7.28nn} and \eqref{7.25n}.\qed
\end{pro}

\begin{lemma}
Let $Q\in C^\infty(K,K^{\mathrm{End}(F)})$. Then
\begin{equation}\label{7.38}
	\begin{split}
\lV\mathrm{ad}_{\Delta^F}Q\rV_{L^2(K,K^{\mathrm{End}(F)})}&\leqslant4d^{1/2}\lV Q\rV_{L^2(K,K^{\mathrm{End}(F)})},\\
\lV\mathrm{ad}_{\Delta^F}Q\rV_{L^\infty(K,K^{\mathrm{End}(F)})}&\leqslant 2d\lV Q\rV_{L^\infty(K,K^{\mathrm{End}(F)})}.
	\end{split}
\end{equation}
\end{lemma}

\begin{pro}
	The first inequality follows from \eqref{7.33nn} and \eqref{7.37}, and the second inequality follows from \eqref{7.7n} and \eqref{7.27n}.\qed
\end{pro}

\subsubsection{Truncation operators}

We define the truncation operator
\begin{equation}\label{7.41}
	\begin{split}
		&\mathrm{T}\colon C^{\infty}(K_{(k)},K^{\mathrm{End}(F)}) \to C^{\infty}(K_{(k+2)},K^{\mathrm{End}(F)}),\\
		&\mathrm{T}(Q)(y_{k+1},\dots,y_{-1})=Q(y_k,\cdots,y_0),
	\end{split}
\end{equation}
and its adjoint
\begin{equation}\label{7.35n}
	\begin{split}
		&\mathrm{T}^*\colon C^{\infty}(K_{(k)},K^{\mathrm{End}(F)}) \to C^{\infty}(K_{(k-2)},K^{\mathrm{End}(F)}),\\
		&\mathrm{T}^*(Q)(y_{k-1},\cdots,y_{1})=\sum_{\substack{y_0\sim y_1, y_0\neq y_2,\\
				y_k\sim y_{k-1},y_k\neq y_{k-2}}}Q(y_k,y_{k-1}\cdots,y_{1},y_0).
	\end{split}
\end{equation}

\begin{lemma}
Let $Q\in C^{\infty}(K_{(k)},K^{\mathrm{End}(F)})$. Then
\begin{equation}\label{7.36n}
	\mathrm{T}(Q)=C_{k+2,k+1}\mathrm{O}C_{k+1,k}\mathrm{O}(Q).
\end{equation}
\end{lemma}
\begin{pro}
This can be verified immediately from \eqref{7.17n}, \eqref{7.19n}, and \eqref{7.41}.\qed
\end{pro}

\begin{lemma}
The truncation operator is a scalar multiple of an isometry. More precisely, for $Q\in C^{\infty}(K_{(k)},K^{\mathrm{End}(F)})$,
\begin{equation}\label{7.41n}
	\mathrm{T}^*\mathrm{T}(Q)=\begin{cases}
		(d-1)^2Q, &\text{if\ }k\neq0,\\
			d(d-1)Q, &\text{if\ }k=0.
	\end{cases}
\end{equation}
\end{lemma}

\begin{pro}
This follows directly from \eqref{7.18n}, \eqref{7.20n}, and \eqref{7.36n}.\qed
\end{pro}

We now express $\nabla^*$ in terms of $\nabla$ and $\mathrm{T}$.
\begin{lemma}
Let $Q\in C^{\infty}(K,K^{\mathrm{End}(F)})$.Then
\begin{equation}\label{7.38n}
\nabla^*(Q)=
-\frac{1}{d-1}\mathrm{T}^*\nabla (Q).
\end{equation}
\end{lemma}

\begin{pro}
By \eqref{7.33n} and \eqref{7.34n},
\begin{equation}
	\begin{split}
		&(\nabla^* Q)(y_k,\cdots,y_1)\\
		&=\frac{1}{d-1}\sum_{\substack{y_0\sim y_1,y_0\neq y_2\\ y_{k+1}\sim y_k,y_{k+1}\neq y_{k-1}
		}}\big(Q(y_k,\cdots,y_0)-Q(y_{k+1},\cdots,y_1)\big)\\
		&=\frac{1}{d-1}\sum_{\substack{y_0\sim y_1,y_0\neq y_2\\ y_{k+1}\sim y_k,y_{k+1}\neq y_{k-1}
		}}-\nabla (Q)(y_{k+1},\cdots,y_{0}),
	\end{split}
\end{equation}
which is exactly \eqref{7.38n} by \eqref{7.35n}.\qed
\end{pro}

 \subsubsection{Nonbacktracking walk operators}\label{s7.4.5}

We define the nonbacktracking walk operator
\begin{equation}\label{7.43n}
	\begin{split}
		&\mathrm{B}\colon C^\infty(K_{(k)},K^{\mathrm{End}(F)}) \to C^\infty(K_{(k)},K^{\mathrm{End}(F)}),\\
		&\mathrm{B}(Q)(y_k,\dots,y_0)=\sum_{y_{-1}\sim y_0,y_{-1}\neq y_1}Q(y_{k-1},\cdots,y_{-1}),
	\end{split}
\end{equation}
and its adjoint
\begin{equation}\label{7.44n}
	\begin{split}
		&\mathrm{B}^*\colon C^\infty(K_{(k)},K^{\mathrm{End}(F)}) \to C^\infty(K_{(k)},K^{\mathrm{End}(F)}),\\
		&\mathrm{B}^*(Q)(y_k,\dots,y_0)=\sum_{y_{k+1}\sim y_{k},y_{k+1}\neq y_{k-1}}Q(y_{k+1},\cdots,y_{1}).
	\end{split}
\end{equation}
In particular, if $Q_{(0)}\in C^\infty(K_{(0)},K^{\mathrm{End}(F)})$, then by \eqref{7.59}, under the identification \eqref{7.10nn}, we have
\begin{equation}\label{7.50nn}
\mathrm{B}(Q_{(0)})=\Delta^{\mathrm{End}(F)}Q_{(0)}.
\end{equation}

\begin{lemma}
Let $Q\in C^\infty(K_{(k)},K^{\mathrm{End}(F)})$. Then
\begin{equation}\label{7.45nn}
	\lV\mathrm{B}(Q)\rV_{L^2(K_{(k)},K^{\mathrm{End}(F)})}\leqslant \begin{cases}
(d-1)\lV Q\rV_{L^2(K_{(k)},K^{\mathrm{End}(F)})},&\text{if }k\neq 0,\\
d\lV Q\rV_{L^2(K_{(k)},K^{\mathrm{End}(F)})},&\text{if }k=0.
	\end{cases}
\end{equation}
\end{lemma}
\begin{pro}
This follows immediately from \eqref{7.7n}, \eqref{7.21n}, Cauchy-Schwarz inequality, and a counting argument similar to that in \eqref{7.28nn} and \eqref{7.25n}.\qed
\end{pro}

The operators $\mathrm{B}$ and $\mathrm{B}^*$ arise naturally in the analysis of $\nabla$. Indeed,
\begin{equation}\label{7.45n}
\lV \nabla Q\rV_{L^2(K,K^{\mathrm{End}(F)})}^2=\langle\nabla^*\nabla Q,Q\rangle_{L^2(K,K^{\mathrm{End}(F)})}.
\end{equation}
Expanding $\nabla^*\nabla$ yields the following identity.

\begin{lemma}
Let $Q\in C^\infty(K_{(k)},K^{\mathrm{End}(F)})$. Then
	\begin{equation}\label{7.46n}
		\nabla^*\nabla (Q)=\begin{cases}
		2(d-1)Q-\mathrm{B}(Q)-\mathrm{B}^*(Q),&\text{if }k\geqslant1,\\
			2\big(dQ-\mathrm{B}(Q)\big),&\text{if }k=0.
		\end{cases}
	\end{equation}
\end{lemma}

\begin{pro}
By \eqref{7.33n} and \eqref{7.34n},
\begin{equation}
	\begin{split}
		(\nabla^*&\nabla Q)(y_k,\cdots,y_0)\\
		=&\sum_{y_{-1}\sim y_0,y_{-1}\neq y_1}(\nabla Q)(y_{k},\cdots,y_{-1})\\
		&-\sum_{y_{k+1}\sim y_{k},y_{k+1}\neq y_{k-1}}(\nabla Q)(y_{k+1},\cdots,y_{0})\\
		=&\sum_{y_{-1}\sim y_0,y_{-1}\neq y_1}Q(y_{k},\cdots,y_{0})-Q(y_{k-1},\cdots,y_{-1})\\
		&-\sum_{y_{k+1}\sim y_{k},y_{k+1}\neq y_{k-1}}Q(y_{k+1},\cdots,y_1)- Q(y_k,\cdots,y_0),
	\end{split}
\end{equation}
which is exactly \eqref{7.46n}.\qed
\end{pro}

The spectral analysis of $\nabla$ via \eqref{7.45n} and \eqref{7.46n} involves both $\mathrm{B}$ and $\mathrm{B}^*$. We now introduce another estimate, which controls $\nabla$ using only $\mathrm{B}$.
\begin{lemma}
	Let $Q\in C^\infty(K_{(k)},K^{\mathrm{End}(F)})$. Then
	\begin{equation}\label{7.48n}
		\begin{split}
			&\lV\nabla Q\rV^2_{L^2(K_{(k+1)},K^{\mathrm{End}(F)})}\\
			&\geqslant\begin{cases}
				\frac{1}{d-1}\bV(d-1)Q-\mathrm{B}(Q)\bV^2_{L^2(K_{(k)},K^{\mathrm{End}(F)})},&\text{if }k\neq 0,\\
				\frac{1}{d}\bV dQ-\mathrm{B}(Q)\bV^2_{L^2(K_{(k)},K^{\mathrm{End}(F)})},&\text{if }k=0.
			\end{cases}
		\end{split}
	\end{equation}
\end{lemma}

\begin{pro}
	This is a direct consequence of \eqref{7.7n}, \eqref{7.33n}, \eqref{7.43n},  and Cauchy-Schwarz inequality.\qed
\end{pro}

The following results show that $\mathrm{B}$ has a filtration structure.
\begin{lemma}
Let $Q\in C^\infty(K_{(k)},K^{\mathrm{End}(F)})$ and $k'>k$. Then
\begin{equation}\label{7.50n}
	\begin{split}
	&\mathrm{B}\mathrm{C}_{k',k}(Q)(y_{k'},\cdots,y_{0})\\
	&=\begin{cases}
	\mathrm{C}_{k',k}\mathrm{B}(Q)(y_{k'},\cdots,y_{0}), &\text{if\ }k\neq 0,\\
		\mathrm{C}_{k',k}\mathrm{B}(Q)(y_{k'},\cdots,y_{0})-Q(y_1), &\text{if\ }k=0.
\end{cases}
	\end{split}
\end{equation}	
	Moreover, if $k\geqslant 2$, then
\begin{equation}\label{7.51n}
	\mathrm{B}\big(C^\infty(K_{(k)},K^{\mathrm{End}(F)})\big)\subseteq \mathrm{C}_{k,k-1}\big(C^\infty(K_{(k-1)},K^{\mathrm{End}(F)})\big).
\end{equation}
\end{lemma}

\begin{pro}
	By \eqref{7.17n} and \eqref{7.34n},
	\begin{equation}
		\begin{split}
			&\mathrm{B}(\mathrm{C}_{k',k}Q)(y_{k'},\cdots,y_{0})\\
			&=\sum_{y_{-1}\sim y_0,y_{-1}\neq y_1}(\mathrm{C}_{k',k}Q)(y_{k'-1},\cdots,y_0,y_{-1})\\
			&=\sum_{y_{-1}\sim y_0,y_{-1}\neq y_1}Q(y_{k-1},\cdots,y_0,y_{-1}),
		\end{split}
	\end{equation}
	which proves \eqref{7.50n}.
	
The inclusion \eqref{7.51n} follows directly from \eqref{7.34n}.\qed
\end{pro}

Note that \eqref{7.51n} fails for $k=1$ since $y_{-1}\neq y_1$ then depends on $y_{1}$.

We now consider the orthogonal decomposition
\begin{equation}\label{7.54n}
	\begin{split}
		&L^2(K_{(k+1)},K^{\mathrm{End}(F)})\\
		&=\mathrm{C}_{k+1,k}\big(L^2(K_{(k)},K^{\mathrm{End}(F)})\big)\oplus \mathrm{C}_{k+1,k}\big(L^2(K_{(k)},K^{\mathrm{End}(F)})\big)^\perp.
	\end{split}
\end{equation}

We add a subscript $k$ to the operators discussed above, for example $(\mathrm{ad}_{\Delta^F},\nabla,\mathrm{T}, \mathrm{B})$, to denote their restriction to $C^\infty(K_{(k)},K^{\mathrm{End}(F)})$.

\begin{lemma}
With respect to the decomposition \eqref{7.54n}, if $k\geqslant1$, then the operator $\mathrm{B}_{k+1}$ has the block form
\begin{equation}\label{7.55n}
	\mathrm{B}_{k+1}=\begin{pmatrix}
		\mathrm{B}_{k+1,11}&  \mathrm{B}_{k+1,12}\\
		\mathrm{B}_{k+1,21}&\mathrm{B}_{k+1,22}
	\end{pmatrix}=\begin{pmatrix}
		\mathrm{B}_{k}&  \mathrm{B}_{k+1,12}\\
		0&0
	\end{pmatrix}.
\end{equation}
Moreover, for $Q\in L^2(K_{(k+1)},K^{\mathrm{End}(F)})$,
\begin{equation}\label{7.56n}
	\lV\mathrm{B}_{k+1,12}(Q)\rV_{L^2(K_{(k+1)},K^{\mathrm{End}(F)})}\leqslant (d-1)\lV Q\rV_{L^2(K_{(k+1)},K^{\mathrm{End}(F)})}.
\end{equation}
\end{lemma}

\begin{pro}
By \eqref{7.50n} and \eqref{7.51n}, we obtain \eqref{7.55n}. Since $\mathrm{B}_{k+1,12}$ is a restriction of $\mathrm{B}_{k+1}$, then \eqref{7.56n} follows directly from \eqref{7.45nn}.\qed
\end{pro}

From \eqref{7.55n}, we see that $\mathrm{B}$ has a useful iterative structure. Thus, its spectral analysis can be reduced inductively to that of $\mathrm{B}_1$ and $\mathrm{B}_0$.

\subsection{Spectral analysis of nonbacktracking operators}\label{s7.4n}

Recalling \eqref{7.59}, \eqref{7.10nn}, and \eqref{7.50nn}, we have $\mathrm{B}_0=\Delta^{\mathrm{End}(F)}$. Then similarly to \eqref{3.27}, we can choose an orthonormal basis of eigensections
\begin{equation}\label{7.60}
\Delta^{\mathrm{End}(F)}(Q_{(0),i})=\lambda_{\mathrm{End}(F),i}Q_{(0),i},\quad \lV Q_{(0),i}\rV_{L^2(K_{0},K^{\mathrm{End}(F)})}^2=1.
\end{equation}

We write $\mathrm{mult}_{\cdot}(\cdot)$ for the algebraic multiplicity. We are now ready to describe the Jordan structure of $\mathrm{B}_1$.
\begin{prop}\label{p7.16}
The operator $\mathrm{B}_1$ acting on $L^2(K_{(1)},K^{\mathrm{End}(F)})$ is unitarily equivalent to a block diagonal operator of the form
\begin{equation}\label{7.61}
\mathrm{diag}\Big((d-1)\cdots,-(d-1)\cdots,1\cdots,-1\cdots,\cdots\big(\begin{smallmatrix}
	\theta_i& a_i\\
	0&\theta_i'
\end{smallmatrix}\big)\cdots\Big).
\end{equation}
Here
\begin{equation}\label{7.62}
\begin{split}
&\mathrm{mult}_{\mathrm{B}_1}\big(\pm(d-1)\big)=\mathrm{mult}_{\Delta^{\mathrm{End}(F)}}(\pm d),\\
&\mathrm{mult}_{\mathrm{B}_1}(\pm1)=(\tfrac{d}{2}-1)\lv X\rv\ell^2+\mathrm{mult}_{\Delta^{\mathrm{End}(F)}}(\mp d),
\end{split}
\end{equation}
and the number of blocks of the form $\big(\begin{smallmatrix}
	\theta_i & a_i\\
	0 & \theta_i'
\end{smallmatrix}\big)$
appearing in \eqref{7.61} is
\begin{equation}\label{7.62n}
\lv X\rv\ell^2-\mathrm{mult}_{\Delta^{\mathrm{End}(F)}}(d)-\mathrm{mult}_{\Delta^{\mathrm{End}(F)}}(-d).
\end{equation}
Moreover, the diagonal entries $\theta_i$ and $\theta_i'$ are the two roots of
\begin{equation}\label{7.63}
	\theta^2-\lambda_{\mathrm{End}(F),i}\theta+d-1=0,
\end{equation}
while the off-diagonal entry $a_i$ satisfies
\begin{equation}\label{7.64}
	\lv a_i\rv_\mathbb{C}=\begin{cases}
		\big(d^2-(\lambda_{\mathrm{End}(F),i})^2\big)^{1/2}, &\text{if } 	(\lambda_{\mathrm{End}(F),i}^{})^2\geqslant4(d-1),\\
		(d-2),&\text{if } 	(\lambda_{\mathrm{End}(F),i})^2\leqslant4(d-1).
	\end{cases}
\end{equation}
\end{prop}

\begin{pro}

First, we have the orthogonal decomposition
	\begin{equation}
		\begin{split}
			L^2(K_{(1)},K^{\mathrm{End}(F)})&=L^2_+(K_{(1)},K^{\mathrm{End}(F)})\oplus L^2_-(K_{(1)},K^{\mathrm{End}(F)}),\\
			L^2_\pm(K_{(1)},K^{\mathrm{End}(F)})&=\big\{Q_{(1)}\mid Q_{(1)}(y_0,y_1)=\pm Q_{(1)}(y_1,y_0) \text{ for any }(y_0,y_1)\in \mathbb{T}_d^2\big\}.
		\end{split}
	\end{equation}
	Define the star subspaces $S_\pm(K_{(1)},K^{\mathrm{End}(F)})\subset L^2_\pm(K_{(1)},K^{\mathrm{End}(F)})$ by
	\begin{equation}\label{7.67}
		\begin{split}
			S_\pm(K_{(1)},K^{\mathrm{End}(F)})&=\mathrm{span}\big\{Q_{(1),y,A}^\pm\mid y\in D, A\in\mathrm{End}(\mathbb{C}^\ell)\big\},\\ Q_{(1),y,A}^\pm(y_1,y_0)&=\begin{cases}\rho(\gamma)A\rho(\gamma)^{-1},&\text{if }y_0=\gamma y\text{ for }\gamma\in \Gamma,\\
				\pm \rho(\gamma)A\rho(\gamma)^{-1},&\text{if }y_1=\gamma y\text{ for }\gamma\in \Gamma,\\
				0,&\text{otherwise}.
			\end{cases}
		\end{split}
	\end{equation}
	
	For any $Q_{(1)}^{\pm}\in L^2_\pm(K_{(1)},K^{\mathrm{End}(F)})$, we compute that
	\begin{equation}
		\big\langle Q_{(1)}^{\pm},Q_{(1),y,A}^\pm\big\rangle_{L^2(K_{(1)},K^{\mathrm{End}(F)})}=2\sum_{y'\sim y}\mathrm{Tr}^{\mathrm{C}^\ell}[A^*Q_{(1)}^\pm(y,y')].
	\end{equation}
	Therefore, if $Q_{(1)}^\pm\perp S_\pm(K_{(1)},K^{\mathrm{End}(F)})$, then
	\begin{equation}
		\sum_{y'\sim y}Q^\pm_{(1)}(y,y')=0.
	\end{equation}
By \eqref{7.43n}, this shows that
	\begin{equation}\label{7.70}
			\mathrm{B}_1(Q_{(1)}^\pm)=\mp Q_{(1)}^\pm,\quad\text{if }Q_{(1)}^\pm\in  L^2_\pm(K_{(1)},K^{\mathrm{End}(F)})\cap S_\pm(K_{(1)},K^{\mathrm{End}(F)})^\perp.
	\end{equation}

	Now we compute the dimensions of these spaces. Suppose that $\{A_y\in\mathrm{End}(\mathbb{C}^\ell) \}_{y\in D}$ satisfies
	\begin{equation}\label{7.71}
		\sum_{y\in D} Q_{(1),y,A_y}^\pm=0.
	\end{equation}
Define $Q_{(0),\{A_y\in\mathrm{End}(\mathbb{C}^\ell) \}_{y\in D}}\in L^2_\pm(K_{(0)},K^{\mathrm{End}(F)})$ by
\begin{equation} Q_{(0),\{A_y\in\mathrm{End}(\mathbb{C}^\ell) \}_{y\in D}}(\gamma y)=\rho(\gamma)A_y\rho(\gamma)^{-1}, \quad\text{for any }y\in D, \gamma\in\Gamma.
\end{equation}
Then by \eqref{7.67}, condition \eqref{7.71} is equivalent to
	\begin{equation}\label{7.73}
Q_{(0),\{A_y\in\mathrm{End}(\mathbb{C}^\ell) \}_{y\in D}}(y_0)\pm Q_{(0),\{A_y\in\mathrm{End}(\mathbb{C}^\ell) \}_{y\in D}}(y_1)=0
	\end{equation}
for any edge $(y_1,y_0)$. In other words, the section $Q_{(0),\{A_y\in\mathrm{End}(\mathbb{C}^\ell) \}_{y\in D}}$ is an eigensection of $\Delta^{\mathrm{End}(F)}$ with eigenvalue $\mp d$. Hence
	\begin{equation}
		\begin{split}
			\dim_{\mathbb{C}}S_\pm(K_{(1)},K^{\mathrm{End}(F)})&=\lv X\rv\ell^2-\mathrm{mult}_{\Delta^{\mathrm{End}(F)}}(\mp d),\\ \dim_{\mathbb{C}}L^2_\pm(K_{(1)},K^{\mathrm{End}(F)})&\cap S_\pm(K_{(1)},K^{\mathrm{End}(F)})^\perp\\
			&=(\tfrac{d}{2}-1)\lv X\rv\ell^2+\mathrm{mult}_{\Delta^{\mathrm{End}(F)}}(\mp d).
		\end{split}
	\end{equation}
Combining this with \eqref{7.70}, we get the $\pm1$ entries in \eqref{7.61} and the second identity of \eqref{7.62}.

By \eqref{7.70}, it remains to analyze the restriction of $\mathrm{B}_1$ to the star spaces \eqref{7.67}. Their sum satisfies
\begin{equation}\label{7.75}
	\begin{split}
		&S_+(K_{(1)},K^{\mathrm{End}(F)})\oplus S_-(K_{(1)},K^{\mathrm{End}(F)})\\
		&=\mathrm{span}\big\{\mathrm{C}_{1,0}Q_{(0)}, \mathrm{O}(\mathrm{C}_{1,0}Q_{(0)})\mid Q_{(0)}\in L^2(K_{(0)},K^{\mathrm{End}(F)})\big\},
	\end{split}
\end{equation}
where $\mathrm{C}_{1,0}$ and $\mathrm{O}$ are defined in \eqref{7.17n} and \eqref{7.19n}. This space can be interpreted as the space generated by functions depending only on the tail or the head. 

We compute the inner product of two such functions. By \eqref{7.18n}, \eqref{7.19n}, and \eqref{7.20n},
\begin{equation}\label{7.76}
	\begin{split}
		&\big\langle \mathrm{O}\mathrm{C}_{1,0}Q_{(0)},\mathrm{O}\mathrm{C}_{1,0}Q'_{(0)} \big\rangle_{L^2(K_{(1)},K^{\mathrm{End}(F)})}\\
		&=\big\langle \mathrm{C}_{1,0}Q_{(0)},\mathrm{C}_{1,0}Q'_{(0)} \big\rangle_{L^2(K_{(1)},K^{\mathrm{End}(F)})}\\
		&=d\big\langle Q_{(0)},Q_{(0)}'\big\rangle_{L^2(K_{(0)},K^{\mathrm{End}(F)})},
	\end{split}
\end{equation}	
and also,
\begin{equation}\label{7.76n}
	\begin{split}
				&\big\langle \mathrm{C}_{1,0}Q_{(0)},\mathrm{O}(\mathrm{C}_{1,0}Q_{(0)}') \big\rangle_{L^2(K_{(1)},K^{\mathrm{End}(F)})}\\
	&=\frac{1}{\dim_{\mathbb{C}}F}\sum_{y_0\in D,\ y_0\sim y_1} \mathrm{Tr}^{\mathbb{C}^\ell}\big[\mathrm{O}(\mathrm{C}_{1,0}Q_{(0)}')(y_1,y_0)^*\mathrm{C}_{1,0}Q_{(0)}(y_1,y_0)\big]\\
	&=\frac{1}{\dim_{\mathbb{C}}F}\sum_{y_0\in D,\ y_0\sim y_1} \mathrm{Tr}^{\mathbb{C}^\ell}\big[Q_{(0)}'(y_1)^*Q_{(0)}(y_0)\big]\\
	&=\big\langle Q_{(0)},\Delta^{\mathrm{End}(F)}Q_{(0)}'\big\rangle_{L^2(K_{(0)},K^{\mathrm{End}(F)})}.
	\end{split}
\end{equation}
Taking together \eqref{7.60}, \eqref{7.75}, \eqref{7.76}, and \eqref{7.76n}, we obtain an orthogonal decomposition
\begin{equation}\label{7.77}
	\begin{split}
	&S_+(K_{(1)},K^{\mathrm{End}(F)})\oplus S_-(K_{(1)},K^{\mathrm{End}(F)})\\
	&=\mathsmaller{\bigoplus}_i\mathrm{span}\big\{\mathrm{C}_{1,0}Q_{(0),i}, \mathrm{O}(\mathrm{C}_{1,0}Q_{(0),i})\big\}.
	\end{split}
\end{equation}

We now determine the one dimensional subspaces in this decomposition. Suppose that $\mathrm{C}_{1,0}Q_{(0),i}=C\mathrm{O}(\mathrm{C}_{1,0}Q_{(0),i})$ for a constant $C$. Applying  this relation to edges $(y_1,y_0)$ and $(y_0,y_1)$, we see that $C^2=1$, hence
\begin{equation}\label{7.78n}
Q_{(0),i}(y_0)=\pm Q_{(0),i}(y_1).
\end{equation}
Then, similarly to \eqref{7.73}, the section $Q_{(0),i}$ must be an eigensection of $\Delta^{\mathrm{End}(F)}$ with eigenvalue $\pm d$. By \eqref{7.50n},
\begin{equation}\label{7.78}
	\mathrm{B}_1\mathrm{C}_{1,0}Q_{(0)}=\mathrm{C}_{1,0}\mathrm{B}_0Q_{(0)}-\mathrm{O}\mathrm{C}_{1,0}Q_{(0)}.
\end{equation}
Substituting $Q_{(0),i}$ for $Q_{(0)}$ in \eqref{7.78}, we get
\begin{equation}\label{7.79n}
\mathrm{B}_1(\mathrm{C}_{1,0}Q_{(0),i})=\pm(d-1)\mathrm{C}_{1,0}Q_{(0),i}.
\end{equation}
This gives the $\pm(d-1)$ entries in \eqref{7.61} and the first identity in \eqref{7.62}. Then \eqref{7.62n} follows from a dimension count.

We now consider the action of $\mathrm{B}_1$ on the two dimensional subspaces in \eqref{7.77}. By \eqref{7.78}, we have
\begin{equation}\label{7.79}
\begin{split}
\mathrm{B}_1\mathrm{C}_{1,0}Q_{(0),i}&=\mathrm{C}_{1,0}\Delta^{\mathrm{End}(F)}Q_{(0),i}-\mathrm{O}\mathrm{C}_{1,0}Q_{(0),i}\\
&=\lambda_{\mathrm{End}(F),i}\mathrm{C}_{1,0}Q_{(0),i}-\mathrm{O}\mathrm{C}_{1,0}Q_{(0),i},
\end{split}
\end{equation}
and also,
\begin{equation}\label{7.80}
\begin{split}
\mathrm{B}_1\mathrm{O}\mathrm{C}_{1,0}Q_{(0)}(y_1,y_0)=&\sum_{y_{-1}\sim y_0,y_{-1}\neq y_1}\mathrm{O}\mathrm{C}_{1,0}Q_{(0)}(y_0,y_{-1})\\
=&\sum_{y_{-1}\sim y_0,y_{-1}\neq y_1}\mathrm{C}_{1,0}Q_{(0)}(y_{-1},y_0)\\
=&(d-1)Q_{(0)}(y_0).
\end{split}
	\end{equation}
Combining \eqref{7.60}, \eqref{7.75}, \eqref{7.76}, \eqref{7.79}, \eqref{7.80}, we see that under the basis
\begin{equation}\label{7.86n}
(e_1,e_2)=(\mathrm{C}_{1,0}Q_{(0),i}/d^{1/2},\mathrm{O}\mathrm{C}_{1,0}Q_{(0),i}/d^{1/2}),
\end{equation}
the operator $\mathrm{B}_1$ and the Gram matrix $G$ take the form
\begin{equation}\label{7.82}
\mathrm{B}_1=\Big(\begin{smallmatrix}
da& (d-1)\\
-1&0
\end{smallmatrix}\Big),\quad G=\Big(\begin{smallmatrix}
1& a\\
a&1
\end{smallmatrix}\Big),\quad a=\lambda_{\mathrm{End}(F),i}/d.
	\end{equation}
	From \eqref{7.82}, we obtain the characteristic equation \eqref{7.63}. We choose an orthonormal basis
	\begin{equation}\label{7.88}
(e_1',e_2')=\big(e_1,\big(e_2-ae_1\big)\big/b\big),\quad b=(1-a^2)^{1/2}.
	\end{equation}
Since $e_2=ae_1'+be_2'$, the matrix of $\mathrm{B}_1$ in this new basis becomes
	\begin{equation}\label{7.89}
\mathrm{B}_1=\Big(\begin{smallmatrix}
	(d-1)a& (d-1)b\\
	-b&a
\end{smallmatrix}\Big).
	\end{equation}
We compute the Hilbert-Schmidt norms of $\mathrm{B}_1$ in \eqref{7.89} and $\big(\begin{smallmatrix}
	\theta_i & a_i\\
	0 & \theta_i'
\end{smallmatrix}\big)$ in \eqref{7.61},
\begin{equation}\label{7.85}
	\begin{split}
\lV\mathrm{B}_1\rV^2_{\mathrm{HS}(\mathrm{span}\{e_1,e_2\})}=&(d-1)^2a^2+(d-1)^2b^2\\
&+a^2+b^2=(d-1)^2+1,\\
\bV\big(\begin{smallmatrix}
	\theta_i & a_i\\
	0 & \theta_i'
\end{smallmatrix}\big)\bV_{\mathrm{HS}(\mathrm{span}\{e_1,e_2\})}^2=&\lv\theta_i\rv_\mathbb{C}^2+\lv\theta_i'\rv_\mathbb{C}^2+\lv a_i\rv_\mathbb{C}^2.
	\end{split}
\end{equation}
By \eqref{7.63}, 
\begin{equation}\label{7.86}
\lv\theta_i\rv_\mathbb{C}^2+\lv\theta_i'\rv_\mathbb{C}^2=\begin{cases}
	(\lambda_{\mathrm{End}(F),i})^2-2(d-1), &\text{if } 	(\lambda_{\mathrm{End}(F),i})^2\geqslant4(d-1),\\
	2(d-1),&\text{if } 	(\lambda_{\mathrm{End}(F),i})^2\leqslant4(d-1).
\end{cases}
\end{equation}
The Hilbert-Schmidt norm is invariant under unitary equivalence, hence, by \eqref{7.85} and \eqref{7.86}, we get \eqref{7.64}.\qed
\end{pro}

\begin{prop}\label{L7.17}
We have
\begin{equation}\label{7.87}
	\begin{split}
\ker(\mathrm{ad}_{\Delta^F})=&\ker(\nabla)\\
=&\mathsmaller{\bigoplus_{k=0}^\infty} \mathrm{C}_{k,0}\ker\big(d\mathrm{Id}-\Delta^{\mathrm{End}(F)}\big)\\
=&\ker\big(d\mathrm{Id}-\mathrm{B}_0\big)\mathsmaller{\bigoplus_{k=1}^\infty}\ker\big((d-1)\mathrm{Id}-\mathrm{B}_k\big)\\
=&\ker\big(d\mathrm{Id}-\mathrm{B}_0\big)\mathsmaller{\bigoplus_{k=1}^\infty}\ker\big((d-1)\mathrm{Id}-\mathrm{B}_k^*\big).
	\end{split}
\end{equation}
\end{prop}

\begin{pro}
	
By \eqref{7.33nn}, \eqref{7.41n}, and \eqref{7.38n}, if $\mathrm{ad}_{\Delta^F}Q=0$, then
	\begin{equation}
		\begin{split}
&\bV\nabla(Q_{(k)})\bV_{L^2(K_{(k+1)},K^{\mathrm{End}(F)})}^2\\
&=
\begin{cases}
		\bV\nabla (Q_{(k+2)})\bV_{L^2(K_{(k+3)},K^{\mathrm{End}(F)})}^2,&\text{if } k\neq0,\\
	\tfrac{d}{d-1}\bV\nabla (Q_{(2)})\bV_{L^2(K_{(3)},K^{\mathrm{End}(F)})}^2,&\text{if }k=0.
\end{cases}
		\end{split}
	\end{equation}
	Since $Q\in L^2(K,K^{\mathrm{End}(F)})$, we have $\lim_{k\to\infty}\Vert\nabla (Q_{(k)})\Vert_{L^2(K_{(k+1)},K^{\mathrm{End}(F)})}^2=0$, and therefore $\nabla Q=0$ and
	\begin{equation}\label{7.93}
\ker(\mathrm{ad}_{\Delta^F})\subseteq\ker(\nabla).
	\end{equation}
By \eqref{7.33n}, we have
	\begin{equation}\label{7.94n}
\ker(\nabla_k)=\mathsmaller{\bigoplus_{k=0}^\infty} \mathrm{C}_{k,0}\ker\big(d\mathrm{Id}-\Delta^{\mathrm{End}(F)}\big).
	\end{equation}
Also, we have
\begin{equation}\label{7.92}
	\mathsmaller{\bigoplus_{k=0}^\infty} \mathrm{C}_{k,0}\ker\big(d\mathrm{Id}-\Delta^{\mathrm{End}(F)}\big)\subseteq \ker(\mathrm{ad}_{\Delta^F}).
\end{equation}	
By \eqref{7.93}, \eqref{7.94n}, and  \eqref{7.92}, we get the first and second identity in \eqref{7.87}.
		
	By \eqref{7.48n},
\begin{equation}\label{7.90}
\ker(\nabla_k)\subseteq\begin{cases}
\ker\big((d-1)\mathrm{Id}-\mathrm{B}_k\big),&\text{if }k\neq0,\\
\ker\big(d\mathrm{Id}-\mathrm{B}_0\big), &\text{if }k=0.
\end{cases}
\end{equation}
By \eqref{7.55n}, \eqref{7.61}, \eqref{7.78n}, and \eqref{7.79n}, we have
\begin{equation}\label{7.91}
	\begin{split}
\mathrm{C}_{k,0}\ker\big(d\mathrm{Id}-\Delta^{\mathrm{End}(F)}\big)=\begin{cases}
			\ker\big((d-1)\mathrm{Id}-\mathrm{B}_k\big),&\text{if }k\neq0,\\
			\ker\big(d\mathrm{Id}-\mathrm{B}_0\big), &\text{if }k=0.
		\end{cases}
	\end{split}
\end{equation}
Combining \eqref{7.46n}, \eqref{7.90}, and \eqref{7.91}, we get the third and fourth identity in \eqref{7.87}.\qed
\end{pro}

By the identity $\ker((d-1)\mathrm{Id}-\mathrm{B}_k)=\ker((d-1)\mathrm{Id}-\mathrm{B}_k^*)$ from \eqref{7.87}, the operator $((d-1)\mathrm{Id}-\mathrm{B}_k)$ is invertible on $\ker(\mathrm{ad}_{\Delta^F})^\perp$. We now estimate the norm of its inverse.

By Proposition \ref{p7.16}, we see that the nontrivial top eigenvalue $\lambda_{\mathrm{End}(F),\mathrm{nt}}$ of $\Delta^{\mathrm{End}(F)}$ defined in \eqref{7.90n} plays an important role in the analysis of $\mathrm{B}$.

\begin{prop}\label{p7.18n}
Let $Q\in L^2(K_{(k)},K^{\mathrm{End}(F)})\cap\ker(\mathrm{ad}_{\Delta^F})^\perp$. Then
\begin{equation}\label{7.94}
	\begin{split}
\bV\big((d-1)\mathrm{Id}-\mathrm{B}_{k}\big)^{-1}Q\bV_{L^2(K_{(k)},K^{\mathrm{End}(F)})}&\\
\leqslant C_{k,\lambda_{\mathrm{End}(F),\mathrm{nt}}}\lV Q\rV_{L^2(K_{(k)},K^{\mathrm{End}(F)})}&,
	\end{split}
	\end{equation}
	where the constant depends only on $k$ and $\lambda_{\mathrm{End}(F),\mathrm{nt}}$.
\end{prop}

\begin{pro}
By \eqref{7.55n}, the resolvent has the block form
\begin{equation}
	\begin{split}
		&\big((d-1)\mathrm{Id}-\mathrm{B}_{k+1}\big)^{-1}\\
		&=\begin{pmatrix}
			\big((d-1)\mathrm{Id}-\mathrm{B}_{k}\big)^{-1}&\big((d-1)\mathrm{Id}-\mathrm{B}_{k}\big)^{-1}(d-1)^{-1}\mathrm{B}_{k,12}\\
			0&(d-1)^{-1}\mathrm{Id}
		\end{pmatrix}.
	\end{split}
\end{equation}
Then by \eqref{7.56n}, we obtain an inductive inequality
\begin{equation}
	\begin{split}
		&\bV\big((d-1)\mathrm{Id}-\mathrm{B}_{k+1}\big)^{-1}\bV_{\mathrm{End}(L^2(K_{(k+1)},K^{\mathrm{End}(F)}))}\\
		&\leqslant 2\bV\big(\mathrm{Id}-(d-1)^{-1}\mathrm{B}_{k}\big)^{-1}\bV_{\mathrm{End}(L^2(K_{(k)},K^{\mathrm{End}(F)}))}+(d-1)^{-1},
	\end{split}
\end{equation}
which, together with \eqref{7.61} and \eqref{7.63}, implies \eqref{7.94}.\qed
\end{pro}

\subsection{Spectral analysis of the gradient and commutator operators}\label{s7.5n}

Let us form the orthogonal projection
\begin{equation}\label{7.99}
P^{\ker(\mathrm{ad}_{\Delta^F})}\colon L^2(K,K^{\mathrm{End}(F)})\rightarrow \ker(\mathrm{ad}_{\Delta^F}).
\end{equation}

\begin{prop}\label{p7.19n}
Let $Q\in L^2(K_{(k)},K^{\mathrm{End}(F)})$. Then
	\begin{equation}\label{7.98}
		\begin{split}
&\bV Q-P^{\ker(\mathrm{ad}_{\Delta^F})}Q\bV_{L^2(K_{(k)},K^{\mathrm{End}(F)})}\\
&\leqslant C_{k,\lambda_{\mathrm{End}(F),\mathrm{nt}}}\lV\nabla Q\rV_{L^2(K_{(k+1)},K^{\mathrm{End}(F)})}.
		\end{split}
	\end{equation}	
\end{prop}

\begin{pro}
By the identity $\ker(\mathrm{ad}_{\Delta^F})=\ker(\nabla)$ from \eqref{7.87}, it suffices to consider the case $Q\in \ker(\mathrm{ad}_{\Delta^F})^\perp$. Under this assumption, by the identity $\ker((d-1)\mathrm{Id}-\mathrm{B}_k)=\ker((d-1)\mathrm{Id}-\mathrm{B}_k^*)$ from \eqref{7.87}, we have $((d-1)\mathrm{Id}-\mathrm{B}_{k})Q\in \ker(\mathrm{ad}_{\Delta^F})^\perp$. Hence, we can write
	\begin{equation}\label{7.100}
Q= \big((d-1)\mathrm{Id}-\mathrm{B}_{k}\big)^{-1}\big((d-1)\mathrm{Id}-\mathrm{B}_{k}\big)Q.
	\end{equation}
Then \eqref{7.48n}, \eqref{7.94}, and \eqref{7.100} imply \eqref{7.98}.\qed
\end{pro}

\begin{prop}\label{p7.20n}
Let $Q\in L^2(K, K^{\mathrm{End}(F)})\cap \ker(\mathrm{ad}_{\Delta^F})^\perp$ satisfy
\begin{equation}\label{7.101}
\lV\mathrm{ad}_{\Delta^F}Q\rV_{L^2(K, K^{\mathrm{End}(F)})}\leqslant\varepsilon \lV Q\rV_{L^2(K, K^{\mathrm{End}(F)})}
\end{equation}
for some $0\leqslant\varepsilon<1$. Then for any $k\in\mathbb{N}$,
	\begin{equation}\label{7.103}
		\lV Q_{(k)}\rV_{L^2(K_{(k)}, K^{\mathrm{End}(F)})}\leqslant C_{k,\lambda_{\mathrm{End}(F),\mathrm{nt}}}\varepsilon^{1/3}\lV Q\rV_{L^2(K,K^{\mathrm{End}(F)})}.
	\end{equation}
\end{prop}

\begin{pro}
By \eqref{7.38n}, we rewrite \eqref{7.101} as
\begin{equation}
		\bV(\mathrm{Id}-\tfrac{\mathrm{T}^*}{d-1})\nabla Q\bV_{L^2(K,K^{\mathrm{End}(F)})}\leqslant\varepsilon \lV Q\rV_{L^2(K,K^{\mathrm{End}(F)})}.
	\end{equation}
Multiplying on the left by $\big(\mathrm{Id}+\tfrac{\mathrm{T}^*}{d-1}+\cdots+(\tfrac{\mathrm{T}^*}{d-1})^{n-1}\big)$ and using \eqref{7.41n}, we obtain
	\begin{equation}\label{7.104}
		\bV\big(\mathrm{Id}-(\tfrac{\mathrm{T}^*}{d-1})^{n}\big)\nabla Q\bV_{L^2(K, K^{\mathrm{End}(F)})}\leqslant n\varepsilon(\tfrac{d}{d-1})^{1/2}\lV Q\rV_{L^2(K,K^{\mathrm{End}(F)})}.
	\end{equation}

For any $Q_{(k)}'\in L^2(K_{(k)}, K^{\mathrm{End}(F)})$,
\begin{equation}\label{7.105}
		\begin{split}
			&\big\langle\nabla Q,Q_{(k)}'\big\rangle_{L^2(K, K^{\mathrm{End}(F)})}\\
			&-\frac{1}{n}\big\langle\nabla Q,\big(\mathrm{Id}+\cdots+(\tfrac{\mathrm{T}}{d-1})^n\big)Q_{(k)}'\big\rangle_{L^2(K, K^{\mathrm{End}(F)})}\\
			&=\frac{1}{n}\Big(\big\langle (\mathrm{Id}-(\tfrac{\mathrm{T}^*}{d-1})^0)\nabla Q,Q_{(k)}'\big\rangle_{L^2(K, K^{\mathrm{End}(F)})}+\\
			&\quad\qquad\cdots+\big\langle \big(\mathrm{Id}-(\tfrac{\mathrm{T}^*}{d-1})^{n-1}\big)\nabla Q,Q_{(k)}'\big\rangle_{L^2(K, K^{\mathrm{End}(F)})}\Big).
		\end{split}
	\end{equation}
By \eqref{7.104}, the right hand side of \eqref{7.105} is bounded by
	\begin{equation}\label{7.106}
		\tfrac{(n-1)\varepsilon}{2}\lV Q\rV_{L^2(K, K^{\mathrm{End}(F)})}\bV Q_{(k)}'\bV_{L^2(K, K^{\mathrm{End}(F)})}.
	\end{equation}
Since $(\tfrac{\mathrm{T}}{d-1})^iQ_{(k)}'\perp (\tfrac{\mathrm{T}}{d-1})^jQ_{(k)}'$ if $i\neq j$, we have
	\begin{equation}\label{7.107}
		\begin{split}
\bV\big(\mathrm{Id}+\cdots+(\tfrac{\mathrm{T}}{d-1})^n\big)Q_{(k)}'\bV_{L^2(K, K^{\mathrm{End}(F)})}&\\
\leqslant \big(\tfrac{dn}{d-1}\big)^{1/2}\bV Q_{(k)}'\bV_{L^2(K, K^{\mathrm{End}(F)})}&.
		\end{split}
	\end{equation}
	
Combining \eqref{7.37}, \eqref{7.105}, \eqref{7.106}, \eqref{7.107}, we get
\begin{equation}\label{7.109}
\begin{split}
&\bv\big\langle\nabla Q,Q_{(k)}'\big\rangle_{L^2(K, K^{\mathrm{End}(F)})}\bv\\
&\leqslant\Big(\tfrac{2d}{((d-1)n)^{1/2}}+\tfrac{(n-1)\varepsilon}{2}\Big)\lV Q\rV_{L^2(K,K^{\mathrm{End}(F)})}\bV Q_{(k)}'\bV_{L^2(K,K^{\mathrm{End}(F)})}.
\end{split}
\end{equation}
Taking $n$ of order $\varepsilon^{-\frac{2}{3}}$ in \eqref{7.109}, we have
\begin{equation}\label{7.113n}
\begin{split}
			&\bv\big\langle\nabla Q,Q_{(k)}'\big\rangle_{L^2(K, K^{\mathrm{End}(F)})}\bv\\
			&\leqslant C\varepsilon^{1/3}\lV Q\rV_{L^2(K,K^{\mathrm{End}(F)})}\bV Q_{(k)}'\bV_{L^2(K,K^{\mathrm{End}(F)})}.
\end{split}
\end{equation}
Taking $Q_{(k)}'=\nabla Q_{(k-1)}$ in \eqref{7.113n}, we get
\begin{equation}
			\bV \nabla Q_{(k-1)}\bV_{L^2(K_{k},K^{\mathrm{End}(F)})}\leqslant C\varepsilon^{1/3}\lV Q\rV_{L^2(K,K^{\mathrm{End}(F)})}.
	\end{equation}
	Then \eqref{7.98} and \eqref{7.101} imply \eqref{7.103}.\qed	
\end{pro}

\begin{prop}\label{p7.21}
Let $Q\in L^2(K_{(k)},K^{\mathrm{End}(F)})\cap\ker(\mathrm{ad}_{\Delta^F})^\perp$. Then,
	\begin{equation}\label{7.112}
		\begin{split}
\bbV\frac{1}{t_0}\int_0^{t_0}e^{\sqrt{-1}t\mathrm{ad}_{\Delta^F}}Qdt&\bbV_{L^2(K,K^{\mathrm{End}(F)})}\\
\leqslant\frac{C_{k,\lambda_{\mathrm{End}(F),\mathrm{nt}}}}{{t_0}^{1/7}}\bV Q&\bV_{L^2(K_{(k)},K^{\mathrm{End}(F)})}.
		\end{split}
	\end{equation}
\end{prop}

\begin{pro}
By \eqref{7.36}, we can	use the spectral decomposition of $\mathrm{ad}_{\Delta^F}$ to write
\begin{equation}\label{7.113}
		Q=\mathbbm{1}_{[-\varepsilon,\varepsilon]}(\mathrm{ad}_{\Delta^F})(Q)+\mathbbm{1}_{(-\infty,\infty)\backslash[-\varepsilon,\varepsilon]}(\mathrm{ad}_{\Delta^F})(Q).
\end{equation}
For the second term, we have
	\begin{equation}\label{7.114}
		\begin{split}
\bbV\frac{1}{{t_0}}\int_0^{t_0}e^{\sqrt{-1}t\mathrm{ad}_{\Delta^F}}\mathbbm{1}_{(-\infty,\infty)\backslash[-\varepsilon,\varepsilon]}(\mathrm{ad}_{\Delta^F})(Q)dt&\bbV_{L^2(K,K^{\mathrm{End}(F)})}\\
\leqslant\frac{2}{\varepsilon {t_0}}\bV Q&\bV_{L^2(K,K^{\mathrm{End}(F)})}.
		\end{split}
	\end{equation}

We also have
	\begin{equation}
		\begin{split}
		\bV\mathrm{ad}_{\Delta^F}\mathbbm{1}_{[-\varepsilon,\varepsilon]}(\mathrm{ad}_{\Delta^F})(Q)&\bV_{L^2(K,K^{\mathrm{End}(F)})}\\
	&\leqslant\varepsilon\lV\mathbbm{1}_{[-\varepsilon,\varepsilon]}(\mathrm{ad}_{\Delta^F})(Q)\rV_{L^2(K,K^{\mathrm{End}(F)})},\\
	\mathbbm{1}_{[-\varepsilon,\varepsilon]}(\mathrm{ad}_{\Delta^F})(Q)&\perp_{L^2(K,K^{\mathrm{End}(F)})}\ker(\mathrm{ad}_{\Delta^F}).
		\end{split}
	\end{equation}
Thus, by \eqref{7.101} and \eqref{7.103},	
	\begin{equation}
		\begin{split}
	\bV\big(\mathbbm{1}_{[-\varepsilon,\varepsilon]}(\mathrm{ad}_{\Delta^F})(Q)\big)_{(k)}&\bV_{L^2(K,K^{\mathrm{End}(F)})}\\
	\leqslant C_{k,\lambda_{\mathrm{End}(F),\mathrm{nt}}}\varepsilon^{1/3}\bV Q&\bV_{L^2(K_{(k)},K^{\mathrm{End}(F)})}.
		\end{split}
	\end{equation}
Using the basic properties of projection-valued measures, we have
	\begin{equation}\label{7.117}
		\begin{split}
			&\lV\mathbbm{1}_{[-\varepsilon,\varepsilon]}(\mathrm{ad}_{\Delta^F})(Q)\rV^2_{L^2(K,K^{\mathrm{End}(F)})}\\
			&=\big\langle \mathbbm{1}_{[-\varepsilon,\varepsilon]}(\mathrm{ad}_{\Delta^F})(Q),Q\big\rangle_{L^2(K_{},K^{\mathrm{End}(F)})}\\
			&=\big\langle \big(\mathbbm{1}_{[-\varepsilon,\varepsilon]}(\mathrm{ad}_{\Delta^F})(Q)\big)_{(k)},Q\big\rangle_{L^2(K_{(k)},K^{\mathrm{End}(F)})}\\
			&\leqslant C_{k,\lambda_{\mathrm{End}(F),\mathrm{nt}}}\varepsilon^{1/3}\lV Q\rV_{L^2(K_{(k)},K^{\mathrm{End}(F)})}^2.
		\end{split}
	\end{equation}
	
Combining \eqref{7.113}, \eqref{7.114}, and \eqref{7.117}, we obtain
	\begin{equation}\label{7.118}
		\begin{split}
&\bbV\frac{1}{{t_0}}\int_0^{t_0}e^{\sqrt{-1}t\mathrm{ad}_{\Delta^F}}Qdt\bbV_{L^2(K,K^{\mathrm{End}(F)})}\\
&\leqslant \Big(C_{k,\lambda_{\mathrm{End}(F),\mathrm{nt}}}\varepsilon^{1/6}+\frac{2}{\varepsilon {t_0}}\Big)\lV Q\rV_{L^2(K_{(k)},K^{\mathrm{End}(F)})}.
		\end{split}
	\end{equation}
Taking $\varepsilon$ of order $t_0^{-6/7}$ in \eqref{7.118}, we get \eqref{7.112}.\qed	
\end{pro}

\subsection{Expander condition}\label{s7.6n}

Recall $\mathscr{F}(X)$ and $\mathscr{F}(X,\varepsilon)$ defined in \eqref{4.1n} and \eqref{n7.3}. We define the class of vector bundles $\mathscr{E}(X,\varepsilon)\subset \mathscr{F}(X)$ by
\begin{equation}\label{8.121}
	\begin{split}
		\mathscr{E}(X,\varepsilon)=\big\{F\in\mathscr{F}(X)\mid&\ \mathrm{End}(F)\in 	\mathscr{F}(X,\varepsilon), \\
		&d \text{ is a simple eigenvalue of } \Delta^{\mathrm{End}(F)}\big\}.
	\end{split}
\end{equation}
We say that $\mathrm{End}(F)$ is an expander, denoted by \textbf{[EXP]}, if $F\in \mathscr{E}(X,\varepsilon)$, or equivalently,
\begin{equation}\label{7.121}
	\begin{split}
		&\mathrm{spec}(\Delta^{\mathrm{End}(F)})\subset [-d,d-\varepsilon]\cup\{d\},\\
		&\ker(\Delta^{\mathrm{End}(F)}-d\mathrm{Id})=\mathrm{span}\{\mathrm{Id}_F\}.
	\end{split}
\end{equation}
Note that \eqref{7.121} is for $\mathrm{End}(F)$, rather than for $F$. If $F=\mathbb{C}$, the trivial line bundle, then $\mathrm{End}(F)=\mathbb{C}$, so the two notions \textbf{[GAP]} and \textbf{[EXP]} coincide.

For $Q\in L^2(K_{(k)},K^{\mathrm{End}(F)})$, we define its normalized average
\begin{equation}\label{7.125}
	\begin{split}
	\langle Q\rangle_{K^{\mathrm{End}(F)}}&=\frac{1}{\lv \pa B_{\mathbb{T}_d}(o,k)\rv}\langle Q,\mathrm{Id}_{(k)}\rangle_{L^2(K_{(k)},K^{\mathrm{End}(F)})}\\
	&=\frac{1}{\dim_{\mathbb{C}}F\cdot \lv \pa B_{\mathbb{T}_d}(o,k)\rv}\sum_{y\in D,y'\in\mathbb{T}_d}\mathrm{Tr}^{\mathbb{C}^\ell}\big[Q(y,y')\big],
	\end{split}
\end{equation}
where $\mathrm{Id}_{(k)}$ is given in \eqref{7.5nn}.

By \eqref{7.6nn}, \eqref{7.87}, and \eqref{7.121}, we get the following result.

\begin{prop}\label{p7.24}
Assume \textbf{[EXP]} \eqref{7.121}. Then
\begin{equation}\label{7.122}
	\ker(\mathrm{ad}_{\Delta^F})=\mathrm{span}\{\mathrm{Id}_{(k)}\mid k\in\mathbb{N}\}=\mathrm{span}\big\{\big(\Delta^{\mathrm{End}(F)}\big)^k\mid k\in\mathbb{N}\big\},
\end{equation}
and the projection in \eqref{7.99} becomes
\begin{equation}\label{7.126}
	P^{\ker(\mathrm{ad}_{\Delta^F})}Q=\langle Q\rangle_{K^{\mathrm{End}(F)}}\mathrm{Id}_{(k)}.
\end{equation}
\end{prop}

Finally, we note that much of the analysis in this section applies not only to the endomorphism bundle, but also to an arbitrary vector bundle. The main difference is that, for a general vector bundle $F$, the kernel of the gradient operator might be trivial, whereas for $\mathrm{End}(F)$ there is an automatic stationary section $\mathrm{Id}_{F}=\mathrm{Id}_{(0)}$.

\section{Geometric Quantization}\label{s8n}

In this section, we introduce Berezin-Toeplitz quantization. In \S\,\ref{s8.1n}, we define Berezin-Toeplitz quantization and establish some basic properties. In \S\,\ref{s8.2n}, we derive a trace formula for Berezin-Toeplitz operators. In \S\,\ref{s8.3n}, we prove off-diagonal exponential decay for Berezin-Toeplitz operators.

Our main references for this section are Ma-Marinescu \cite{MR2339952,MR2393271}.

\subsection{Berezin-Toeplitz quantization}\label{s8.1n}

Recall the projection operator $P_p$ defined in \eqref{3.3n}. 

The \emph{Berezin-Toeplitz quantization} is the family of maps
\begin{equation}\label{8.1}
	\begin{split}
&T_{\cdot,p}\colon C^\infty(N)\to \mathrm{End}\big(L^2(N,L^p)\big),\\
&T_{f,p}=P_pfP_p.
	\end{split}
\end{equation}
Equivalently, for any $s,s'\in H^{(0,0)}(N,L^p)$, we have
\begin{equation}\label{c3}
	\langle T_{f,p}s,s'\rangle_{L^2(N,L^p)}=\int_{N}f(z)\langle s(z),s'(z)\rangle_{h^{L^p}}dv_{N}(z).
\end{equation}

The operator $T_{f,p}$ has a smooth Schwartz kernel $T_{f,p}(z,z')\in L_z^p\otimes (L_{z'}^*)^p$ given by
\begin{equation}\label{8.3}
	T_{f,p}(z,z')=\int_NP_p(z,z'')f(z'')P_p(z'',z')dv_N(z'').
\end{equation}
Indeed, by \eqref{3.3},
\begin{equation}
	\begin{split}
		(T_{f,p}s)(z)=&(P_pfP_ps)(z)\\
		=&\int_N\int_NP_p(z,z'')f(z'')P_p(z'',z')s(z')dv_N(z'')dv_N(z').
	\end{split}
\end{equation}

\begin{lemma}
Let $f\in C^\infty(N)$. Then for any $p\in\mathbb{N}$ and $s\in L^2(N,L^p)$,
\begin{equation}\label{c4}
	\lV T_{f,p}s\rV_{L^2(N,L^p)}\leqslant \lV f\rV_{C^0(N)}\lV s\rV_{L^2(N,L^p)}.
\end{equation}
\end{lemma}

\begin{pro}
This follows directly from \eqref{c3}.\qed
\end{pro}

\subsection{Trace formula}\label{s8.2n}

On the diagonal, we have $T_{f,p}(z,z)\in \mathrm{End}(L_z^p)=\mathbb{C}$, and by \eqref{8.3},
\begin{equation}\label{8.6}
	\begin{split}
\tro^{H^{(0,0)}(N,L^p)}[T_{f,p}]&=\int_NT_{f,p}(z,z)dv_N(z)\\
&=\int_N\int_Nf(z'')P_p(z'',z)P_p(z,z'')dv_N(z)dv_N(z'')\\
&=\int_Nf(z'')P_p(z'',z'')dv_N(z''),
	\end{split}
\end{equation}
where we use the fact that $P_p$ is a projection, so $P_p^2=P_p$, and hence
\begin{equation}
\int_NP_p(z,z'')P_p(z'',z')dv_N(z'')=P_p(z,z').
\end{equation}

\begin{prop}\label{p8.2n}
There is $C>0$ such that for any $p\in\mathbb{N}$ and $f\in C^\infty(N)$,
\begin{equation}\label{8.8}
	\begin{split}
		\Bv\frac{1}{\dim_{\mathbb{C}} H^{(0,0)}(N,L^p)}&\tro^{H^{(0,0)}(N,L^p)}[T_{f,p}]\\
		&-\frac{1}{\mathrm{Vol}(N)}\int_Nfdv_N\Bv_{\mathbb{C}}\leqslant \frac{C}{p}\int_{N}\lv f\rv_{\mathbb{C}}dv_{N},\\
		\frac{1}{\dim_{\mathbb{C}} H^{(0,0)}(N,L^p)}&\lV T_{f,p}\rV_{\mathrm{HS}(H^{(0,0)}(N,L^p))}^2\leqslant C\int_{N}\lv f\rv_{\mathbb{C}}^2dv_{N}.
	\end{split}
\end{equation}
\end{prop}

\begin{pro}
The first inequality follows directly from \eqref{3.6n} and \eqref{8.6}.
	
For the second inequality, we use the elementary estimate
	\begin{equation}
		\lV P_pfP_ps\rV_{L^2(N,L^p)}^2\leqslant\lV fP_ps\rV_{L^2(N,L^p)}^2.
	\end{equation}
	Rewriting this in terms of \eqref{c3}, we obtain
	\begin{equation}
	\big\langle T_{f,p}^*T_{f,p}s,s\big\rangle_{L^2(N,L^p)}\leqslant \big\langle T_{\lv f\rv_{\mathbb{C}}^2,p}s,s\big\rangle_{L^2(N,L^p)}.
\end{equation}
The second inequality then follows from the first.\qed
\end{pro}

Note that the right hand side of the second inequality in \eqref{8.8} is uniform for $p\in\mathbb{N}$, with no $O(p^{-1})$ remainder term. By contrast, if we apply the product formula \cite[Theorem 7.4.1]{MR3368102} directly, an additional $O(p^{-1})$ remainder term appears.

\subsection{Off-diagonal exponential decay}\label{s8.3n}

Similarly to Proposition \ref{p3.1}, we have the following estimate.

\begin{prop}\label{8.3n}
There exist $C,c>0$ such that for any $f\in C^0(N)$, $p\in\mathbb{N}$, and $z,z'\in N$, we have
	\begin{equation}\label{8.11}
		\lV T_{f,p}(z,z')\rV_{L_z^p\otimes (L_{z'}^*)^p}\leqslant C\lV f\rV_{C^0(N)}p^{n}e^{-c\sqrt{p}d_N(z,z')},
	\end{equation}
	where $d_N(\cdot,\cdot)$ denotes the Riemannian geodesic distance on $N$.
\end{prop}

\begin{pro}
By \eqref{3.5} and \eqref{8.3}, we have
\begin{equation}\label{8.12}
	\begin{split}
		&\lV T_{f,p}(z,z')\rV_{L_z^p\otimes (L_{z'}^*)^p}\\
		&\leqslant Cp^{2n}\lV f\rV_{C^0(N)} \int_N e^{-c\sqrt{p}(d_N(z,z'')+d_N(z'',z'))}dv_N(z'')\\
		&\leqslant \big(Cp^{n}e^{-c\sqrt{p}d_N(z,z')/2}\lV f\rV_{C^0(N)}\big)\\
		&\quad\cdot \Big(p^n\int_N e^{-c\sqrt{p}(d_N(z,z'')+d_N(z'',z'))/2}dv_N(z'')\Big).
	\end{split}
\end{equation}
Thus, it suffices to show that the final integral in \eqref{8.12} is bounded uniformly for $p\in\mathbb{N}$.

Let $\mathrm{inj}_N>0$ denote the injective radius of $N$. Let $B_N(z,\mathrm{inj}_N/2)$ be the geodesic ball in $N$ centered at $z$ with radius $\mathrm{inj}_N/2$, and let $B_N(z,\mathrm{inj}_N/2)^c$ denote its complement. We split the integral in \eqref{8.12} accordingly. 

First, we have
\begin{equation}\label{8.13}
	\begin{split}
		&p^{n}\int_{B_N(z,\mathrm{inj}_N/2)^c} e^{-c\sqrt{p}d_N(z,z'')/2}dv_N(z'')\\
		&\leqslant p^{n}\int_{N} e^{-c\sqrt{p}\mathrm{inj}_N/4}dv_N(z'')=p^{n}\mathrm{Vol}(N) e^{-c\sqrt{p}\mathrm{inj}_N/4}.
	\end{split}
\end{equation}
Next, using geodesic normal coordinates centered at $z$, we obtain
\begin{equation}\label{8.14}
	\begin{split}
		&p^{n}\int_{B_N(z,\mathrm{inj}_N/2)} e^{-c\sqrt{p}d_N(z,z'')/2}dv_N(z'')\\
		&=p^{n}\int_{B_{T_zN}(0,\mathrm{inj}_N/2)} e^{-c\sqrt{p}Z''/2}\big(\det g_{ij}(Z'')\big)^{1/2}dZ''\\
		&=\int_{B_{T_zN}(0,\mathrm{inj}_N/2)} e^{-c\sqrt{p}Z''/2}\big(\det g_{ij}(Z'')\big)^{1/2}d(\sqrt{p}Z'')\\
		&\leqslant\Big(\sup_{Z''\in B_N(z,\mathrm{inj}_N/2) }\big(\det g_{ij}(Z'')\big)^{1/2}\Big)\int_{\mathbb{R}^{2n}}e^{-cZ''/2}dZ''.
	\end{split}
\end{equation}
Here $B_{T_zN}(0,\mathrm{inj}_N/2)$ is the ball in the tangent space $T_zN$, centered at $0$ with radius $\mathrm{inj}_N/2$, the functions $g_{ij}$ denote the local components of the Riemannian metric, and the factor $p^n$ is absorbed by the change of variables.

Combining \eqref{8.13} and \eqref{8.14}, we complete the proof.\qed
\end{pro}

The product $T_{f,p}T_{f',p}$ of two Berezin-Toeplitz operators is again an operator with smooth kernel
\begin{equation}
	(T_{f,p}T_{f',p})(z,z')=\int_NT_{f,p}(z,z'')T_{f',p}(z'',z')dv_N(z'').
\end{equation}
Repeating the proof of Proposition \ref{8.3n}, we obtain the following estimate.

\begin{prop}\label{8.4n}
There exist $C,c>0$ such that for any any $f,f'\in C^0(N)$ $p\in\mathbb{N}$ and $z,z'\in N$, we have
\begin{equation}\label{8.16}
	\lV (T_{f,p}T_{f',p})(z,z')\rV_{L_z^p\otimes (L_{z'}^*)^p}\leqslant C\lV f\rV_{C^0(N)}\lV f'\rV_{C^0(N)}p^{n}e^{-c\sqrt{p}d_N(z,z')}.
\end{equation}
\end{prop}

\section{Geometric Quantization of Kernel Functions}\label{s9n}

In this section, we combine \S\S\,\ref{s7n} and \ref{s8n} to discuss mixed quantization. In \S\,\ref{s9.1}, we extend some of the results of \S\,\ref{s7n} to an infinite dimensional setting. In \S\,\ref{s9.2}, we introduce mixed quantization and establish some useful properties. In \S\,\ref{s9.3}, we discuss the spin-scale-frequency correspondence and the pseudolocality of mixed quantization operators.

For mixed quantization, we refer to Ma-Ma \cite[\S\,4]{MR4808253} and Ben Ovadia-Ma-Rodriguez-Hertz \cite[\S\,3]{ovadia2025mixedquantizationpartialhyperbolicity} for the manifold case. For pseudolocality in semiclassical analysis on manifolds, we refer to Zworski \cite[\S\,8.4.5]{MR2952218}.

\subsection{Kernel functions and expander condition}\label{s9.1}

Throughout this section, we work in the setting of \eqref{4.4n} and \eqref{4.5n}.

Similarly to \eqref{2.19n}, \eqref{2.20n}, \eqref{7.58}, and \eqref{7.59}, we define
\begin{equation}\label{9.1}
	\begin{split}
		X^N=&\Gamma\backslash\big(\mathbb{T}_d\times N\big)\\
		=&\big\{(y,z)\in \mathbb{T}_d\times N\big\}\big/(y,z)\sim(\gamma y,\gamma z) \text{ for any }\gamma\in\Gamma.
	\end{split}
\end{equation}
Then the space of functions and the Laplacian are given by
\begin{equation}\label{9.2}
	\begin{split}
	&C^\infty(X^N)\\
	&=C^\infty_\Gamma(\mathbb{T}_d\times N)\\
		&=\big\{\mathscr{U}\mid \mathscr{U}(y,z)=\mathscr{U}(\gamma y,\gamma z) \text{ for any }(y,z)\in \mathbb{T}_d\times N \text{ and }\gamma\in\Gamma\big\},\\
		&(\Delta^{X^N}\mathscr{U})(y,z)=\sum_{y'\sim y}\mathscr{U}(y',z).
	\end{split}
\end{equation}

Similarly to \eqref{7.1n} and \eqref{7.2n}, we define
 \begin{equation}
 	\begin{split}
K^N=&\Gamma\backslash\big(\mathbb{T}_d^2\times N\big)\\
=&\big\{(y,y',z)\in \mathbb{T}_d^2\times N\big\}\big/(y,y',z)\sim(\gamma y,\gamma y',\gamma z) \text{ for any }\gamma\in\Gamma.
 	\end{split}
 \end{equation}
 Then the space of functions is given by
 \begin{equation}
\begin{split}
	&C^\infty(K^N)\\
	&=C^\infty_\Gamma\big(\mathbb{T}_d^2\times N\big)\\
&=\big\{\mathscr{Q}\mid \mathscr{Q}(y,y',z)=\mathscr{Q}(\gamma y,\gamma y',\gamma z) \text{ for any }(y,y',z)\in \mathbb{T}_d^2\times N \text{ and }\gamma\in\Gamma\big\},
\end{split}
 \end{equation}
 and we call its elements \emph{kernel functions}. We define $K^N_{(k)},K^N_{(\leqslant k)}$, and $\mathscr{Q}_{(k)}$ analogously to \eqref{7.4} and \eqref{7.5nnn}.

Similarly to \eqref{7.5nn}, define $\mathbbm{1}_{(k)}\in C^\infty(K_{(k)}^N)$ by
 \begin{equation}\label{9.5}
\mathbbm{1}_{(k)}(y,y',z)=\mathbbm{1}_{\{d_{\mathbb{T}_d}(y,y')=k\}}(y,y'),
 \end{equation}
 which depends only on the distance in the base $K$.
 
Similarly to \eqref{7.5n}, let $\mathscr{Q}\in C^\infty(K^N)$ act on $\mathscr{U}\in L^2(X^N)$ by
 \begin{equation}
 	\begin{split}
 		(\mathscr{QU})(y,z)=&\sum_{y'\in\mathbb{T}_d}\mathscr{Q}(y,y',z)\mathscr{U}(y',z)\\
 		=&\sum_{y'\in D,\gamma\in\Gamma}\mathscr{Q}(y,\gamma y',z)\mathscr{U}(\gamma y',z)\\
 		=&\sum_{y'\in D,\gamma\in\Gamma}\mathscr{Q}(y,\gamma y')\mathscr{U}\big(y',\gamma^{-1}z\big).
 	\end{split}
 \end{equation}

Similarly to \eqref{2.21n} and \eqref{7.7n}, we define the $L^2$-norm on $L^2(X^N)$, the normalized $L^{2}$- and the normalized mixed $L^{\infty,2}$- norms on $C^\infty(K^N)$ by
\begin{equation}\label{9.7}
	\begin{split}
		&\lV \mathscr{U}\rV_{L^2(X^N)}^2=\sum_{y\in D}\lV \mathscr{U}(y,\cdot)\rV_{L^2(N)}^2,\\
		&\lV\mathscr{Q}\rV_{L^2(K^N)}^2=\frac{1}{\lv D\rv}\sum_{y\in D,y'\in \mathbb{T}_d}\lV \mathscr{Q}(y,y',\cdot)\rV_{L^2(N)}^2,\\
		&\lV\mathscr{Q}\rV_{L^{\infty,2}(K^N)}^2=\max_{y\in D,y'\in \mathbb{T}_d}\lV \mathscr{Q}(y,y',\cdot)\rV_{L^2(N)}^2.
	\end{split}
\end{equation}
Here $\lV\cdot\rV_{L^{\infty,2}(K^N)}$ denotes the mixed norm with $L^\infty$ in the base variables $(y,y')$ and $L^2$ in the fibre variable $z$.

Similarly to Lemma \ref{L7.2}, we have the following estimate.
\begin{lemma}
	Let $\mathscr{Q}\in C^\infty(K_{(\leqslant k)}^N)$. Then
	\begin{equation}\label{9.8}
		\lV\mathscr{Q}\rV_{L^2(K^N)}^2\leqslant	\lv B_{ \mathbb{T}_d}(o,k)\rv\lV\mathscr{Q}\rV_{L^{\infty,2}(K^N)}^2.
	\end{equation}
\end{lemma}

Similarly to \eqref{7.27n} and \eqref{7.33n}, we define
\begin{equation}\label{9.9}
	\begin{split}
		&\mathrm{ad}_{\Delta^{X^N}}\colon C^{\infty}(K_{(k)}^N) \to C^{\infty}(K_{(k-1)}^N\cup K_{(k+1)}^N),\\
		&\mathrm{ad}_{\Delta^{X^N}}(\mathscr{Q})(y,y',z)=\sum_{y''\sim y}\mathscr{Q}(y'',y',z)-\sum_{y''\sim y'}\mathscr{Q}(y,y'',z),\\
			&\nabla\colon C^{\infty}(K_{(k)}^N) \to C^{\infty}(K_{(k+1)}^N),\\
		&\nabla(\mathscr{Q})(y_{k+1},\cdots,y_{0},z)=Q(y_{k+1},\cdots,y_{1},z)-Q(y_k,\cdots,y_0,z).
	\end{split}
\end{equation}

Similarly to \eqref{7.125}, for $\mathscr{Q}\in L^2(K_{(k)}^N)$, we define its normalized average
\begin{equation}\label{9.12}
	\begin{split}
		\langle \mathscr{Q}\rangle_{K_{}^N}&=\frac{1}{\lv \pa B_{\mathbb{T}_d}(o,k)\rv\cdot\mathrm{Vol}(N)}\big\langle \mathscr{Q},\mathbbm{1}_{(k)}\big\rangle_{L^2(K_{(k)}^N)}\\
		&=\frac{1}{\lv D\rv\cdot \lv \pa B_{\mathbb{T}_d}(o,k)\rv}\sum_{y\in D,y'\in\mathbb{T}_d}\frac{1}{\mathrm{Vol}(N)}\int_N\mathscr{Q}(y,y',z)dv_N(z),
	\end{split}
\end{equation}
where $\mathbbm{1}_{k}$ is given in \eqref{9.5}.

Similarly to Proposition \ref{p7.24}, we have the following result.
\begin{prop}
Assume \textbf{[DEN]} \eqref{4.9n}. Then
	\begin{equation}\label{9.10}
			\ker(\mathrm{ad}_{\Delta^{X^N}})=\mathrm{span}\{\mathbbm{1}_{(k)}\mid k\in\mathbb{N}\},
	\end{equation}
	and we can express the orthogonal projection
	\begin{equation}\label{9.13}
		\begin{split}
			P^{\ker(\mathrm{ad}_{\Delta^{X^N}})}&\colon L^2(K_{}^N)\to \ker(\mathrm{ad}_{\Delta^{X^N}}),\\
			P^{\ker(\mathrm{ad}_{\Delta^{X^N}})}&\mathscr{Q}=\langle \mathscr{Q}\rangle_{K_{}^N}\mathbbm{1}_{(k)}.
		\end{split}
	\end{equation}
\end{prop}

\begin{pro}
It remains to studying $\ker(\nabla)$, kernel functions $\mathscr{Q}_{(k)}$ of the form
\begin{equation}
\begin{split}
\mathscr{Q}_{(k)}(y,y',z)&=f(z),\\
\mathscr{Q}_{(k)}(\gamma y,\gamma y',\gamma z)&=\mathscr{Q}_{(k)}(y,y',z)\text{ for any }\gamma\in\Gamma,
\end{split}
\end{equation}
for some $f\in L^2(N)$. These force $f(z)=f(\gamma z)$. By \textbf{[DEN]} \eqref{4.9n}, the only $\Gamma$-invariant functions on $N$ are constant functions.\qed
\end{pro}

Similarly to Proposition \ref{p7.21}, but in a nonquantitative form, we have the following result.

\begin{prop}\label{p9.3}
Assume \textbf{[DEN]} \eqref{4.9n}. Let $\mathscr{Q}\in L^2(K^N)\cap\ker(\mathrm{ad}_{\Delta^{X^N}})^\perp$. Then
	\begin{equation}\label{9.14}
		\lim_{{t_0}\to\infty}\bbV\frac{1}{{t_0}}\int_0^{t_0}e^{\sqrt{-1}t\mathrm{ad}_{\Delta^{X^N}}}\mathscr{Q}dt\bbV_{L^2(K^N)}=0.
	\end{equation}
\end{prop}

\begin{pro}
This follows from the von Neumann mean ergodic theorem.\qed
\end{pro}

No quantitative result is obtained here because infinite dimensional operators may have continuous spectrum. This motivates us to introduce the following condition, analogous to \eqref{7.121}.

We say that $X^N$ is an expander, denoted by \textbf{[EXP]}, if there is $\varepsilon>0$ such that
\begin{equation}\label{9.15}
\begin{split}
&\mathrm{spec}(\Delta^{X^N})\subseteq [-d,d-\varepsilon]\cup \{d\},\\
&\ker\big(\Delta^{X^N}-d\mathrm{Id}\big)=\mathrm{span}\{\mathbbm{1}_{X^N}\}.
\end{split}
\end{equation}
In other words, $d$ is a simple discrete eigenvalue of $\Delta^{X^N}$.

In the Borel-Weil-Bott case $(N,L)=(\mathcal{O}_\alpha,L_\alpha)$ as in \eqref{3.18n} and \eqref{4.10n}, we may formulate a spectral gap condition uniformly for $\alpha\in \Lambda^*\cap\overline{C^+}$, again denoted by \textbf{[EXP]}, by requiring that
\begin{equation}\label{9.15n}
	\begin{split}
		&\mathrm{spec}(\Delta^{X^G})\subseteq [-d,d-\varepsilon]\cup \{d\},\\
		&\ker\big(\Delta^{X^G}-d\mathrm{Id}\big)=\mathrm{span}\{\mathbbm{1}_{X^G}\},
	\end{split}
\end{equation}
where we replace $N$ with $G$ in \eqref{9.1}, \eqref{9.2}, and \eqref{9.15}.

In general, \textbf{[EXP]} in \eqref{9.15} and \eqref{9.15n} is much stronger than \textbf{[DEN]} in \eqref{4.9n} and \eqref{4.10n}. 

We define the infinite dimensional analogue of $\mathrm{C}_{1,0}, \mathrm{O}$, and $\mathrm{B}_1$ similarly to \eqref{7.17n}, \eqref{7.19n}, and \eqref{7.43n} by
\begin{equation}
	\begin{split}
		&\mathrm{C}_{1,0}(\mathscr{Q}_{(0)})(y_1,y_0,z)=\mathscr{Q}_{(0)}(y_0,z),\\
		&\mathrm{O}(\mathscr{Q}_{(1)})(y_1,y_0,z)=\mathscr{Q}_{(1)}(y_0,y_1,z),\\
		&\mathrm{B}_1(\mathscr{Q}_{(1)})(y_1,y_0,z)=\sum_{y_{-1}\sim y_0,y_{-1}\neq y_1}\mathscr{Q}_{(1)}(y_{0},y_{-1},z).\\
	\end{split}
\end{equation}

Similarly to Proposition \ref{p7.16}, the Jordan structure of $\mathrm{B}_1$ is as follows.
\begin{prop}\label{p9.4}
Assume \textbf{[EXP]} \eqref{9.15}. The operator $\mathrm{B}_1$ acting on $L^2(K_{(1)}^N)$ is unitarily equivalent to a block diagonal operator of the form
	\begin{equation}\label{9.18n}
		\mathrm{diag}\bigg((d-1)\mathrm{Id},-(d-1)\mathrm{Id},\mathrm{Id},-\mathrm{Id},\bigg(\begin{smallmatrix}
			\theta(\Delta^{X^N})& a(\Delta^{X^N})\\
			0&\theta(\Delta^{X^N})
		\end{smallmatrix}\bigg)\bigg).
	\end{equation}
Here $\theta(\Delta^{X^N}), \theta'(\Delta^{X^N})$, and $a(\Delta^{X^N})$ are defined by functional calculus of $\Delta^{X^N}$, in which $\theta(\lambda)$ and $\theta'(\lambda)$ are the two roots of
	\begin{equation}
		\theta^2-\lambda\theta+d-1=0, 
	\end{equation}
	while $a(\lambda)$ satisfies
	\begin{equation}
a(\lambda)=\begin{cases}
			\big(d^2-\lambda^2\big)^{1/2}, &\text{if } 	\lambda^2\geqslant4(d-1),\\
			(d-2),&\text{if } 	\lambda^2\leqslant4(d-1).
		\end{cases}
	\end{equation}
\end{prop}

\begin{pro}
	We follow the proof of Proposition \ref{p7.16}. Although the dimension counting argument is no longer available in the present infinite dimensional setting, the rest of the proof can be adapted without essential modifications. We summarize the main steps.
	
	First, we have the orthogonal decomposition
	\begin{equation}\label{9.19n}
		\begin{split}
			L^2(K_{(1)}^N)&=L^2_+(K_{(1)}^N)\oplus L^2_-(K_{(1)}^N),\\
			L^2_\pm(K_{(1)}^N)&=\big\{\mathscr{Q}_{(1)}\mid \mathscr{Q}_{(1)}(y_0,y_1,z)=\pm \mathscr{Q}_{(1)}(y_1,y_0,z) \text{ for any }(y_0,y_1,z)\in \mathbb{T}_d^2\times N\big\}.
		\end{split}
	\end{equation}
	Define the star subspaces $S_\pm(K_{(1)}^N)\subset L^2_\pm(K_{(1)}^N)$ by
	\begin{equation}\label{9.22n}
		\begin{split}
			S_\pm(K_{(1)}^N)&=\overline{\mathrm{span}\big\{\mathscr{Q}_{(1),y,f}^\pm\mid y\in D, f\in L^2(N)\big\}},\\ \mathscr{Q}_{(1),y,f}^\pm(y_1,y_0,z)&=\begin{cases}f(\gamma z),&\text{if }y_0=\gamma y\text{ for }\gamma\in \Gamma,\\
				\pm f(\gamma z),&\text{if }y_1=\gamma y\text{ for }\gamma\in \Gamma,\\
				0,&\text{otherwise},
			\end{cases}
		\end{split}
	\end{equation}
where $D\subset\mathbb T_d$ is the fundamental domain from \eqref{2.19n}. Then, we have
	\begin{equation}
		\mathrm{B}_1(\mathscr{Q}_{(1)}^\pm)=\mp \mathscr{Q}_{(1)}^\pm,\quad\text{if }\mathscr{Q}_{(1)}^\pm\in  L^2_\pm(K_{(1)}^N)\cap S_\pm(K_{(1)}^N)^\perp.
	\end{equation}
	
	Next, we analyze the restriction of $\mathrm{B}_1$ to the star spaces \eqref{9.22n}. 
	
	Similarly to \eqref{7.76} and \eqref{7.76n}, we have
	\begin{equation}\label{9.23}
		\begin{split}
			\big\langle \mathrm{O}\mathrm{C}_{1,0}\mathscr{Q}_{(0)},\mathrm{O}\mathrm{C}_{1,0}\mathscr{Q}'_{(0)} \big\rangle_{L^2(K_{(1)}^N)}&=d\big\langle \mathscr{Q}_{(0)},\mathscr{Q}_{(0)}'\big\rangle_{L^2(K_{(0)}^N)},\\
			\big\langle \mathrm{C}_{1,0}\mathscr{Q}_{(0)},\mathrm{O}(\mathrm{C}_{1,0}\mathscr{Q}_{(0)}') \big\rangle_{L^2(K_{(1)}^N)}&=\big\langle \mathscr{Q}_{(0)},\Delta^{X^N}\mathscr{Q}_{(0)}'\big\rangle_{L^2(K_{(0)}^N)}.
		\end{split}
	\end{equation}	
	
	We replace \eqref{7.60} with the spectral decomposition
	\begin{equation}
		L^2(K_{(0)}^N)=\mathrm{im}\big(\mathbbm{1}_{\{\pm d\}}(\Delta^{X^N})\big)\oplus\mathrm{im}\big(\mathbbm{1}_{(-d,d-\varepsilon]}(\Delta^{X^N})\big).
	\end{equation}
	Then by \eqref{9.23} and arguing as in \eqref{7.77}, \eqref{7.86n}, \eqref{7.88}, and \eqref{7.89}, we obtain the orthogonal decomposition
	\begin{equation}\label{9.25}
		\begin{split}
			S_+(K_{(1)}^N)\oplus S_-(K_{(1)}^N)=&\mathrm{C}_{1,0}\mathrm{im}\big(\mathbbm{1}_{\{\pm d\}}(\Delta^{X^N})\big)\\
			&\oplus\mathrm{C}_{1,0}\mathrm{im}\big(\mathbbm{1}_{(-d,d-\varepsilon]}(\Delta^{X^N})\big)\\
			&\oplus \big(\mathrm{O}\mathrm{C}_{1,0}-\tfrac{\mathrm{C}_{1,0}\Delta^{X^N}}{d}\big)\mathrm{im}\big(\mathbbm{1}_{(-d,d-\varepsilon]}(\Delta^{X^N})\big).
		\end{split}
	\end{equation}
	Moreover, if we identify the last two summands in \eqref{9.25} with two copies of $\mathrm{im}\big(\mathbbm{1}_{(-d,d-\varepsilon]}(\Delta^{X^N})\big)$, then $\mathrm{B}_1$ takes the form
	\begin{equation}
		\mathrm{B}_1=\begin{pmatrix}
			\tfrac{d-1}{d}\Delta^{X^N}& (d-1)\big(1-\big(\tfrac{\Delta^{X^N}}{d}\big)^2\big)^{1/2}\\
			-\big(1-\big(\tfrac{\Delta^{X^N}}{d}\big)^2\big)^{1/2}&\frac{\Delta^{X^N}}{d}
		\end{pmatrix}.
	\end{equation}
	Since the entries of $\mathrm{B}_1$ are obtained by functional calculus from $\Delta^{X^N}$, we may formally regard $\Delta^{X^N}$ as a real number and apply the Schur decomposition to obtain \eqref{9.18n}.\qed
\end{pro}

Similarly to Proposition \ref{p7.21}, we then have the following quantitative form of Proposition \ref{p9.3}.
\begin{prop}
Assume \textbf{[EXP]} \eqref{9.15}. Let $\mathscr{Q}\in L^2(K_{(k)}^N)\cap\ker(\mathrm{ad}_{\Delta^{X^N}})^\perp$. Then,
	\begin{equation}\label{9.16}
			\bbV\frac{1}{{t_0}}\int_0^{t_0}e^{\sqrt{-1}t\mathrm{ad}_{\Delta^{X^N}}}\mathscr{Q}dt\bbV_{L^2(K_{}^N)}\leqslant\frac{C_{k,\varepsilon}}{{t_0}^{1/7}}\bV\mathscr{Q}\bV_{L^2(K_{(k)}^N)}.
	\end{equation}
\end{prop}

\begin{pro}
With \eqref{9.18n} in hand, the arguments leading to \eqref{7.87}, \eqref{7.94}, \eqref{7.98}, \eqref{7.103}, and \eqref{7.112} extend to the infinite dimensional setting without essential modifications.\qed
\end{pro}

\subsection{Mixed quantization}\label{s9.2}

By taking fibrewise Berezin-Toeplitz quantization as in \eqref{8.1}, we obtain the map
\begin{equation}\label{9.17}
	\begin{split}
&T_{\cdot,p}\colon C^\infty(K_{(k)}^N)\to C^\infty\big(K_{(k)},K^{\mathrm{End}(F_p)}\big),\\
&T_{\mathscr{Q},p}(y,y')=T_{\mathscr{Q}(y,y',\cdot),p},
\end{split}
\end{equation}
where $\mathscr{Q}(y,y',\cdot)\in C^\infty(N)$. Since kernel operators serve as a discrete analogue of semiclassical quantization on manifolds, the fibrewise geometric quantization $T_{\mathscr{Q},p}$ of a kernel function $\mathscr{Q}$ provides a natural discrete model for \emph{mixed quantization}.

The following result describes the relation between $\mathrm{ad}_{\Delta^{X^N}\mathscr{Q}}$ and $\mathrm{ad}_{\Delta^{F_p}}$.
\begin{lemma}
	Let $\mathscr{Q}\in C^\infty(K^N)$. Then for any $p\in\mathbb{N}$,
	\begin{equation}\label{9.18}
		\mathrm{ad}_{\Delta^{F_p}}T_{\mathscr{Q},p}=T_{\mathrm{ad}_{\Delta^{X^N}\mathscr{Q}},p}.
	\end{equation}
\end{lemma}

\begin{pro}
This is immediate from \eqref{7.27n}, \eqref{9.9}, and \eqref{9.17}.\qed
\end{pro}

\begin{lemma}\label{l9.6n}
	There is $C>0$ such that for any $p\in\mathbb{N}$ and $\mathscr{Q}\in C^\infty(K^N)$,
	\begin{equation}\label{9.19}
		\begin{split}
\lV T_{\mathscr{Q},p}\rV_{L^2(K,K^{\mathrm{End}(F_p)})}^2&\leqslant C\lV \mathscr{Q}\rV_{L^2(K^N)}^2,\\
\lV T_{\mathscr{Q},p}\rV_{L^\infty(K,K^{\mathrm{End}(F_p)})}^2&\leqslant C\lV \mathscr{Q}\rV_{L^{\infty,2}(K^N)}^2.
		\end{split}
	\end{equation}
\end{lemma}

\begin{pro}
This is a direct consequence of \eqref{7.7n}, \eqref{8.8}, and \eqref{9.7}.\qed
\end{pro}

The right hand side of \eqref{9.19} is uniform in $p\in\mathbb{N}$, with no $O(p^{-1})$ remainder term. This will be a key point later in large scale QE uniform in the spin parameter. It is also the key point emphasized in \cite[\S\,5]{MR4808253}.

The following result will be crucial later in large spin QE.
\begin{prop}\label{p9.7}
Assume \textbf{[FREE]} \eqref{4.8n}. Let $\mathscr{Q}\in C^\infty(K_{(k)}^N)$. As $p\to\infty$,
\begin{equation}\label{9.20}
		\lV T_{\mathscr{Q},p}\rV_{\mathrm{HS}(X,F_p)}^2=\lV T_{\mathscr{Q},p}\rV_{L^2(K_{(k)},K^{\mathrm{End}(F_p)})}^2+o(1).
\end{equation}
\end{prop}

\begin{pro}
	
By \eqref{7.10n}, it suffices to show that
	\begin{equation}\label{9.21}
		\lim_{p\to\infty}\sum_{\substack{y,y'\in D,\  \gamma,\gamma'\in\Gamma\\ \gamma\neq \gamma'}}\frac{\mathrm{Tr}^{H^{(0,0)}(N,L^p)}}{\dim_{\mathbb{C}}F_p}\Big[(\gamma')^{-1}T_{\overline{\mathscr{Q}},p}(y,\gamma'y')T_{\mathscr{Q},p}(y,\gamma y')\gamma\Big]=0.
	\end{equation}
	Since \eqref{9.21} is a finite sum, the proof reduces to establishing
	\begin{equation}\label{9.22}
		\lim_{p\to \infty}\frac{\mathrm{Tr}^{H^{(0,0)}(N,L^p)}}{\dim_{\mathbb{C}}H^{(0,0)}(N,L^p)}\Big[(\gamma')^{-1}T_{\overline{\mathscr{Q}},p}(y,\gamma'y')T_{\mathscr{Q},p}(y,\gamma y')\gamma\Big]=0.
	\end{equation}
Similarly to \eqref{3.10n}, the operator $(\gamma')^{-1}T_{\overline{\mathscr{Q}},p}(y,\gamma'y')T_{\mathscr{Q},p}(y,\gamma y')\gamma$ has kernel
	\begin{equation}
		\begin{split}
			&\Big((\gamma')^{-1}T_{\overline{\mathscr{Q}}(y,\gamma'y'),p}T_{\mathscr{Q}(y,\gamma y'),p}\gamma\Big)(z,z')\\
			&=(\gamma')^{-1}\Big((T_{\overline{\mathscr{Q}}(y,\gamma'y'),p}T_{\mathscr{Q}(y,\gamma y'),p})(\gamma'z,\gamma z')\Big)\gamma.
		\end{split}
	\end{equation}
Hence the trace in \eqref{9.22} can be written as
	\begin{equation}
		\frac{1}{\dim_{\mathbb{C}}H^{(0,0)}(N,L^p)}\int_N(\gamma')^{-1}T_{\overline{\mathscr{Q}}(y,\gamma'y'),p}T_{\mathscr{Q}(y,\gamma y'),p}(\gamma'z,\gamma z)\gamma dv_N(z).
	\end{equation}
	
	By \eqref{3.2} and \eqref{8.16},
	\begin{equation}
		\begin{split}
			\lv\frac{\mathrm{Tr}^{H^{(0,0)}(N,L^p)}}{\dim_{\mathbb{C}}H^{(0,0)}(N,L^p)}\Big[(\gamma')^{-1}T_{\overline{\mathscr{Q}},p}(y,\gamma'y')T_{\mathscr{Q},p}(y,\gamma y')\gamma\Big]\rv&\\
			\leqslant C\int_Ne^{-\sqrt{p}d_N(\gamma'z,\gamma z)}dv_N(z).&
		\end{split}
	\end{equation}
	Similar to the proof of Proposition \ref{p3.2n}, under \textbf{[FREE]} \eqref{4.8n}, the subset $\{z\mid \gamma z=\gamma'z\}\subset N$ has measure zero. This implies \eqref{9.22}.\qed
\end{pro}

\subsection{Spin-scale-frequency correspondence and pseudolocality}\label{s9.3}

We now explain why the mixed quantization $T_{\mathscr{Q},p}$ of $\mathscr{Q}\in C^\infty(K_{(k)}^N)$ provides a natural formalism for handling both large scale and large spin limits.

First, the discussion in \S\,\ref{s4.4n} may be summarized as the following broader \emph{spin-scale-frequency correspondence}, relating semiclassical analysis on a manifold $M$, harmonic analysis on a graph $X$, and geometric quantization on a Kähler manifold $N$,
\begin{equation}
	\begin{split}
		&\big(\text{quantization scheme}, \text{Hilbert space},\text{phase space}, \text{evolution}, \text{limit}\big)\\
		&=\begin{cases}
			\big(\text{Weyl},L^2(M),T^*M,\text{geodesic flow},\text{frequency }\lambda\to\infty\big),\\
			\big(\text{kernel operators}, L^2(X),\text{NB paths},\text{NB walk},\text{scale } m\to\infty\big),\\
			\big(\text{Berezin-Toeplitz}, H^{(0,0)}(N,L^p),N,\text{Hamiltonian flow}, \text{spin }p\to\infty\big).
		\end{cases}
	\end{split}
\end{equation}
Here, NB stands for nonbacktracking.

A key feature of Weyl quantization is \emph{pseudolocality}, and see \cite[\S\,8.4.5]{MR2952218} for details. This feature is also reflected in kernel operators, through Proposition \ref{p7.3}, and in Berezin-Toeplitz quantization, through Propositions \ref{8.3n} and \ref{8.4n}. In all settings, only contributions near the diagonal are significant.

On the one hand, if $k$ is fixed, then $Q\in C^\infty(K_{(k)},K^{\mathrm{End}(F)})$ only involves pairs $(y,y')$ whose distance is at most $k$. In the large scale limit, this becomes effectively near diagonal, and this is the key point underlying Proposition \ref{p7.3}. On the other hand, by Propositions \ref{8.3n} and \ref{8.4n}, the Berezin-Toeplitz quantization $T_{f,p}$ of $f\in C^\infty(N)$, as well as products of such operators, has kernels that are negligible away from the diagonal. This is the key point behind Proposition \ref{p9.7}. If we work with the full spaces $\big(C^\infty(K,K^{\mathrm{End}(F)}),\mathrm{End}(H^{(0,0)}(N,L^p))\big)$ instead of the pseudolocal pair $\big(C^\infty(K_{(k)},K^{\mathrm{End}(F)}),T_{f,p}\big)$, then this essential pseudolocality is lost, and we should not expect analogous results to hold.

\section{Quantum Ergodicity}\label{s10}

In this section, we prove several versions of QE. In \S\,\ref{s10.1}, we introduce the auxiliary quantity called quantum variance and establish some useful properties. In \S\,\ref{s10.2}, we establish uniform large scale QE. In \S\,\ref{s10.3}, we prove large spin QE. In \S\,\ref{s10.4}, we derive strong uniform QE. In \S\,\ref{s10.5}, we interpret these results in the Borel-Weil-Bott setting.

For large scale QE, we refer to Anantharaman-Le Masson \cite[Theorem 1]{MR3322309}. For large spin QE, we refer to Brooks-Le Masson-Lindenstrauss \cite[Theorem 1]{MR3567266}. For uniform QE in the manifold setting, we refer to Ma-Ma \cite[Theorem 1.1]{MR4808253}.

\subsection{Quantum variance}\label{s10.1}

Recall the setting of \eqref{2.19n}, \eqref{2.20n}, \eqref{2.21n}, \eqref{3.27}, \eqref{7.1n}, \eqref{7.2n}, and \eqref{7.5n}. For a kernel operator $Q\in C^\infty(K, K^{\mathrm{End}(F)})$, we define its quantum variance by
\begin{equation}
	\mathrm{Var}(Q)=\frac{1}{\dim_{\mathbb{C}}F}\sum_{i=1}^{\dim_{\mathbb{C}}F}\bv\big\langle Qu_{F,i},u_{F,i}\big\rangle_{L^2(X,F)}\bv^2.
\end{equation}

\begin{lemma}
Let $Q\in C^\infty(K, K^{\mathrm{End}(F)})$. Then
\begin{equation}\label{10.2}
	\begin{split}
		\mathrm{Var}(Q)\leqslant\lV Q\rV_{\mathrm{HS}(X,F)}^2.
	\end{split}
\end{equation}	
	Moreover, for any $j\in\mathbb{N}$ and $t_0>0$, we have
\begin{equation}\label{10.4}
	\begin{split}
		\mathrm{Var}(Q)&=\mathrm{Var}\bigg(\sum_{i\leqslant j}\frac{1}{{t_0}}\int_0^{t_0}\frac{(t\sqrt{-1}\mathrm{ad}_{\Delta^F})^i}{i!}Qdt\bigg)\\
		&=\mathrm{Var}\bigg(\frac{1}{{t_0}}\int_0^{t_0}e^{t\sqrt{-1}\mathrm{ad}_{\Delta^F}}Qdt\bigg).
	\end{split}
\end{equation}
\end{lemma}

\begin{pro}
Since $\Delta^F$ is selfadjoint and $(u_{F,i})$ are orthonormal eigensections,
\begin{equation}
	\begin{split}
		\mathrm{Var}(Q)\leqslant\frac{1}{\dim_{\mathbb{C}}F}\sum_{i=1}^{\dim_{\mathbb{C}}F}\bV&Qu_{F,i}\bV_{L^2(X,F)}^2=\lV Q\rV_{\mathrm{HS}(X,F)}^2,\\
\mathrm{Var}\big(\mathrm{ad}_{\Delta^F}Q\big)=\frac{1}{\dim_{\mathbb{C}}F}\sum_{i=1}^{\dim_{\mathbb{C}}F}\bv&\langle Qu_{F,i},\Delta^Fu_{F,i}\big\rangle_{L^2(X,F)}\\
&-\big\langle Q\Delta^Fu_{F,i},u_{F,i}\big\rangle_{L^2(X,F)}\bv^2=0,
	\end{split}
\end{equation}
from which \eqref{10.2} and \eqref{10.4} follow.\qed
\end{pro}

We shall be particularly interested in quantum variances of the form
\begin{equation}
\mathrm{Var}\big(Q_{(k)}-\langle Q_{(k)}\rangle_{K^{\mathrm{End}(F)}}\mathrm{Id}_{(k)}\big),
\end{equation}
where $\langle Q_{(k)}\rangle_{K^{\mathrm{End}(F)}}$ is given in \eqref{7.125}. By \eqref{7.6nn},
\begin{equation}
	\begin{split}
	&\mathrm{Var}\big(Q_{(k)}-\langle Q_{(k)}\rangle_{K^{\mathrm{End}(F)}}\mathrm{Id}_{(k)}\big)\\
	&=\frac{1}{\dim_{\mathbb{C}}F}\sum_{i=1}^{\dim_{\mathbb{C}}F}\bv\big\langle Q_{(k)}u_{F,i},u_{F,i}\big\rangle_{L^2(X,F)}-\langle Q_{(k)}\rangle h_k(\lambda_{F,i})\bv^2.
	\end{split}
\end{equation}
To emphasize its action on each eigensection, we shall also denote this quantum variance by
\begin{equation}
\mathrm{Var}\big(Q_{(k)}-\langle Q_{(k)}\rangle_{K^{\mathrm{End}(F)}} h_k(\lambda_{F,i})\big).
\end{equation}

\subsection{Uniform large scale QE}\label{s10.2}

We now state the uniform large scale QE.
\begin{theo}
		In the setting of \eqref{4.2n}, assume \textbf{[BST]} \eqref{4.3nn} and \textbf{[EXP]} \eqref{8.121}. Then, for any $k\in\mathbb{N}$ and $C>0$,
	\begin{equation}\label{10.5}
		\begin{split}
			\lim_{m\to\infty}\sup_{\substack{F_m\in \mathscr{E}(X_m,\varepsilon)\\ \lV Q_{m,(k)}\rV_{L^{\infty,2}(K_{m,(k)},K^{\mathrm{End}(F_m)})}\leqslant C}}&\\
			\mathrm{Var}\big(Q_{m,(k)}-&\langle Q_{m,(k)}\rangle_{K^{\mathrm{End}(F_m)}}h_k(\lambda_{F_m,i})\big)=0,
		\end{split}
	\end{equation}
	where, for each $m\in\mathbb{N}$, the supremum is taken first over all vector bundles $F_m$ with $\textbf{[EXP]}$, and then over all kernel operators $Q_{m}\in L^{\infty,2}(K_{m,(k)},K^{\mathrm{End}(F_m)})$ satisfying $\lV Q_{m}\rV_{L^{\infty,2}(K_{m,(k)},K^{\mathrm{End}(F_m)})}\leqslant C$.
\end{theo}

\begin{pro}
By \eqref{7.14n}, \eqref{10.2}, and \eqref{10.4},
\begin{equation}\label{10.6}
	\begin{split}
		\mathrm{Var}(Q_{(k)})&\leqslant\bbV\sum_{i\leqslant j}\frac{1}{{t_0}}\int_0^{t_0}\frac{(t\sqrt{-1}\mathrm{ad}_{\Delta^F})^i}{i!}Q_{(k)}dt\bbV_{\mathrm{HS}(X,F)}^2\\
		&\leqslant\bbV\sum_{i\leqslant j}\frac{1}{{t_0}}\int_0^{t_0}\frac{(t\sqrt{-1}\mathrm{ad}_{\Delta^F})^i}{i!}Q_{(k)}dt\bbV_{L^2(K_{},K^{\mathrm{End}(F)})}^2\\
		&\quad+\frac{\lv\{x\in X\mid \mathrm{inj}_x\leqslant k+j\}\rv\lv B_{ \mathbb{T}_d}(o,k+j)\rv^2}{\lv X\rv}\\
	&\qquad\cdot\bbV\sum_{i\leqslant j}\frac{1}{{t_0}}\int_0^{t_0}\frac{(t\sqrt{-1}\mathrm{ad}_{\Delta^F})^i}{i!}Q_{(k)}dt\bbV_{L^\infty(K,K^{\mathrm{End}(F)})}^2.
	\end{split}
\end{equation}
From \eqref{7.38}, we have
\begin{equation}\label{10.7}
	\begin{split}
		\bbV\sum_{i\leqslant j}\frac{1}{{t_0}}\int_0^{t_0}\frac{(t\sqrt{-1}\mathrm{ad}_{\Delta^F})^i}{i!}Q_{(k)}dt&\bbV_{L^\infty(K,K^{\mathrm{End}(F)})}\\
		\leqslant\Big(\sum_{i\leqslant j}\frac{1}{{t_0}}\int_0^{t_0}\frac{(2dt)^i}{i!}dt\Big)\bV Q_{(k)}&\bV_{L^\infty(K,K^{\mathrm{End}(F)})}.
	\end{split}
\end{equation}

Combining \eqref{7.126}, \eqref{10.6}, and \eqref{10.7}, we see that, under \textbf{[BST]} \eqref{4.3nn} and \textbf{[EXP]} \eqref{7.121}, the limit \eqref{10.5} is bounded by
\begin{equation}\label{10.8}
	\begin{split}
\lim_{m\to\infty}\sup_{\substack{\mathrm{End}(F_m)\in \mathscr{F}(X_m,\varepsilon), Q_{m,(k)}\perp\mathrm{Id}_{(k)},\\ \lV Q_{m,(k)}\rV_{L^{\infty,2}(K_{m,(k)},K^{\mathrm{End}(F_m)})}\leqslant 1}}&\\
\bbV\sum_{i\leqslant j}\frac{1}{{t_0}}\int_0^{t_0}\frac{(t\sqrt{-1}\mathrm{ad}_{\Delta^{F_m}})^i}{i!}&Q_{m,(k)}dt\bbV_{L^2(K_{m},K^{\mathrm{End}(F_m)})}^2.
	\end{split}
\end{equation}

By \eqref{7.38}, we have
\begin{equation}\label{10.9}
	\begin{split}
		&\bbV\sum_{i\leqslant j}\frac{1}{{t_0}}\int_0^{t_0}\frac{(t\sqrt{-1}\mathrm{ad}_{\Delta^F})^i}{i!}Q_{(k)}dt\bbV_{L^2(K,K^{\mathrm{End}(F)})}\\
		&\leqslant \bbV \frac{1}{{t_0}}\int_0^{t_0}e^{t\sqrt{-1}\mathrm{ad}_{\Delta^F}}Q_{(k)}dt\bbV_{L^2(K,K^{\mathrm{End}(F)})}\\
		&\quad+\Big(\sum_{i>j}\frac{1}{{t_0}}\int_0^{t_0}\frac{(4td^{1/2})^i}{i!}dt\Big)\bV Q_{(k)}\bV_{L^2(K_{(k)},K^{\mathrm{End}(F)})}.
	\end{split}
\end{equation}
Taking together \eqref{7.17}, \eqref{7.112}, and \eqref{10.9}, the limit \eqref{10.8} is bounded by
\begin{equation}\label{10.10}
	\bigg(\frac{C_{k,\varepsilon}}{{t_0}^{1/7}}+\sum_{i>j}\frac{1}{{t_0}}\int_0^{t_0}\frac{(4td^{1/2})^i}{i!}dt\bigg)^2\lv B_{\mathbb{T}_d}(o,k)\rv.
\end{equation}
Letting $j\to\infty$ and then ${t_0}\to\infty$ in \eqref{10.10}, we obtain \eqref{10.5}.\qed
\end{pro}

\subsection{Large spin QE}\label{s10.3}

Recall \eqref{9.12} and \eqref{9.17}. We now state the large spin QE.
\begin{theo}
In the setting of \eqref{4.5n}, assume \textbf{[FREE]} \eqref{4.8n} and \textbf{[DEN]} \eqref{4.9n}. Then for any $\mathscr{Q}_{(k)}\in C^\infty(K^N_{(k)})$,
	\begin{equation}\label{10.11}
			\lim_{p\to\infty}\mathrm{Var}\big(T_{\mathscr{Q}_{(k)},p}-\langle \mathscr{Q}_{(k)}\rangle_{K^N} h_k(\lambda_{F_p,i})\big)=0.
	\end{equation}
\end{theo}

\begin{pro}
By \eqref{9.13}, under \textbf{[DEN]} \eqref{4.9n}, we may suppose that $\mathscr{Q}_{(k)}\perp \mathbbm{1}_{(k)}$. From \eqref{9.18}, \eqref{10.2}, and \eqref{10.4},
\begin{equation}\label{10.12}
	\begin{split}
	&\lim_{p\to\infty}	\mathrm{Var}(T_{\mathscr{Q}_{(k)},p})\\
	&\leqslant\lim_{p\to\infty} \bbV\sum_{i\leqslant j}\frac{1}{{t_0}}\int_0^{t_0}\frac{(t\sqrt{-1}\mathrm{ad}_{\Delta^{F_p}})^i}{i!}T_{\mathscr{Q}_{(k)},p}dt\bbV_{\mathrm{HS}(X,F_p)}^2\\
		&=\lim_{p\to\infty}\BV T_{\sum_{i\leqslant j}\frac{1}{{t_0}}\int_0^{t_0}\frac{(t\sqrt{-1}\mathrm{ad}_{\Delta^{X^N}})^i}{i!}\mathscr{Q}_{(k)}dt,p}\BV_{\mathrm{HS}(X,F_p)}^2.
	\end{split}
\end{equation}
By \eqref{9.19} and \eqref{9.20}, under \textbf{[FREE]} \eqref{4.8n}, the limit \eqref{10.12} is bounded by \begin{equation}\label{10.12n}
	\begin{split}
&\lim_{p\to\infty}\BV T_{\sum_{i\leqslant j}\frac{1}{{t_0}}\int_0^{t_0}\frac{(t\sqrt{-1}\mathrm{ad}_{\Delta^{X^N}})^i}{i!}\mathscr{Q}_{(k)}dt,p}\BV_{L^2(K,K^{\mathrm{End}(F_p)})}^2\\
		&\leqslant C\bbV\sum_{i\leqslant j}\frac{1}{{t_0}}\int_0^{t_0}\frac{(t\sqrt{-1}\mathrm{ad}_{\Delta^{X^N}})^i}{i!}\mathscr{Q}_{(k)}dt\bbV_{L^2(K^N)}^2.
	\end{split}
\end{equation}
By \eqref{9.14}, we take $j\to\infty$ and then $t_0\to \infty$ in \eqref{10.12n} to get \eqref{10.11}.\qed
\end{pro}

Note that we cannot simply replace the finite sum in \eqref{10.12} with $e^{t\sqrt{-1}\mathrm{ad}_{\Delta^{F_p}}}T_{\mathscr{Q},p}$, because this would produce a kernel function $e^{t\sqrt{-1}\mathrm{ad}_{\Delta^{X^N}}}\mathscr{Q}$ with infinite propagation speed, and then \eqref{9.20} would fail.

\subsection{Strong uniform QE}\label{s10.4}

We now state the strong uniform QE.
\begin{theo}\label{t10.3}
In the setting of \eqref{4.11n}, assume \textbf{[BST]} \eqref{4.3nn}, and assume \textbf{[FREE]} \eqref{4.8n} and \textbf{[EXP]} \eqref{9.15} for every $\Gamma_m$. Then, for any series $(\mathscr{Q}_{m,(k)}\in C^\infty(K^N_{m,(k)}))_{m\in\mathbb{N}}$ satisfying $\lV\mathscr{Q}_{m,(k)}\rV_{L^{\infty,2}(K^N_{m,(k)})}\leqslant C$,	
\begin{equation}\label{10.13}
		\lim_{\max(m,p)\to\infty}\mathrm{Var}\big(T_{\mathscr{Q}_{m,(k)},p}-\langle \mathscr{Q}_{m,(k)}\rangle_{K^N_m} h_k(\lambda_{F_{m,p},i})\big)=0.
\end{equation}
\end{theo}

\begin{pro}
First, we show that
\begin{equation}\label{10.18}
	\langle T_{\mathscr{Q}_{m,(k)},p}\rangle_{K^{\mathrm{End}(F_{m,p})}}=\langle \mathscr{Q}_{m,(k)}\rangle_{K_m^N}.
\end{equation}
By \eqref{7.125} and \eqref{9.12}, it suffices to prove the following exact version of \eqref{8.8}, for any $p\in\mathbb{N}$ and $f\in C^\infty(N)$,
\begin{equation}\label{10.17}
\frac{1}{\dim_{\mathbb{C}} H^{(0,0)}(N,L^p)}\tro^{H^{(0,0)}(N,L^p)}[T_{f,p}]=\frac{1}{\mathrm{Vol}(N)}\int_{N}fdv_{N}.
\end{equation}
We claim that this holds under \textbf{[DEN]} \eqref{4.9n}, weaker than \textbf{[EXP]} \eqref{9.15}. Indeed, the Bergman kernel given in \eqref{3.5n} is equivariant, that is,
\begin{equation}
P_p(\gamma z,\gamma z')=\gamma P_p(z,z')\gamma^{-1}.
\end{equation}
Applying this to the diagonal $(z,z)$, we see that the function $z\in N\mapsto P_p(z,z)\in\mathbb{C}$ is constant. Since
\begin{equation}
\int_NP_p(z,z)dv_N(z)=\dim_\mathbb{C}H^{(0,0)}(N,L^p),
\end{equation}
we get the following identity, analogous to \eqref{3.33nn},
\begin{equation}\label{10.20}
P_p(z,z)=\frac{\dim_\mathbb{C}H^{(0,0)}(N,L^p)}{\mathrm{Vol}(N)}.
\end{equation}
By \eqref{8.6} and \eqref{10.20}, we get \eqref{10.17}.

By \eqref{9.8}, \eqref{9.16}, \eqref{9.18}, \eqref{9.19}, and \eqref{10.18}, we repeat the proof of \eqref{10.5} to obtain
\begin{equation}\label{10.14}
	\lim_{m\to\infty}\sup_{p\in\mathbb{N}}\mathrm{Var}\big(T_{\mathscr{Q}_{m,(k)},p}-\langle \mathscr{Q}_{m,(k)}\rangle_{K^N_m}\mathrm{Id}_{(k)}\big)=0.
\end{equation}

Replacing $(F_p,T_{\mathscr{Q}_{(k)},p})_{p\in\mathbb{N}}$ with $(F_{m,p},T_{\mathscr{Q}_{m,(k)},p})_{p\in\mathbb{N}}$ in \eqref{10.11}, we get
\begin{equation}\label{10.15}
	\sup_m\lim_{p\to\infty}\mathrm{Var}(T_{\big(\mathscr{Q}_{m,(k)}-\langle \mathscr{Q}_{m,(k)}\rangle_{K^N_m}\mathbbm{1}_{(k)}\big),p})=0.
\end{equation}

Combining \eqref{10.14} and \eqref{10.15}, we get \eqref{10.13} from \eqref{4.18n}.\qed
\end{pro}

As noted after Lemma  \ref{l9.6n}, a key point in the proof of \eqref{10.14} is that, in \eqref{9.19}, the right hand sides of the inequalities are uniform in $p\in\mathbb{N}$, with no $O(p^{-1})$ remainder term.

Moreover, if we replace \textbf{[EXP]} \eqref{9.15} with the stronger assumption that \eqref{7.121} holds for all $(F_{m,p})_{m,p\in\mathbb{N}}$, then \eqref{10.14} follows directly from \eqref{10.5} by substituting $(F_m,Q_{m,(k)})$ with $(F_{m,p},T_{\mathscr{Q}_{m,(k)},p})$. This avoids the infinite dimensional analysis in \S\,\ref{s9.1}. However, because of the pseudolocality discussed in \S\,\ref{s9.3}, even under this stronger assumption we cannot extend the large spin QE from $T_{\mathscr{Q}_{m,(k)},p}$ to the space $C^\infty(K_{m,(k)},K^{\mathrm{End}(F_{m,p})})$. For this reason, we formulate the main result in its present form under \eqref{9.15}, rather than imposing the stronger assumption that \eqref{7.121} holds for all $(F_{m,p})_{m,p\in\mathbb{N}}$.

\subsection{QE in the Borel-Weil-Bott case}\label{s10.5}

We now state the main results of this section in the Borel-Weil-Bott setting.

\begin{theo}
	In the setting of \eqref{4.20nn}, assume \textbf{[BST]} \eqref{4.3nn} and \textbf{[EXP]} \eqref{9.15n}, and suppose that $G$ is compact connected. For any $k\in\mathbb{N}$ and $C>0$,
		\begin{equation}
			\begin{split}
	\lim_{m\to\infty}&\sup_{\substack{\alpha\in \Lambda^*\cap\overline{C^+},\\ \Vert Q_{m,\alpha,(k)}\Vert_{L^{\infty,2}(K_{m,(k)},K^{\mathrm{End}(F_{m,\alpha})})}\leqslant C}}\\
	&\mathrm{Var}\big(Q_{m,\alpha,(k)}-\langle Q_{m,\alpha,(k)}\rangle_{K^{\mathrm{End}(F_{m,\alpha})}}h_k(\lambda_{F_{m,\alpha},i})\big)=0.
			\end{split}
	\end{equation}
\end{theo}

\begin{theo}
	In the setting of \eqref{4.11nn}, assume \textbf{[FREE]} and \textbf{[DEN]} \eqref{4.10n}. Suppose, in addition, that either $G$ is compact connected simple and $\alpha$ is an arbitrary highest weight, or that $G$ is compact connected and $\alpha$ is regular. Then for any $\mathscr{Q}_{(k)}\in C^\infty(K^{\mathcal{O}_\alpha}_{(k)})$,
	\begin{equation}
		\lim_{p\to\infty}\mathrm{Var}\big(T_{\mathscr{Q}_{(k)},p}-\langle \mathscr{Q}_{(k)}\rangle_{K^{\mathcal{O}_\alpha}} h_k(\lambda_{F_{p\alpha},i})\big)=0.
	\end{equation}
\end{theo}

\begin{theo}\label{t10.6}
	In the setting of \eqref{4.20nn}, assume \textbf{[BST]} \eqref{4.3nn}, and assume \textbf{[FREE]} \eqref{4.10n} and \textbf{[EXP]} \eqref{9.15n} for every $\Gamma_m$. Suppose, in addition, that either $G$ is compact connected simple and $\alpha$ is an arbitrary highest weight, or that $G$ is compact connected and $\alpha$ is regular. Then for any series $(\mathscr{Q}_{m,(k)}\in C^\infty(K^{\mathcal{O}_\alpha}_{m,(k)}))_{m\in\mathbb{N}}$ satisfying $\lV\mathscr{Q}_{m,(k)}\rV_{L^{\infty,2}(K^{\mathcal{O}_\alpha}_{m,(k)})}\leqslant C$,	
	\begin{equation}
		\lim_{\max(m,p)\to\infty}\mathrm{Var}\big(T_{\mathscr{Q}_{m,(k)},p}-\langle \mathscr{Q}_{m,(k)}\rangle_{K^{\mathcal{O}_\alpha}_m} h_k(\lambda_{F_{m,p\alpha},i})\big)=0.
	\end{equation}
\end{theo}

\section{Zero Divisor Equidistribution}\label{s11}

In this section, we prove the zero divisor equidistribution. In \S\,\ref{s11.1}, we recall a result of Nonnenmacher-Voros \cite[Theorem 1]{MR1649013} and Shiffman-Zelditch \cite[Lemma 1.4]{MR1675133}. In \S\,\ref{s11.2}, we establish the zero divisor equidistribution. In \S\,\ref{s11.3}, we interpret the result in the Borel-Weil-Bott setting.

Note that zero divisor equidistribution only involves the large spin limit, at present, no large scale analogue is known.

\subsection{QE and zero divisor equidistribution}\label{s11.1}

Recall the setting of \S\,\ref{s3.1}. In particular, let $N$ be a complex manifold, $(L,h^L)$ a positive Hermitian holomorphic line bundle over $N$, and $c_1(L,h^L)$ the first Chern form of $(L,h^L)$. For any holomorphic section $u_p\in H^{(0,0)}(N,L^p)$, its zero divisor $\mathrm{Div}(u_p)$ is a formal sum of analytic subvarieties of $N$ of complex codimension one.

We now state an important result on the zero divisors of large spin holomorphic sections, established in \cite[Theorem 1]{MR1649013} for the theta bundle over an elliptic curve and in \cite[Lemma 1.4]{MR1675133} in the general case, both by means of potential theory.
\begin{prop}
	Let $(u_p\in H^{(0,0)}(N,L^p))_{p\in\mathbb{N}}$ be
	a series of sections such that $\lim_{p\to\infty}\lV u_p(z)\rV_{L^p}^2=1$
	in the weak star sense, namely, for any $f\in C(N)$,
	\begin{equation}\label{11.1}
		\lim_{p\to\infty}\int_N\lV u_p(z)\rV^2_{L^p}f(z)dv_N(z)=\int_Nf(z)dv_N(z).
	\end{equation} 
	Then $\lim_{p\to\infty}[\mathrm{Div}(u_{p})]/p=c_1(L,h^L)$ weakly in the sense of currents, that is, for any continuous $(n-1,n-1)$-form $\phi^{(n-1,n-1)}$ on $N$,
	\begin{equation}\label{11.2}
		\lim_{p\to\infty}\frac{1}{p}\int_{\mathrm{Div}(u_{p})}\phi^{(n-1,n-1)}=\int_{N}\phi^{(n-1,n-1)}\wedge c_1(L,h^L).
	\end{equation}
\end{prop}

\subsection{Zero divisor equidistribution}\label{s11.2}

For any eigensection $u_{p,i}\in C^\infty(X,F_p)$, its divisor $\mathrm{Div}(u_{p,i})$ can be viewed as a current on $X^N$, defined in \eqref{9.1}.

We now state the large spin zero divisor equidistribution.
\begin{theo}\label{t11.2}
	In the setting of \eqref{4.5n}, assume \textbf{[FREE]} \eqref{4.8n} and \textbf{[DEN]} \eqref{4.9n}. Then there exists a series of subsets $\big(I_p\subseteq \{1,\cdots, \dim_{\mathbb{C}}F_p\}\big)_{p\in\mathbb{N}}$ of asymptotic density one, in the sense that
	\begin{equation}\label{11.3}
		\lim_{p\to\infty}\lv I_p\rv\big/\dim_{\mathbb{C}}F_p=1,
	\end{equation}
	such that for any series $(i_p\in I_p)_{p\in\mathbb{N}}$, and for any continuous $(n-1,n-1)$-form $\mathscr{Q}^{(n-1,n-1)}_{(0)}$ on $X^N$, we have
	\begin{equation}\label{11.4}
		\begin{split}
			\lim_{p\to\infty}\frac{1}{p}&\sum_{y\in D}\int_{\mathrm{Div}(u_{p,i_p}(y,\cdot))}\mathscr{Q}^{(n-1,n-1)}_{(0)}(y,z)\\
			=&\sum_{y\in D}\int_{N}\mathscr{Q}^{(n-1,n-1)}_{(0)}(y,z)\wedge c_1(L,h^L).
		\end{split}
	\end{equation}
\end{theo}

\begin{pro}
By the classical formulation of QE on manifolds, there exists a series of subsets $(I_p\subseteq \{1,\cdots, \dim_{\mathbb{C}}F_p\})_{p\in\mathbb{N}}$ satisfying \eqref{11.3}, such that for any series $(i_p)_{p\in\mathbb{N}}$ satisfying $i_p\in I_p$ for any $p\in\mathbb{N}$, and for any $\mathscr{Q}_{(0)}\in C^\infty(X^N)$,
\begin{equation}
	\begin{split}
\lim_{p\to\infty}\sum_{y\in D}\int_N\mathscr{Q}_{(0)}(y,z)\Vert u_{p,i_p}(y,z)\Vert_{L^p}^2dv_N(z)&\\
=\frac{1}{\lv X\rv\mathrm{Vol}(N)}\sum_{y\in D}\int_N\mathscr{Q}_{(0)}(y,z)dv_N(z)&.
	\end{split}
\end{equation}
Then \eqref{11.4} follows from \eqref{11.1} and \eqref{11.2}.
\end{pro}

\subsection{Zero divisor equidistribution in the Borel-Weil-Bott case}\label{s11.3}

We now state the main result of this section in the Borel-Weil-Bott setting.

\begin{theo}\label{t11.3}
	In the setting of \eqref{4.11nn}, assume \textbf{[FREE]} and \textbf{[DEN]} \eqref{4.10n}. Suppose, in addition, that either $G$ is compact connected simple and $\alpha$ is an arbitrary highest weight, or that $G$ is compact connected and $\alpha$ is regular. Then there is a series of subsets $I_{p\alpha}\subseteq \{1,\cdots, \dim_{\mathbb{C}}F_{p\alpha}\}$ of asymptotic density one, in the sense that
	\begin{equation}
		\lim_{p\to\infty}\lv I_{p\alpha}\rv\big/\dim_{\mathbb{C}}F_{p\alpha}=1,
	\end{equation}
	such that for any series $(i_p\in I_{p\alpha})_{p\in\mathbb{N}}$, and for any continuous $(n-1,n-1)$-form $\mathscr{Q}^{(n-1,n-1)}_{(0)}(y,z)$ on $X^{\mathcal{O}_\alpha}$,
	\begin{equation}
		\begin{split}
			\lim_{p\to\infty}\frac{1}{p}&\sum_{y\in D}\int_{\mathrm{Div}(u_{p,i_p}(y,\cdot))}\mathscr{Q}^{(n-1,n-1)}_{(0)}(y,z)\\
			=&\sum_{y\in D}\int_{\mathcal{O}_\alpha}\mathscr{Q}^{(n-1,n-1)}_{(0)}(y,z)\wedge c_1(L_\alpha,h^{L_\alpha}).
		\end{split}
	\end{equation}
\end{theo}

\section{Ramanujan Vector Bundles and Questions}\label{s12}

In \S\,\ref{s12.1}, we introduce general Ramanujan vector bundles. In \S\,\ref{s12.2}, we describe various questions and conjectures.

\subsection{Ramanujan vector bundles}\label{s12.1}

Let $a,b\in \mathbb{Q}$ with $a,b<0$, and let $\mathbb{D}_{a,b}$ be the associated quaternion algebra over $\mathbb{Q}$, defined by
\begin{equation}
\mathbb{D}_{a,b}=\{x_0+x_1\mathrm{i}+x_2\mathrm{j}+x_3\mathrm{ij}\mid x_0,\cdots,x_3\in \mathbb{Q},\ \mathrm{i}^2=a, \mathrm{j}^2=b, \mathrm{ij}=-\mathrm{ji}\},
\end{equation}
which ramifies at $\infty$. The reduced norm $\mathrm{Nrd}$ on $\mathbb{D}_{a,b}$ is defined by
\begin{equation}
	\mathrm{Nrd}(x_0+x_1\mathrm{i}+x_2\mathrm{j}+x_3\mathrm{ij})=x_0^2-ax_1^2-bx_2^2+abx_3^2.
\end{equation}

Let $p_0$ be a prime at which $\mathbb{D}_{a,b}$ splits, that is,
\begin{equation}\label{12.3}
\mathbb{D}_{a,b}(\mathbb{Q}_{p_0})\cong \mathrm{M}_2(\mathbb{Q}_{p_0}).
\end{equation}
This is equivalent to the condition that the equation
\begin{equation}
ax_0^2+bx_1^2=x_2^2
\end{equation}
has a nontrivial solution $(x_0,x_1,x_2)\in\mathbb{Q}_{p_0}^3$, see Serre \cite[\S\S\,3, 4]{MR344216}. Then the isomorphism \eqref{12.3} can be written explicitly by
\begin{equation}
\begin{cases}
\mathrm{i}\to \big(\begin{smallmatrix}
\sqrt{a}&0\\0&-\sqrt{a}
\end{smallmatrix}\big),\quad\mathrm{j}\to \big(\begin{smallmatrix}
0&b\\1&0
\end{smallmatrix}\big),&\text{if } a \text{ is a square in }\mathbb{Q}_{p_0},\\
\mathrm{i}\to \big(\begin{smallmatrix}
0&a\\1&0
\end{smallmatrix}\big),\quad\mathrm{j}\to \big(\begin{smallmatrix}
	x_2&ax_0\\-x_0&-x_2
\end{smallmatrix}\big),&\text{if }  ax_0^2+b=x_2^2 \text{ for } (x_0,x_2)\in\mathbb{Q}_{p_0}^2.
\end{cases}
\end{equation}

Let $\mathbb{G}_{a,b}$ be the algebraic group of the invertible elements of $\mathbb{D}_{a,b}$, in other words, $\mathrm{Nrd}\neq0$. Then $\Gamma=\mathbb{G}_{a,b}\big(\mathbb{Z}\big[\tfrac{1}{p_0}\big]\big)$ can be embedded to a discrete cocompact lattice of $\mathrm{PGL}_2(\mathbb{Q}_{p_0})$, and we have a $\mathrm{SU}_2(\mathbb{C})$-representation of $\Gamma$ by composing
\begin{equation}\label{12.6}
	\rho\colon\Gamma=\mathbb{G}_{a,b}\big(\mathbb{Z}\big[\tfrac{1}{p_0}\big]\big)\to \mathbb{G}_{a,b}(\mathbb{R})\cong\mathrm{SU}_2(\mathbb{C})\times \mathbb{R}^*\to \mathrm{SU}_2(\mathbb{C}).
\end{equation}
For another prime $p_1\neq p_0$, let $\Gamma_{(p_0,p_1)}\subset\Gamma$ be the congruence subgroup
\begin{equation}\label{12.7}
	\begin{split}
\Gamma_{(p_0,p_1)}=\ker\Big(&\mathbb{G}_{a,b}\big(\mathbb{Z}\big[\tfrac{1}{p_0}\big]\big)\\
&\to\mathbb{G}_{a,b}\big(\mathbb{Z}\big[\tfrac{1}{p_0}\big]\big/p_1\mathbb{Z}\big[\tfrac{1}{p_0}\big]\big)\cong \mathbb{G}_{a,b}\big(\mathbb{Z}/p_1\mathbb{Z}\big)\big).
	\end{split}
\end{equation}

We now define the \emph{Ramanujan vector graph} and \emph{Ramanujan vector bundles} for $p\in\mathbb{N}$ by
\begin{equation}\label{12.8}
	\begin{split}
		X_{(p_0,p_1)}&=\Gamma_{(p_0,p_1)}\big\backslash\big(\mathrm{PGL}_2(\mathbb{Q}_{p_0})\big/\mathrm{PGL}(\mathbb{Z}_{p_0})\big),\\
			F_{(p_0,p_1),p}&=\Gamma_{(p_0,p_1)}\big\backslash\Big(\big(\mathrm{PGL}_2(\mathbb{Q}_{p_0})/\mathrm{PGL}(\mathbb{Z}_{p_0})\big)\times \mathrm{Sym}^p(\mathbb{C}^2)\Big)
	\end{split}
\end{equation}
where $\mathrm{PGL}_2(\mathbb{Q}_{p_0})/\mathrm{PGL}(\mathbb{Z}_{p_0})$ is the \emph{Bruhat-Tits tree},
\begin{equation}
	\mathrm{PGL}_2(\mathbb{Q}_{p_0})/\mathrm{PGL}(\mathbb{Z}_{p_0})\cong \mathbb{T}_{p_0+1}.
\end{equation}
Note that \eqref{1.30n} is a special case of \eqref{12.8} when the class number of $\mathbb{D}_{a,b}$ is one.

The following result is by Lubotzky-Phillips-Sarnak \cite[Theorem 4.1]{MR890171}.
\begin{theo}
	The vector bundles $F_{(p_0,p_1),p}$ are Ramanujan vector bundles, that is, if $p\neq 0$ or if $p=0$ and $\lambda_{(p_0,p_1),0,i}\neq\pm(p_0+1)$, then
	\begin{equation}\label{12.10}
		\lv\lambda_{(p_0,p_1),p,i}\rv\leqslant 2\sqrt{p_0}.
	\end{equation}
\end{theo}

We now show that our main results apply to Ramanujan vector bundles.
\begin{theo}
Theorems \ref{t5.6}, \ref{t6.6}, \ref{nt6.6}, \ref{t10.6}, and \ref{t11.2} all apply to the Ramanujan vector bundles $(F_{(p_0,p_1),p})_{p_1 \text{prime}, p\in\mathbb{N}}$ defined in \eqref{12.8}. Here, the prime $p_1$ replaces the large scale parameter $m$, and the complex geometric data are as in \eqref{1.26n}.
\end{theo}

\begin{pro}
It is enough to check \textbf{[BST]} \eqref{4.3nn} for $(X_{(p_0,p_1)})_{p_1\text{ prime}}$, [\textbf{GAP}] \eqref{n7.3} for each $F_{(p_0,p_1),p}$, and \textbf{[FREE]} \eqref{4.10n} and \textbf{[EXP]} \eqref{9.15n} for each $(\Gamma_{(p_0,p_1)})_{p_1\text{ prime}}$, since \textbf{[FREE]} \eqref{4.10n} is stronger than \textbf{[DEN]} \eqref{4.10n}.

\textbf{[BST]} \eqref{4.3nn} for $(X_{(p_0,p_1)})_{p_1\text{ prime}}$ follows from Lubotzky-Phillips-Sarnak \cite[Theorem 3.4]{MR963118}, where they proved a stronger result that the girth goes to infinity as $p_1\to\infty$.

\textbf{[FREE]} \eqref{4.10n} is a direct consequence of the fact that $\Gamma_{(p_0,p_1)}$, being the fundamental group of a graph, is free, together with \eqref{12.6}. It suffices to show that the $\Gamma_{(p_0,p_1)}$-action on $\big(\mathrm{PGL}_2(\mathbb{Q}_{p_0})\big/\mathrm{PGL}(\mathbb{Z}_{p_0})\big)$ has no vertex stabilizer. If not, then $\gamma g_{p_0}\mathrm{PGL}(\mathbb{Z}_{p_0})=g_{p_0}\mathrm{PGL}(\mathbb{Z}_{p_0})$ for some $g_{p_0}\in\mathrm{PGL}(\mathbb{Q}_{p_0})$. Hence $\gamma\in \Gamma\cap g_{p_0}\mathrm{PGL}(\mathbb{Z}_{p_0})g_{p_0}^{-1}$, the intersection of a discrete and a compact group, must be a finite group. Therefore, $\gamma$ is of finite order. On the other hand, since $\Gamma_{(p_0,p_1)}$ is a congruence subgroup, Selberg's lemma \cite[Theorem 4.8.2]{MR3307755} implies that it contains no nontrivial finite order element.

[\textbf{GAP}] \eqref{n7.3} and \textbf{[EXP]} \eqref{9.15n} are more difficult. Fortunately, they follow from \eqref{12.10}, the Peter-Weyl theorem
\begin{equation}
L^2(\mathrm{SU}_2(\mathbb{C}))=\mathsmaller{\bigoplus}_{p=0}^\infty (p+1)\mathrm{Sym}^p(\mathbb{C}^2),
\end{equation}
and the fact that multiplicities of nontrivial representations do not affect the expander condition.\qed
\end{pro}

\subsection{Questions and conjectures}\label{s12.2}

We now discuss various questions and conjectures.

\subsubsection{Cayley vector bundles}

The construction in \eqref{1.30} is based on Cayley graphs and Cayley vector bundles, which we now define.

Let $H$ be a finite group and $S\subset H$ a symmetric subset, meaning that $h\in S$ if and only if $h^{-1}\in S$. The Cayley graph of $H$ with respect to $S$ is the graph whose vertices are the elements of $H$, and in which two vertices $(h,h')$ are connected by an edge if $h=h''h'$ for some $h''\in S$. Let $G$ be a compact connected Lie group. For any $h\in S$, assign an element $\phi(h)\in G$ such that $\phi(h^{-1})=(\phi(h))^{-1}$. Then we define an adjacency operator $\Delta^{H\times G}_{\phi}$ by
\begin{equation}
	\begin{split}
&\Delta^{H\times G}_{\phi}\colon L^2(H\times G)\to L^2(H\times G),\\
&(\Delta^{H\times G}_{\phi}f)(h,g)=\sum_{h'\in S}f(h'h,\phi(h')g).
	\end{split}
\end{equation}
Formally, this gives an infinite graph with vertex set $H\times G$, in which two vertices $(h,g)$ and $(h',g')$ are adjacent if there exists $h''\in S$ such that $(h,g)=(h''h',\phi(h'')g')$.

Naor-Sah-Sawhney-Zhao \cite{MR4651019} and Magee-Thomas-Zhao \cite{MR4630491} studied QE for Cayley graphs. It would be interesting to study QE for Cayley vector bundles.

\subsubsection{Spectral gap and Ramanujan vector bundles}

Friedman \cite{MR1137767,MR2437174} showed that many families of random regular graphs are weakly Ramanujan. Huang-McKenzie-Yau \cite{huang2025ramanujanpropertyedgeuniversality} proved edge universality and consequently that approximately 69\% of regular graphs are Ramanujan. It would be interesting to formulate a suitable probabilistic model, for instance by independently sampling a gauge element on each edge according to the Haar measure on a compact connected Lie group, and then to study the proportion of Ramanujan vector bundles.

Gamburd-Jakobson-Sarnak \cite{MR1677685} and Bourgain-Gamburd \cite{MR2358056,MR2966656} obtained spectral gap results for averaging operators on $\mathrm{SU}_k(\mathbb{C})$ under certain arithmetic Diophantine assumptions. In view of \eqref{1.41} and \eqref{1.41n}, it would be interesting to generalize these results to vector bundles over more general graphs.

Marcus-Spielman-Srivastava \cite{MR3374962} and Hall-Puder-Sawin \cite{MR3725881} constructed Ramanujan coverings of a Ramanujan graph. This raises the question of whether every Ramanujan graph admits a Ramanujan vector bundle.

\subsubsection{(A)UQUE}

Rudnick-Sarnak \cite{MR1266075} formulated the quantum unique ergodicity (QUE) conjecture. Lindenstrauss \cite{MR2195133} proved arithmetic QUE (AQUE), and Soundararajan \cite{MR2680500} ruled out escape of mass in the noncompact case.

In the discrete setting, we may formulate uniform QUE (UQUE) as follows, for any series $(\mathscr{Q}_{m,(k)}\in C^\infty(K^N_{m,(k)}))_{m\in\mathbb{N}}$ satisfying $\lV\mathscr{Q}_{m,(k)}\rV_{L^{\infty,2}(K^N_{m,(k)})}\leqslant C$, do we have
\begin{equation}
	\begin{split}
\lim_{\max(m,p)\to\infty}\sup_{1\leqslant i\leqslant\dim_{\mathbb{C}}F_{m,p}}\Bv\big\langle T_{\mathscr{Q}_{m,(k)},p}u_{F_{m,p},i},u_{F_{m,p},i}\big\rangle_{L^2(X_m,F_{m,p})}&\\
-\langle \mathscr{Q}_{m,(k)}\rangle_{K_m^N} h_k(\lambda_{F_{m,p},i})\Bv=&0.
	\end{split}
\end{equation}
In particular, for the arithmetic Ramanujan vector bundles in \eqref{12.8}, we may call this arithmetic UQUE (AUQUE).

\subsubsection{Nonregular graphs}

Anantharaman-Sabri \cite{MR3961083} extended the work of Anantharaman-Le Masson \cite{MR3322309} to disordered systems, that is, to nonregular graphs, possibly with a potential on the vertices or weights on the edges. This applies in particular to graphs converging to the Anderson model.

This raises the question of extending the main results of this paper to the nonregular setting.

\def\cprime{$'$} \def\cprime{$'$}

\begin{thebibliography}{10}
	
	\bibitem{MR875835}
	N.~Alon.
	\newblock Eigenvalues and expanders.
	\newblock volume~6, pages 83--96. 1986.
	\newblock Theory of computing (Singer Island, Fla., 1984).
	
	\bibitem{MR3649482}
	N.~Anantharaman.
	\newblock Quantum ergodicity on regular graphs.
	\newblock {\em Comm. Math. Phys.}, 353(2):633--690, 2017.
	
	\bibitem{MR4477342}
	N.~Anantharaman.
	\newblock {\em Quantum ergodicity and delocalization of {S}chr\"{o}dinger
		eigenfunctions}.
	\newblock Zurich Lectures in Advanced Mathematics. European Mathematical
	Society (EMS), Z\"{u}rich, [2022] \textcopyright 2022.
	
	\bibitem{MR3322309}
	N.~Anantharaman and E.~Le~Masson.
	\newblock Quantum ergodicity on large regular graphs.
	\newblock {\em Duke Math. J.}, 164(4):723--765, 2015.
	
	\bibitem{MR3961083}
	N.~Anantharaman and M.~Sabri.
	\newblock Quantum ergodicity on graphs: from spectral to spatial
	delocalization.
	\newblock {\em Ann. of Math. (2)}, 189(3):753--835, 2019.
	
	\bibitem{MR236951}
	M.~F. Atiyah and G.~B. Segal.
	\newblock The index of elliptic operators. {II}.
	\newblock {\em Ann. of Math. (2)}, 87:531--545, 1968.
	
	\bibitem{ovadia2025mixedquantizationpartialhyperbolicity}
	S.~{Ben Ovadia}, Q.~Ma, and F.~Rodriguez-Hertz.
	\newblock Mixed quantization and partial hyperbolicity.
	\newblock {\em arXiv: 2409.13660}, 2025.
	
	\bibitem{MR1873300}
	I.~Benjamini and O.~Schramm.
	\newblock Recurrence of distributional limits of finite planar graphs.
	\newblock {\em Electron. J. Probab.}, 6:no. 23, 13, 2001.
	
	\bibitem{MR3028790}
	N.~Bergeron and A.~Venkatesh.
	\newblock The asymptotic growth of torsion homology for arithmetic groups.
	\newblock {\em J. Inst. Math. Jussieu}, 12(2):391--447, 2013.
	
	\bibitem{MR2273508}
	N.~Berline, E.~Getzler, and M.~Vergne.
	\newblock {\em Heat kernels and {D}irac operators}.
	\newblock Grundlehren Text Editions. Springer-Verlag, Berlin, 2004.
	\newblock Corrected reprint of the 1992 original.
	
	\bibitem{MR1188532}
	J.-M. Bismut and G.~Lebeau.
	\newblock Complex immersions and {Q}uillen metrics.
	\newblock {\em Inst. Hautes \'{E}tudes Sci. Publ. Math.}, (74):ii+298 pp.
	(1992), 1991.
	
	\bibitem{MR2838248}
	J.-M. Bismut, X.~Ma, and W.~Zhang.
	\newblock Op\'{e}rateurs de {T}oeplitz et torsion analytique asymptotique.
	\newblock {\em C. R. Math. Acad. Sci. Paris}, 349(17-18):977--981, 2011.
	
	\bibitem{MR3615411}
	J.-M. Bismut, X.~Ma, and W.~Zhang.
	\newblock Asymptotic torsion and {T}oeplitz operators.
	\newblock {\em J. Inst. Math. Jussieu}, 16(2):223--349, 2017.
	
	\bibitem{MR1016875}
	J.-M. Bismut and E.~Vasserot.
	\newblock The asymptotics of the {R}ay-{S}inger analytic torsion associated
	with high powers of a positive line bundle.
	\newblock {\em Comm. Math. Phys.}, 125(2):355--367, 1989.
	
	\bibitem{MR1096593}
	J.-M. Bismut and E.~Vasserot.
	\newblock The asymptotics of the {R}ay-{S}inger analytic torsion of the
	symmetric powers of a positive vector bundle.
	\newblock {\em Ann. Inst. Fourier (Grenoble)}, 40(4):835--848 (1991), 1990.
	
	\bibitem{MR89473}
	R.~Bott.
	\newblock Homogeneous vector bundles.
	\newblock {\em Ann. of Math. (2)}, 66:203--248, 1957.
	
	\bibitem{MR2358056}
	J.~Bourgain and A.~Gamburd.
	\newblock On the spectral gap for finitely-generated subgroups of {$\rm
		SU(2)$}.
	\newblock {\em Invent. Math.}, 171(1):83--121, 2008.
	
	\bibitem{MR2966656}
	J.~Bourgain and A.~Gamburd.
	\newblock A spectral gap theorem in {${\rm SU}(d)$}.
	\newblock {\em J. Eur. Math. Soc. (JEMS)}, 14(5):1455--1511, 2012.
	
	\bibitem{MR0781344}
	T.~Br\"{o}cker and T.~tom Dieck.
	\newblock {\em Representations of compact {L}ie groups}, volume~98 of {\em
		Graduate Texts in Mathematics}.
	\newblock Springer-Verlag, New York, 1985.
	
	\bibitem{MR3567266}
	S.~Brooks, E.~Le~Masson, and E.~Lindenstrauss.
	\newblock Quantum ergodicity and averaging operators on the sphere.
	\newblock {\em Int. Math. Res. Not. IMRN}, (19):6034--6064, 2016.
	
	\bibitem{MR3961523}
	F.~Calegari and A.~Venkatesh.
	\newblock A torsion {J}acquet-{L}anglands correspondence.
	\newblock {\em Ast\'erisque}, (409):x+226, 2019.
	
	\bibitem{cekić2024semiclassicalanalysisprincipalbundles}
	M.~Ceki\'{c} and T.~Lefeuvre.
	\newblock Semiclassical analysis on principal bundles.
	\newblock {\em arXiv: 2405.14846}, 2024.
	
	\bibitem{MR818831}
	Y.~Colin~de Verdi{\`e}re.
	\newblock Ergodicit\'{e} et fonctions propres du laplacien.
	\newblock {\em Comm. Math. Phys.}, 102(3):497--502, 1985.
	
	\bibitem{MR1438595}
	J.~B. Conrey, W.~Duke, and D.~W. Farmer.
	\newblock The distribution of the eigenvalues of {H}ecke operators.
	\newblock {\em Acta Arith.}, 78(4):405--409, 1997.
	
	\bibitem{MR4922234}
	A.~Drewitz, B.~Liu, and G.~Marinescu.
	\newblock Gaussian holomorphic sections on noncompact complex manifolds.
	\newblock {\em J. Inst. Math. Jussieu}, 24(4):1197--1262, 2025.
	
	\bibitem{MR4905035}
	A.~Drewitz, B.~Liu, and G.~Marinescu.
	\newblock Large deviations for zeros of holomorphic sections on punctured
	{R}iemann surfaces.
	\newblock {\em Michigan Math. J.}, 75(2):381--421, 2025.
	
	\bibitem{MR757256}
	V.~G. Drinfel\cprime~d.
	\newblock Finitely-additive measures on {$S\sp{2}$}\ and {$S\sp{3}$}, invariant
	with respect to rotations.
	\newblock {\em Funktsional. Anal. i Prilozhen.}, 18(3):77, 1984.
	
	\bibitem{MR3969938}
	S.~Dyatlov and M.~Zworski.
	\newblock {\em Mathematical theory of scattering resonances}, volume 200 of
	{\em Graduate Studies in Mathematics}.
	\newblock American Mathematical Society, Providence, RI, 2019.
	
	\bibitem{Arizona}
	M.~Einsiedler and T.~Wald.
	\newblock Quantum unique ergodicity on {$\Gamma\backslash\mathbb{H}$}.
	\newblock {\em Arizona Winter School, Tucson}, 2010.
	
	\bibitem{etingof2024representationsliegroups}
	P.~Etingof.
	\newblock Representations of {L}ie groups.
	\newblock {\em arXiv: 2401.01446}, 2024.
	
	\bibitem{MR3864507}
	S.~Finski.
	\newblock On the full asymptotics of analytic torsion.
	\newblock {\em J. Funct. Anal.}, 275(12):3457--3503, 2018.
	
	\bibitem{MR1137767}
	J.~Friedman.
	\newblock On the second eigenvalue and random walks in random {$d$}-regular
	graphs.
	\newblock {\em Combinatorica}, 11(4):331--362, 1991.
	
	\bibitem{MR2437174}
	J.~Friedman.
	\newblock A proof of {A}lon's second eigenvalue conjecture and related
	problems.
	\newblock {\em Mem. Amer. Math. Soc.}, 195(910):viii+100, 2008.
	
	\bibitem{MR1677685}
	A.~Gamburd, D.~Jakobson, and P.~Sarnak.
	\newblock Spectra of elements in the group ring of {${\rm SU}(2)$}.
	\newblock {\em J. Eur. Math. Soc. (JEMS)}, 1(1):51--85, 1999.
	
	\bibitem{MR3725881}
	C.~Hall, D.~Puder, and W.~F. Sawin.
	\newblock Ramanujan coverings of graphs.
	\newblock {\em Adv. Math.}, 323:367--410, 2018.
	
	\bibitem{MR1867354}
	A.~Hatcher.
	\newblock {\em Algebraic topology}.
	\newblock Cambridge University Press, Cambridge, 2002.
	
	\bibitem{MR2247919}
	S.~Hoory, N.~Linial, and A.~Wigderson.
	\newblock Expander graphs and their applications.
	\newblock {\em Bull. Amer. Math. Soc. (N.S.)}, 43(4):439--561, 2006.
	
	\bibitem{huang2025ramanujanpropertyedgeuniversality}
	J.~Huang, T.~McKenzie, and H.-T. Yau.
	\newblock Ramanujan property and edge universality of random regular graphs.
	\newblock {\em arXiv: 2412.20263}, 2025.
	
	\bibitem{MR2093043}
	D.~Huybrechts.
	\newblock {\em Complex geometry}.
	\newblock Universitext. Springer-Verlag, Berlin, 2005.
	\newblock An introduction.
	
	\bibitem{MR2884879}
	R.~Kenyon.
	\newblock Spanning forests and the vector bundle {L}aplacian.
	\newblock {\em Ann. Probab.}, 39(5):1983--2017, 2011.
	
	\bibitem{MR109367}
	H.~Kesten.
	\newblock Symmetric random walks on groups.
	\newblock {\em Trans. Amer. Math. Soc.}, 92:336--354, 1959.
	
	\bibitem{MR909698}
	S.~Kobayashi.
	\newblock {\em Differential geometry of complex vector bundles}, volume~15 of
	{\em Publications of the Mathematical Society of Japan}.
	\newblock Princeton University Press, Princeton, NJ; Princeton University
	Press, Princeton, NJ, 1987.
	\newblock Kan\^o{} Memorial Lectures, 5.
	
	\bibitem{MR3245884}
	E.~Le~Masson.
	\newblock Pseudo-differential calculus on homogeneous trees.
	\newblock {\em Ann. Henri Poincar\'e}, 15(9):1697--1732, 2014.
	
	\bibitem{MR4055707}
	W.-C.~W. Li.
	\newblock The {R}amanujan conjecture and its applications.
	\newblock {\em Philos. Trans. Roy. Soc. A}, 378(2163):20180441, 14, 2020.
	
	\bibitem{MR2195133}
	E.~Lindenstrauss.
	\newblock Invariant measures and arithmetic quantum unique ergodicity.
	\newblock {\em Ann. of Math. (2)}, 163(1):165--219, 2006.
	
	\bibitem{MR4746871}
	B.~Liu.
	\newblock On full asymptotics of real analytic torsions for compact locally
	symmetric orbifolds.
	\newblock {\em Anal. PDE}, 17(4):1261--1329, 2024.
	
	\bibitem{MR3558308}
	E.~Lubetzky and Y.~Peres.
	\newblock Cutoff on all {R}amanujan graphs.
	\newblock {\em Geom. Funct. Anal.}, 26(4):1190--1216, 2016.
	
	\bibitem{MR861487}
	A.~Lubotzky, R.~Phillips, and P.~Sarnak.
	\newblock Hecke operators and distributing points on the sphere. {I}.
	\newblock volume~39, pages S149--S186. 1986.
	\newblock Frontiers of the mathematical sciences: 1985 (New York, 1985).
	
	\bibitem{MR890171}
	A.~Lubotzky, R.~Phillips, and P.~Sarnak.
	\newblock Hecke operators and distributing points on {$S^2$}. {II}.
	\newblock {\em Comm. Pure Appl. Math.}, 40(4):401--420, 1987.
	
	\bibitem{MR963118}
	A.~Lubotzky, R.~Phillips, and P.~Sarnak.
	\newblock Ramanujan graphs.
	\newblock {\em Combinatorica}, 8(3):261--277, 1988.
	
	\bibitem{MR1990480}
	W.~Luo and P.~Sarnak.
	\newblock Mass equidistribution for {H}ecke eigenforms.
	\newblock volume~56, pages 874--891. 2003.
	\newblock Dedicated to the memory of J\"urgen K. Moser.
	
	\bibitem{MR2160416}
	R.~Lyons.
	\newblock Asymptotic enumeration of spanning trees.
	\newblock {\em Combin. Probab. Comput.}, 14(4):491--522, 2005.
	
	\bibitem{MR4808253}
	M.~Ma and Q.~Ma.
	\newblock Semiclassical analysis, geometric representation and quantum
	ergodicity.
	\newblock {\em Comm. Math. Phys.}, 405(11):Paper No. 259, 28, 2024.
	
	\bibitem{MR4665497}
	Q.~Ma.
	\newblock Toeplitz operators and the full asymptotic torsion forms.
	\newblock {\em J. Funct. Anal.}, 286(3):Paper No. 110210, 74pp, 2024.
	
	\bibitem{MR2339952}
	X.~Ma and G.~Marinescu.
	\newblock {\em Holomorphic {M}orse inequalities and {B}ergman kernels}, volume
	254 of {\em Progress in Mathematics}.
	\newblock Birkh\"{a}user Verlag, Basel, 2007.
	
	\bibitem{MR2393271}
	X.~Ma and G.~Marinescu.
	\newblock Toeplitz operators on symplectic manifolds.
	\newblock {\em J. Geom. Anal.}, 18(2):565--611, 2008.
	
	\bibitem{MR3368102}
	X.~Ma and G.~Marinescu.
	\newblock Exponential estimate for the asymptotics of {B}ergman kernels.
	\newblock {\em Math. Ann.}, 362(3-4):1327--1347, 2015.
	
	\bibitem{MR4630491}
	M.~Magee, J.~Thomas, and Y.~Zhao.
	\newblock Quantum unique ergodicity for {C}ayley graphs of quasirandom groups.
	\newblock {\em Comm. Math. Phys.}, 402(3):3021--3044, 2023.
	
	\bibitem{MR3374962}
	A.~W. Marcus, D.~A. Spielman, and N.~Srivastava.
	\newblock Interlacing families {I}: {B}ipartite {R}amanujan graphs of all
	degrees.
	\newblock {\em Ann. of Math. (2)}, 182(1):307--325, 2015.
	
	\bibitem{MR484767}
	G.~A. Margulis.
	\newblock Explicit constructions of expanders.
	\newblock {\em Problemy Pereda\v ci Informacii}, 9(4):71--80, 1973.
	
	\bibitem{MR596890}
	G.~A. Margulis.
	\newblock Some remarks on invariant means.
	\newblock {\em Monatsh. Math.}, 90(3):233--235, 1980.
	
	\bibitem{MR939574}
	G.~A. Margulis.
	\newblock Explicit group-theoretic constructions of combinatorial schemes and
	their applications in the construction of expanders and concentrators.
	\newblock {\em Problemy Peredachi Informatsii}, 24(1):51--60, 1988.
	
	\bibitem{MR629617}
	B.~D. McKay.
	\newblock The expected eigenvalue distribution of a large regular graph.
	\newblock {\em Linear Algebra Appl.}, 40:203--216, 1981.
	
	\bibitem{MR3307755}
	D.~W. Morris.
	\newblock {\em Introduction to arithmetic groups}.
	\newblock Deductive Press, [place of publication not identified], 2015.
	
	\bibitem{MR3220447}
	W.~M\"{u}ller.
	\newblock The asymptotics of the {R}ay-{S}inger analytic torsion of hyperbolic
	3-manifolds.
	\newblock In {\em Metric and differential geometry}, volume 297 of {\em Progr.
		Math.}, pages 317--352. Birkh\"{a}user/Springer, Basel, 2012.
	
	\bibitem{MR4651019}
	A.~Naor, A.~Sah, M.~Sawhney, and Y.~Zhao.
	\newblock Cayley graphs that have a quantum ergodic eigenbasis.
	\newblock {\em Israel J. Math.}, 256(2):599--617, 2023.
	
	\bibitem{MR1124768}
	A.~Nilli.
	\newblock On the second eigenvalue of a graph.
	\newblock {\em Discrete Math.}, 91(2):207--210, 1991.
	
	\bibitem{MR1649013}
	S.~Nonnenmacher and A.~Voros.
	\newblock Chaotic eigenfunctions in phase space.
	\newblock {\em J. Statist. Phys.}, 92(3-4):431--518, 1998.
	
	\bibitem{MR858831}
	A.~Perelomov.
	\newblock {\em Generalized coherent states and their applications}.
	\newblock Texts and Monographs in Physics. Springer-Verlag, Berlin, 1986.
	
	\bibitem{MR363209}
	A.~M. Perelomov.
	\newblock Coherent states for arbitrary {L}ie group.
	\newblock {\em Comm. Math. Phys.}, 26:222--236, 1972.
	
	\bibitem{MR4611826}
	M.~Puchol.
	\newblock The asymptotics of the holomorphic analytic torsion forms.
	\newblock {\em J. Lond. Math. Soc. (2)}, 108(1):80--140, 2023.
	
	\bibitem{MR1266075}
	Z.~Rudnick and P.~Sarnak.
	\newblock The behaviour of eigenstates of arithmetic hyperbolic manifolds.
	\newblock {\em Comm. Math. Phys.}, 161(1):195--213, 1994.
	
	\bibitem{MR1018385}
	P.~Sarnak.
	\newblock Statistical properties of eigenvalues of the {H}ecke operators.
	\newblock In {\em Analytic number theory and {D}iophantine problems
		({S}tillwater, {OK}, 1984)}, volume~70 of {\em Progr. Math.}, pages 321--331.
	Birkh\"auser Boston, Boston, MA, 1987.
	
	\bibitem{MR0995750}
	R.~Schrader and M.~E. Taylor.
	\newblock Semiclassical asymptotics, gauge fields, and quantum chaos.
	\newblock {\em J. Funct. Anal.}, 83(2):258--316, 1989.
	
	\bibitem{MR344216}
	J.-P. Serre.
	\newblock {\em A course in arithmetic}, volume No. 7 of {\em Graduate Texts in
		Mathematics}.
	\newblock Springer-Verlag, New York-Heidelberg, 1973.
	\newblock Translated from the French.
	
	\bibitem{MR1396897}
	J.-P. Serre.
	\newblock R\'epartition asymptotique des valeurs propres de l'op\'erateur de
	{H}ecke {$T_p$}.
	\newblock {\em J. Amer. Math. Soc.}, 10(1):75--102, 1997.
	
	\bibitem{MR1675133}
	B.~Shiffman and S.~Zelditch.
	\newblock Distribution of zeros of random and quantum chaotic sections of
	positive line bundles.
	\newblock {\em Comm. Math. Phys.}, 200(3):661--683, 1999.
	
	\bibitem{MR2680500}
	K.~Soundararajan.
	\newblock Quantum unique ergodicity for {${\rm SL}_2(\Bbb Z)\backslash\Bbb H$}.
	\newblock {\em Ann. of Math. (2)}, 172(2):1529--1538, 2010.
	
	\bibitem{MR590825}
	D.~Sullivan.
	\newblock For {$n>3$}\ there is only one finitely additive rotationally
	invariant measure on the {$n$}-sphere defined on all {L}ebesgue measurable
	subsets.
	\newblock {\em Bull. Amer. Math. Soc. (N.S.)}, 4(1):121--123, 1981.
	
	\bibitem{MR0402834}
	A.~I. \v{S}nirel$'$man.
	\newblock Ergodic properties of eigenfunctions.
	\newblock {\em Uspehi Mat. Nauk}, 29(6(180)):181--182, 1974.
	
	\bibitem{MR916129}
	S.~Zelditch.
	\newblock Uniform distribution of eigenfunctions on compact hyperbolic
	surfaces.
	\newblock {\em Duke Math. J.}, 55(4):919--941, 1987.
	
	\bibitem{MR2952218}
	M.~Zworski.
	\newblock {\em Semiclassical analysis}, volume 138 of {\em Graduate Studies in
		Mathematics}.
	\newblock American Mathematical Society, Providence, RI, 2012.
	
\end{thebibliography}

\end{document}